\documentstyle[12pt]{article}
\input pstricks

\def\marginpar#1{}
\input epsf
\newcommand\vol{\mathrm {vol}}

\newcommand\card{\mathrm {card}}
\newcommand\Var{\mathrm {Var}}
\newcommand\sg{\mathrm {sgn}}
\newtheorem{theo}{Theorem}

\newtheorem{defi}{Definition}
\newtheorem {lemma}{Lemma}
\newtheorem {coro}{Corollary}
\newtheorem {pro}{Proposition}
\newdimen\AAdi%
\newbox\AAbo%
\def\AArm{\fam0 }
\def\AAk#1#2{\setbox\AAbo=\hbox{#2}\AAdi=\wd\AAbo\kern#1\AAdi{}}%
\def\AAr#1#2#3{\setbox\AAbo=\hbox{#2}\AAdi=\ht\AAbo\raise#1\AAdi\hbox{#3}}%
\def\BBe{{\AArm I\!E}}%
\def\BBn{{\AArm I\!N}}%
\def\BBr{{\AArm I\!R}}%
\def\BBz{{\AArm Z\!\!Z}}%
\def\BBone{{\AArm 1\AAk{-.8}{I}I}}%

\begin{document}
\title{The initial drift of a 2D droplet at zero temperature}
\date{}
\author{Rapha\"el
Cerf \footnote{Universit\'e de Paris-Sud, Probabilit\'es,
statistique et mod\'elisation, B\^at. 425, 91405 Orsay Cedex,
France. E-mail: rcerf@math.u-psud.fr}\,\,\, and Sana Louhichi
\footnote{Corresponding author. Universit\'e de Paris-Sud,
Probabilit\'es, statistique et mod\'elisation, B\^at. 425, 91405
Orsay Cedex, France. E-mail: sana.louhichi@math.u-psud.fr
\newline
 {{\bf{Acknowledgments.}} The
authors wish to thank Patrice Assouad and Sophie Lemaire for
useful discussions.} The first version of this work dealt only
with the deterministic initial condition; we are very grateful to
Herbert Spohn for explaining to us the relevance of the randomness
in the initial condition.} } \maketitle
\vspace{4cm}
\begin{abstract}
We consider the 2D stochastic Ising model evolving according to
the Glauber dynamics at zero temperature. We compute the initial
drift for droplets which are suitable approximations
of smooth domains. A
specific spatial average of the derivative at time~$0$ of the
volume variation of a droplet close to a boundary point is equal
to its curvature multiplied by a direction dependent coefficient.
We compute the explicit value of this coefficient.
\end{abstract}

\newpage

\centerline{
\epsfxsize 14.25cm
\epsfbox{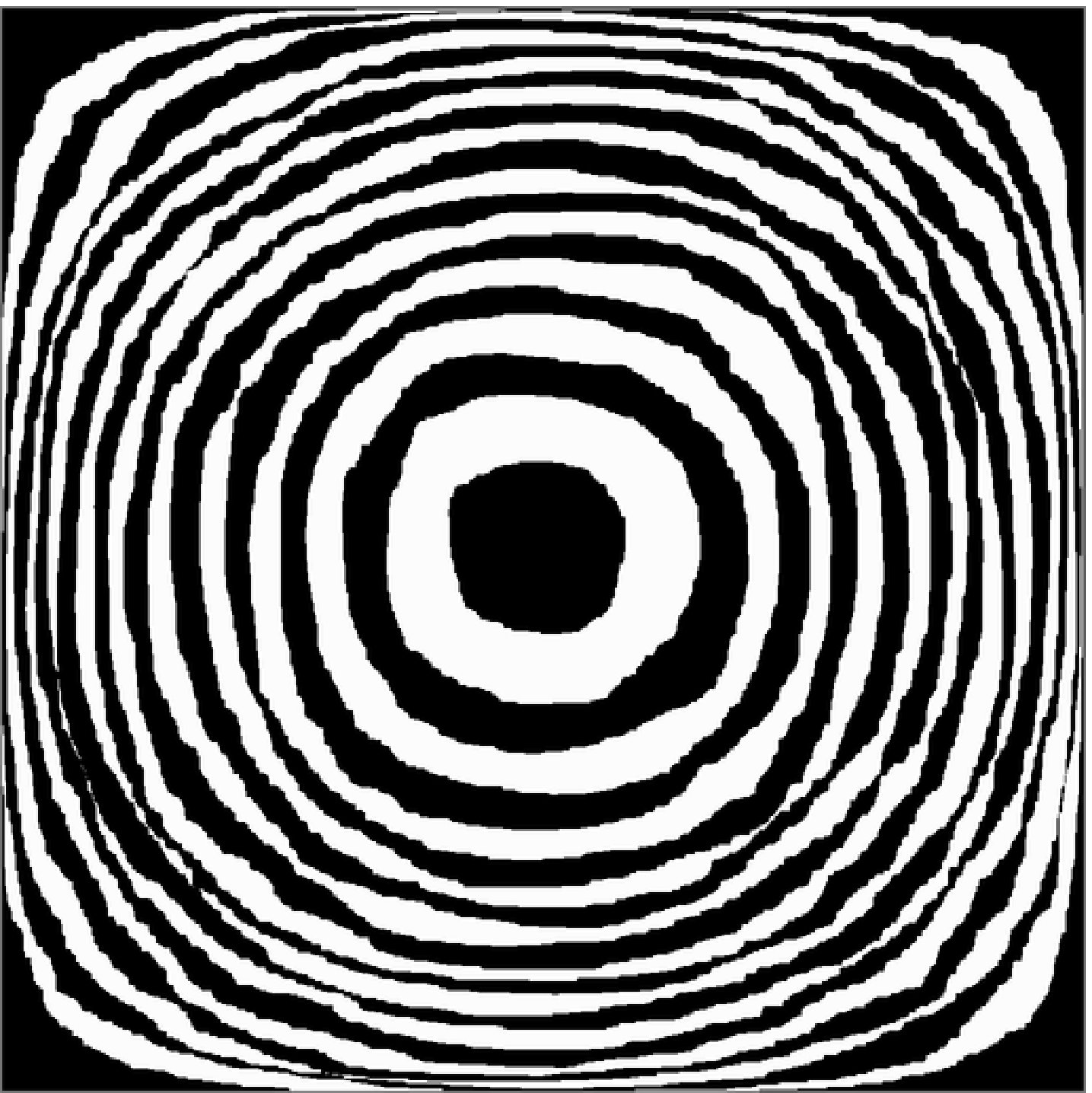}
}
\vspace{2cm}
\centerline{Evolution of a square droplet}
\newpage
\section*{Introduction}
The phenomenological theory asserts that the evolution of the
shape of a droplet of one phase immersed in another phase is
governed by the motion by mean curvature. We are still far from
being able to verify this assertion starting from a genuine
microscopic dynamics. Very interesting results have been obtained
in a series of works in the context of the Ising model with Ka\'c
potentials \cite{DOPT1, DOPT2, KS, KSP}. However, motion  by mean
curvature is recovered in some scaling limit where the range of
the interactions diverges to infinity: the model becomes somehow
close to a mean--field model and the ensuing motion is isotropic.
For the true Ising model with only nearest--neighbour
interactions, it is expected that an interface between the minus
and the plus phase evolves according to an anisotropic motion  by
mean curvature, that is, each point~$x$ of the interface has
velocity
$$v(x)\,=\,-c(\nu(x))\,\xi\,\nu(x)$$
where $\nu(x)$ is the vector normal to the interface at~$x$, $\xi$
is the curvature of the interface at~$x$ and $c(\nu)$ is a
coefficient depending on the direction of~$\nu$. This anisotropy
stems from the anisotropy of the cubic lattice.

In this paper, we consider the zero temperature Glauber dynamics
for the 2D Ising model. Although we do not succeed in deriving the
full motion by mean curvature, we manage to compute the initial
drift for droplets which approximate suitably
smooth domains and
we believe this is a crucial step. Four works are directly
relevant. In \cite{Sp}, Spohn claims to establish rigorously the
mean curvature motion in the context of the 2D Ising model at zero
temperature for interfaces which can be represented as the graph
of a function. 
Although his results do not apply directly to the case of a
full droplet, he succeeds in deriving an explicit formula for
the coefficient $c(\nu)$. We recover this result here with a
different approach.
The
computation we present here can be considered to be a refinement
of the observation of \cite{CSS}. Chayes, Schonmann and Swindle
proved a Lifshitz law for the volume of a two-dimensional droplet
at zero temperature. Instead of looking at the total volume of the
droplet, we shall concentrate here on the volume variation of the
droplet in a small ball attached to its boundary. In \cite{CS}, by
interpreting the interface as a one dimensional exclusion process,
Chayes and Swindle manage to prove that, starting from a square
droplet, the evolution of the shape of one corner is described in
the hydrodynamical limit by an appropriate Stefan problem.
Finally, Sowers develops in \cite{So} a framework of geometric
measure theory to obtain the hydrodynamical limit. His convergence
theorem is conditional on the verification of several assumptions,
some of them concerning the structure of the interface. It might
be that these estimates are the missing pieces to complete the
picture.

Let us turn now to the description of our result. We work with the
stochastic Ising model evolving according to the Glauber dynamics
at zero temperature. We consider the diffusive limit where space
is rescaled by a factor~$N$ and time is speeded up by a
factor~$N^2$. We start with a plus droplet immersed in
the minus phase, whose boundary is a
${\cal C}_1$ simple Jordan curve~$\gamma$: the initial
configuration at step~$N$ is a suitable approximation
of the smooth
droplet, drawn on the square lattice ${\BBz}^2/N$.
We consider two cases:
\medskip

\noindent
{\bf{Deterministic initial condition.}}
The approximating set
at step~$N$
consists of the squares of the lattice ${\BBz}^2/N$
which intersect the interior of~$\gamma$.
\medskip

\noindent
{\bf{Spohn's initial condition.}}
The approximating set
at step~$N$ is random.
Its boundary converges in probability towards~$\gamma$ as
$N$ goes to~$\infty$ and its law~$\mu_N$ is given by
the invariant measure of the associated zero range process.
\medskip

\noindent
The droplet is immersed
in the minus phase, hence all the sites of the approximating set
are initially set to plus, while the other sites of the lattice
are set to minus.
We then look at the process
$(\sigma_{N^2t},t\geq 0)$ and we denote by ${\cal A}^N_\sigma(t)$ the
plus droplet at time~$N^2t$ starting from~$\sigma$.
Let $x$ be a point of~$\gamma$. We
study the variation of the magnetization inside the ball $B(x,r)$
centered at $x$ with radius $r$, for $r$ small. Equivalently, we
look at the volume $\vol\big(B(x,r)\cap {\cal A}^N_\sigma(t)\big)$ of the
plus droplet in this ball and we aim at computing its derivative
$$\lim_{t\to 0}
\,\,{1\over t} \Big(\vol\big(B(x,r)\cap {\cal A}^N_\sigma(t)\big)-
\vol\big(B(x,r)\cap {\cal A}^N_\sigma(0)\big)\Big)\,.$$
Several problems
arise. Since the dynamics proceeds by jumps, we have to take the
expectation to get a differentiable quantity. Next we wish to link
the infinitesimal volume variation with the curvature of the
droplet's boundary at~$x$. To achieve this, we need to recover
the slope of the continuous curve from its
approximation.
We perform a spatial averaging. Letting
$x_0,x_1$ be the two points of $\gamma$ which belong to the sphere
$\partial B(x,r)$, we consider the domain
 $$
{\cal{S}}(x,r,\alpha_1,\alpha_2)={B}(x,r)\cup { B}(x_{0},\alpha_1)
\cup {B}(x_1,\alpha_2),
$$
and we denote by ${\cal{S}}_N$ its discretization at step~$N$. The
quantity of primary interest to link the volume variation and the
curvature is
$${\bf A}_N^{\sigma,\gamma}(x,r,\delta)=\frac{1}{\delta^2}\int_{0}^{\delta}\int_{0}^{\delta}
 \lim_{t\rightarrow 0}\,\,
\frac{1}{t}\,\BBe\left(({\vol}({\cal{A}}^{N}_\sigma(t)\cap
{\cal{S}}_N))-{\vol}({\cal{A}}^{N}_\sigma(0)\cap
{\cal{S}}_N)\right) \,d\alpha_1\,d\alpha_2.$$ Let $\theta$ be the
angle of the tangent to $\gamma$ at~$x$ and let $\xi_{\gamma}(x)$
be the curvature of $\gamma$ at~$x$. Our main result states that,
for the deterministic initial condition,
$$
 \lim_{r\rightarrow 0}\,
 \lim_{\delta\rightarrow 0}\,
 \liminf_{N\rightarrow \infty}\,
\frac{1}{2r}{\bf A}_N^{\sigma,\gamma}(x,r,\delta)=\lim_{r\rightarrow 0}\,
 \lim_{\delta\rightarrow 0}\,
 \limsup_{N\rightarrow \infty}\,
\frac{1}{2r}{\bf A}_N^{\sigma,\gamma}(x,r,\delta)=-\frac{1 }{2}
|\cos(2\theta)|\xi_{\gamma}(x)$$
while for
Spohn's initial condition,
$$\displaylines{
 \lim_{r\rightarrow 0}\,
 \lim_{\delta\rightarrow 0}\,
 \liminf_{N\rightarrow \infty}\,
\frac{1}{2r}\mu_N({\bf A}_N^{\sigma,\gamma}(x,r,\delta))
\,=\,\lim_{r\rightarrow 0}\,
 \lim_{\delta\rightarrow 0}\,
 \limsup_{N\rightarrow \infty}\,
\frac{1}{2r}\mu_N({\bf A}_N^{\sigma,\gamma}(x,r,\delta))
\hfill\cr
= -\frac{ \xi_{\gamma}(x) }{2(|\cos\theta| + |\sin\theta|)^2}\,.}
$$
In fact, we compute the above limits for a more general class of
initial conditions, which includes the two cases above.
The physically relevant case should be the one studied by Spohn,
it corresponds to the equilibrium state of the zero range process.
This
indicates that the limit $({\cal{A}}(t),t\geq 0)$ of any decently
converging subsequence of the stochastic motion
$({\cal{A}}^{N}(t),t\geq 0)$ should satisfy the equation, for any
$s>0$ and for any $x\in\partial {\cal A}(s)$,
 $$\displaylines{\lim_{r\rightarrow 0}\,
 \lim_{\delta\rightarrow 0}\,
\frac{1}{2r\delta^2}\int_{0}^{\delta}\int_{0}^{\delta}
 \lim_{\scriptstyle t\rightarrow s\atop \scriptstyle t>s }\,
\frac{1}{t-s}\,\BBe\left({\vol}({\cal{A}}(t)\cap
{\cal{S}})-{\vol}({\cal{A}}(s)\cap {\cal{S}})\right)
\,d\alpha_1\,d\alpha_2 \hfill\cr = -\frac{
\xi_{\partial{\cal A}(s)}(x)
}{2(|\cos\theta| + |\sin\theta|)^2}\,}$$
 or at least a weaker
variant of it. Here $({\cal{A}}(t),t\geq 0)$ is a random process
describing the evolution of the shape of the droplet. A standard
computation shows that the deterministic motion by mean curvature
satisfies this equation. However we do not know whether it is the
only solution to this equation; we have not investigated the
corresponding theory so far. For instance, can one get rid of the
expectation? Anyway, we are still far from establishing that the
hydrodynamical limit of the droplet process satisfies the above
equation.
An important issue
is to control
dynamically the proportion of the corners in a microscopic random
interface when its average slope is known. This would probably
require some additional probabilistic input.

\section{The model}
We consider a zero-temperature 2D-stochastic Ising model. More
precisely it is a continuous time Markov process
$(\sigma_t)_{t\geq 0}$ taking values in $\{-1,+1\}^{\BBz^2}$ with
generator $L$ which acts on each local function
$f:\{-1,+1\}^{\BBz^2}\rightarrow \BBr$ as
$$(Lf)(\sigma)=\sum_{x\in
\BBz^2}c(x,\sigma)(f(\sigma^{x})-f(\sigma)). $$
Here,
for
$\smash{\sigma\in \{-1,+1\}^{\BBz^2}}$ and $x\in\BBz^2$, we define
 $$
\forall y\in{\BBz^2}\qquad
 \sigma^x(y) = \left\{
\begin{array}{rl}
\sigma(y)   & {\mbox{if $y\neq x$,}} \\ -\sigma(y) & {\mbox{if $y=
x$,}}
\\
\end{array}
  \right.
$$ and $c(x,\sigma)$ is the rate with which the spin at site $x$
flips when the configuration is $\sigma$. The rates $c(x,\sigma)$
define the dynamics. For the zero-temperature 2D--Ising model, the
rates $c(x,\sigma)$ are given by
 $$
 c(x,\sigma) = \left\{
\begin{array}{rl}

  1      & {\mbox{if more than 2 neighbors of $x$ have a spin opposite to $x$,}} \\
  \alpha & {\mbox{if exactly 2 neighbors of $x$ have a spin opposite to $x$,}} \\
  0      & {\mbox{otherwise,}}
\end{array}
  \right.
$$
where $0<\alpha\leq 1$ is a fixed parameter. For technical
reasons, we will take $\alpha=\smash{\frac{\textstyle 1}{\textstyle 2}}$
in the sequel.
\section{Notation}
 Let $N$ be a fixed positive integer. We denote by $\BBz_N^2$ the
 grid $\frac{\textstyle \BBz^2}{\textstyle N}$.
For $x=(x_1,x_2)\in\BBz^2$, $\Lambda_{x /  N}$ is the box defined
as

 \begin{figure}[hbt]
\vbox{ \centerline{ \psset{unit=2.5cm}
\pspicture(-0.5,-0.5)(1.5,1.5)
\psline[linewidth=2pt](0,1)(0,0)(1,0)
\psline[linestyle=dashed,linewidth=2pt](1,0)(1,1)(0,1)
\psline[linestyle=dashed,linewidth=0.5pt](-0.2,0.5)(1.5,0.5)
\psline[linestyle=dashed,linewidth=0.5pt](0.5,-0.2)(0.5,1.5)
\psline[linewidth=2pt]{<->}(-0.2,0)(-0.2,1)
\rput(-0.3,0.5){$1\over N$}
\psline[linewidth=2pt]{<->}(0,-0.2)(1,-0.2)
\rput(0.5,-0.3){$1\over N$} \psdots(0.5,0.5)
\rput(0.35,0.35){$u_N$} \psdots[dotstyle=x,dotscale=1.2](0.8,0.7)
\rput(0.7,0.7){$u$}
\endpspicture
} \centerline{A point $u$ and the box $\Lambda_{u_N}$ ($Nu_N\in
\BBz^2$).} }
\end{figure}

\begin{equation}\label{boite}
\Lambda_{{x}/ {N}}=\left\{(u_1,u_2)\in\BBr^2,\ \ \
-\frac{1}{2N}\leq u_1-\frac{x_1}{N}<\frac{1}{2N}; \ \
  -\frac{1}{2N}\leq u_2-\frac{x_2}{N}< \frac{1}{2N} \right\}.
\end{equation}
 The family of boxes $(\Lambda_{x},\,
x\in {\BBz_N^2})$, as defined by (\ref{boite}), forms a partition
of $\BBr^2$:
$$\BBr^2={\bigcup}_{x\in
{\BBz_N^2}}\Lambda_{x},\qquad \forall\, x,y\in\,\BBz_N^2
\quad x\neq y\Rightarrow \Lambda_{x}\cap \Lambda_{y}=\emptyset.$$ Hence, for each
$u=(u_1,u_2)\in \BBr^2$ there exists a
  unique $u_N \in \BBz_N^2$ such that
$u\in \Lambda_{u_N}$. Moreover
  $
  \|u-u_N\|_{\infty}\leq \frac{1}{2N},
  $
  where $\|u\|_{\infty}=\max(|u_1|,|u_2|)$.
\\
To each bounded set ${\cal{S}}$ of $\BBr^2$, we associate
 the set ${\cal{S}}_N$ defined by $${\cal{S}}_N= {\bigcup}_{x\in {\BBz_N^2}:\,\Lambda_x\cap
{\cal{S}}\neq \emptyset} \Lambda_{x}.$$

\begin{figure}[hbt]
\vbox{ \centerline{ \psset{unit=0.5cm} \pspicture(0,-1)(10,9)
\pscurve[fillstyle=solid,fillcolor=gray,linewidth=1.75pt](1.2,2.3)(2.3,2.2)(4.5,3.4)(5.7,5.4)(5.2,7.1)(1.2,2.3)
\pspolygon[fillcolor=lightgray,
linewidth=2.5pt](0.75,1.5)(1.5,1.5)(2.25,1.5)(3,1.5)(3,2.25)(3.75,2.25)
(4.5,2.25)(4.5,3)(5.25,3)(5.25,3.75)(6,3.75)(6,4.5)(6,5.25)(6,6)(6,6.75)
(6,7.5)(3.75,7.5)(3.75,6.75)(3,6.75)(3,6)(2.25,6)(2.25,5.25)(1.5,5.25)(1.5,4.5)
(1.5,3.75)(0.75,3.75)(0.75,3.5)(0.75,1.5)
\pspolygon[fillstyle=solid,fillcolor=gray,
linewidth=0pt](-2,3)(-1.5,3)(-1.5,3.5)(-2,3.5)
\rput(-0.65,3.25){${\cal S}$}
\psgrid[subgriddiv=0,gridlabels=0,griddots=10,unit=0.75](10,11)
\endpspicture
}
\centerline{The set ${\cal S}$ is included in the set ${\cal S}_N$
with polygonal boundary.} }
\end{figure}

\noindent For $\sigma \in \{-1,+1\}^{\BBz^2}$ and for $x\in
\BBz^2$, we denote by $s({\sigma},x)$, the number of the neighbors
of $x$ having a spin opposite to $x$ in the configuration
$\sigma$:
$$ s({\sigma},x)=\frac{1}{2}\sum_{y\in\BBz^2,\,|x-y|=1
}\kern-20pt|\sigma(x)-\sigma(y)|, $$ where
$|x|=\sqrt{x_1^2+x_2^2}$ for $x=(x_1,x_2)$. \\
Let $N$ be a fixed
positive integer, we define the set
$${\cal{A}}^{\sigma}_{N}=\kern-3pt\bigcup_{x \in \BBz^2,\
\sigma(x)=+1}\kern-3pt\Lambda_{x/ N}.$$

\begin{figure}[hbt]
\vbox{ \centerline{ \psset{unit=0.75cm} \pspicture(0,-1)(10,9)
\pspolygon[fillstyle=solid, fillcolor=lightgray,
linewidth=2.5pt](0.75,6)(1.5,6)(1.5,8.25)(0.75,8.25)
\pspolygon[fillstyle=solid, fillcolor=lightgray, linewidth=2.5pt]
(6,0.75)(6.75,0.75)(6.75,3)(6,3)(6,2.25)(5.25,2.25)(5.25,1.5)(6,1.5)
\pspolygon[fillstyle=solid, fillcolor=lightgray,
linewidth=2.5pt](0.75,1.5)(1.5,1.5)(2.25,1.5)(3,1.5)(3,2.25)(3.75,2.25)
(4.5,2.25)(4.5,3)(5.25,3)(5.25,3.75)(6,3.75)(6,4.5)(6,5.25)(6,6)(6,6.75)
(6,7.5)(3.75,7.5)(3.75,6.75)(3,6.75)(3,6)(2.25,6)(2.25,5.25)(1.5,5.25)(1.5,4.5)
(1.5,3.75)(0.75,3.75)(0.75,3.5)(0.75,1.5)
\pspolygon[fillstyle=solid, fillcolor=lightgray,
linewidth=2.5pt](6.75,6)(7.5,6)(7.5,6.75)(6.75,6.75)
\pspolygon[fillstyle=solid,fillcolor=lightgray,
linewidth=2.5pt](-2,5)(-1.25,5)(-1.25,5.75)(-2,5.75)
\rput(-2.9,5.5){${\cal A}^{\sigma}_N$}
\psgrid[subgriddiv=0,griddots=10,unit=0.75,gridlabels=0](12,12)
\endpspicture
} \vskip -0.59truecm \vspace{5mm} \centerline{ For $x\in
\BBz^2$, $\sigma(x)=+1$ if and only if $x\in N{\cal
A}^{\sigma}_{N}$.} }
\end{figure}
\vspace{1cm} \noindent Let $\gamma$ be a curve of $\BBr^2$. We
define for
 $s\in\gamma$ and for $r,\alpha_1,\alpha_2$  positive real numbers,
 the set
 $$
{\cal{S}}(s,r,\alpha_1,\alpha_2)={B}(s,r)\cup { B}(x_{0},\alpha_1)
\cup {B}(x_1,\alpha_2),
 $$
 where
  ${ B}(s,r)$ is the closed ball centered at $s$ with radius
  $r$ chosen sufficiently small, so that
  ${\partial{B}}(s,r)\cap \gamma$ contains exactly 2 points
  $x_{0}$ and $x_1$. We suppose that $x_0$, $s$ and $x_1$ are
  arranged counterclockwise.
  \\
   Let $$L_N^{\sigma,\gamma}(s,r,\alpha_1,\alpha_2)= \lim_{t\rightarrow
0}\,\frac{1}{t}\,\left(\BBe_{\sigma}({\vol}({\cal{A}}^{\sigma_{tN^2}}_{N}\cap
{\cal{S}}_N))-{\vol}({\cal{A}}^{\sigma}_{N}\cap
{\cal{S}}_N)\right),$$ where
${\cal{S}}_N=\left({\cal{S}}(s,r,\alpha_1,\alpha_2)\right)_N=\left(
{B}(s,r)\cup B(x_{0},\alpha_1) \cup { B}(x_1,\alpha_2)\right)_N$
and ${\vol}$ denotes the planar Lebesgue measure.
 \newpage
\begin{figure}[hbt]
\vspace{5 cm} \vbox{ \centerline{ \psset{unit=0.85cm}
\pspicture(0,-1)(10,9) \psline[linewidth=2.5pt](8.6,7.1)(6.3,
10.2)(4.2,7.7)
\pspolygon[fillstyle=solid,fillcolor=lightgray,linewidth=1.75pt](9.75,7.5)(9.75,9)(9,9)(9,9.75)(9,12)
(8.25,12)(8.25,12.75)(4.5,12.75)(4.5,12)(3.75,12)(3.75,7.5)(6,7.5)
(6,6.75)(9,6.75)(9,7.5)
\pscurve[linewidth=0.75pt,showpoints=true](10.3,4.7)(8.6,7.1)(7.74,8.2)(6.3,
10.2) (4.7,8.3)(4.2,7.7)(2.4,5.6) (3,3)(7,3.4)(10.3,4.7)
\rput(5.2,8.4){$x_1$} \rput(7.2,8.3){$x_0$}
\rput(6.3, 9.9){$s$}
\psline[linewidth=1.5pt]{->}(6.3, 10.2)(8.5,10.95)
\rput(7.3,10.3){$r$}
\psline[linewidth=1.5pt]{->}(7.74,8.2)(9.09,8.2)
\psline[linewidth=1.5pt]{->}(4.7,8.3)(4.7,7.55)
\rput(8.6,7.8){$\alpha_1$} \rput(5,7.8){$\alpha_2$}
\pscircle[linewidth=1.5pt](7.72,8.21){1.44}
\pscircle[linewidth=1.5pt](4.7,8.3){0.78}
\pscircle[linewidth=1.5pt](6.3,10.2){2.4}
\rput(2.65,7.7){$\gamma$}
\psline[linewidth=1pt]{->}(2.8,7.7)(3,6.43) \psgrid[subgriddiv=0,
unit=0.75, griddots=10,gridlabels=0](1,3)(16,18)
\endpspicture
} \vskip -1.5truecm \centerline{The set $\big(B(s,r)\cup
B(x_0,\alpha_1) \cup B(x_1,\alpha_2)\big)_N$.} }
\end{figure}

\noindent Finally, we define the average
  $${\bf{A}}_N^{\sigma,\gamma}(s,r,\delta)=
  \frac{1}{\delta^2}\int_{0}^{\delta}\int_{0}^{\delta}
  L^{\sigma,\gamma}_N(s,r,\alpha_1,\alpha_2)\,d\alpha_1\,d\alpha_2.$$
\section{Results}
 We first control the quantity ${\bf{A}}_N^{\sigma,\gamma}(s,r,\delta)$ for
 deterministic
sets ${\cal A}^{\sigma}_N$ defined as follows.
\medskip

\noindent
 {\bf{Deterministic initial condition.}} Let $\gamma$ be a Jordan curve of
$\BBr^2$. Suppose that $\gamma$ encloses a connected, compact and
bounded set $\Omega$ of $\BBr^2$, so that $\gamma=\partial
\Omega$. Let $N$ be a fixed positive integer. We define the spin
configuration $\sigma$ at time $0$ as :
$$\forall\ x\in \BBz^2\qquad\sigma(x) = \left\{
 \begin{array}{rl}
 +1\ \ &{\mbox{if}}\ \ \ \Lambda_{x / N}\cap \Omega \neq \emptyset,\\
-1\ \ &{\mbox{otherwise,}} \\
\end{array} \right.
 $$
where, for $x\in\BBz^2$ and $N\in\BBn^*$, $\Lambda_{x /  N}$ is
the box as defined by (\ref{boite}). We will say that $\sigma$ is
the spin configuration associated to the curve $\gamma$ at step
$N$.
\newpage
\noindent Having both the initial condition and the generator, the
Markov process $(\sigma_t)_{t\geq 0}$ at step $N$ is well defined.

\begin{figure}[hbt]
\vskip -1truecm \vbox{ \centerline{ \psset{unit=0.6cm}
\pspicture(0,-1)(10,9)
\pscurve[fillstyle=solid,fillcolor=gray,linewidth=1.75pt](1.2,2.3)(2.3,2.2)(4.5,3.4)(5.7,5.4)(5.2,7.1)(1.2,2.3)
\pspolygon[fillcolor=lightgray,
linewidth=2.5pt](0.75,1.5)(1.5,1.5)(2.25,1.5)(3,1.5)(3,2.25)(3.75,2.25)
(4.5,2.25)(4.5,3)(5.25,3)(5.25,3.75)(6,3.75)(6,4.5)(6,5.25)(6,6)(6,6.75)
(6,7.5)(3.75,7.5)(3.75,6.75)(3,6.75)(3,6)(2.25,6)(2.25,5.25)(1.5,5.25)(1.5,4.5)
(1.5,3.75)(0.75,3.75)(0.75,3.5)(0.75,1.5)
\pspolygon[fillstyle=solid,fillcolor=gray,
linewidth=0pt](-2,3)(-1.5,3)(-1.5,3.5)(-2,3.5)
\pspolygon[fillstyle=solid,fillcolor=white,
linewidth=0pt](-2,4)(-1.5,4)(-1.5,4.5)(-2,4.5)
\rput(-0.65,4.25){${\cal A}_N^\sigma$}
\rput(-0.65,3.25){${\Omega}$}
\endpspicture
} \vskip -1truecm \centerline{The curve ${\gamma}=\partial \Omega$
and the set ${\cal A}_N^\sigma$.} }
\end{figure}
\begin{pro}\label{theo1heu1} Let $\gamma$ be a Jordan curve of $\BBr^2$ of
class ${\cal C}_2$. Suppose that $\gamma$ encloses a connected,
compact and bounded set $\Omega$ of $\BBr^2$.
Let $s$ be a point of $\gamma$. Let $\sigma$ be the spin configuration
associated to the curve $\gamma$ at step $N$.
  Then,
  $$
 \lim_{r\rightarrow 0}\, \,\lim_{\delta\rightarrow
  0}\,\,\liminf_{N\rightarrow +\infty}\, \,\frac{1}{2r}{\bf{A}}^{\sigma, \gamma}_N(s,r,\delta)=\lim_{r\rightarrow 0}\, \,\lim_{\delta\rightarrow
  0}\, \,\limsup_{N\rightarrow +\infty}\, \,\frac{1}{2r}{\bf{A}}^{\sigma, \gamma}_N(s,r,\delta)=
  -\frac{\textstyle 1}{\textstyle 2}\left|\cos
  2\theta\right|\,\xi_{\gamma}(s),
$$
where $\xi_{\gamma}(s)$ is the curvature of $\gamma$ at
$s$ and $\theta$ is the angle between the horizontal axis and the
tangent to the curve $\gamma$ at $s$.
\end{pro}
\noindent We suppose next that the sets
${\cal A}^{\sigma}_N$ are random and that locally the  height
function associated to $\partial {\cal A}^{\sigma}_N$ obeys to
Spohn's initial condition described as follows.
\medskip

\noindent {\bf{Spohn's initial condition.}} Let $\gamma$ be a
Jordan curve of $\BBr^2$. Suppose that $\gamma$ encloses a
connected, compact and bounded set $\Omega$ of $\BBr^2$, so that
$\gamma=\partial \Omega$. Let $s$ be a point of $\gamma$. Suppose
that, on a neighborhood $V_s$ of $s$, the contour $\gamma$ is the
graph of a monotone differentiable function $f$ defined on a
segment $[a,b]$. For each positive integer $N$, and for each
random boundary $\partial {\cal A}^{\sigma}_N \cap V_s$, let
$\Phi_N$ be the random height function associated to $\partial
{\cal A}^{\sigma}_N \cap V_s$ above $\frac{\BBz}{N}\cap [a,b]$,
defined by
$$\forall u\in
\frac{\BBz}{N}\cap [a,b]\qquad
\Phi_N(u)=\sup\{v:\, (u,v)\in \partial {\cal A}^{\sigma}_N\}.
$$
Let $\mu_N$ be the initial distribution of $\Phi_N$. We suppose
that, under $\mu_N$, the increments
$$\Phi_N(\frac{k+1}{N})-\Phi_N(\frac{k}{N})\,,\qquad
\frac{k}{N}\in
[a,b]\cap\frac{ \BBz}{N}\,,$$
are independent and their laws are such
that
\medskip

\noindent $\bullet$ If $f$ is nondecreasing, then for $l\in \BBz$
$$
{\mu_N}\kern-4pt\left(\Phi_N(\frac{k+1}{N})-\Phi_N(\frac{k}{N})=\frac{l}{N}\right)=\left\{
\begin{array}{cccc}
  (f'(\frac{k}{N}))^{l}\,(1+f'(\frac{k}{N}))^{-l-1} & {\mbox{if}}\,\, l\geq 0\\
  0                                               & {\mbox{if}}\,\, l < 0

\end{array}
 \right.
$$
$\bullet$ If $f$ is noninceasing, then for $l\in \BBz$
$$
{\mu_N}\kern-4pt\left(\Phi_N(\frac{k+1}{N})-\Phi_N(\frac{k}{N})=\frac{l}{N}\right)=\left\{
\begin{array}{cccc}
  (|f'|(\frac{k}{N}))^{|l|}\,(1+|f'|(\frac{k}{N}))^{-|l|-1} & {\mbox{if}}\,\, l\leq 0\\
  0                                               & {\mbox{if}}\,\, l >0

\end{array}
 \right.
$$

\begin{pro}\label{P1}
Let $\gamma$ be a Jordan curve of $\BBr^2$ of class ${\cal C}_2$.
Let $s$ be a point of $\gamma$.
Suppose that, for any positive real numbers $r$ and
$\delta$ sufficiently small, the curve $\gamma\cap {\cal
S}(s,r,\delta,\delta)$ is the graph of a monotone
function $f$ defined on a segment $[a,b]$ of $\BBr$.  Let $\mu_N$
be the measure as defined above. Suppose that,
\begin{equation}\label{(i)}
\forall\varepsilon>0\qquad \lim_{N\rightarrow
+\infty}{\mu_N}{\kern-2pt}\left(|\Phi_N(a_N)-f(a)|\geq
\varepsilon\right)=0,
\end{equation}
 where $a_N$ is a point of
$[a,b]\cap\frac{\BBz}{N}$ such that $|a-a_N|\leq \frac{1}{N}$.
Then
  $$\displaylines{
 \lim_{r\rightarrow 0}\, \,\lim_{\delta\rightarrow
  0}\,\,\liminf_{N\rightarrow +\infty}\, \,\frac{1}{2r}
\mu_N\big({\bf{A}}^{\sigma, \gamma}_N(s,r,\delta)\big)
\,=\,\lim_{r\rightarrow 0}\, \,\lim_{\delta\rightarrow
  0}\, \,\limsup_{N\rightarrow +\infty}\, \,\frac{1}{2r}
\mu_N\big({\bf{A}}^{\sigma, \gamma}_N(s,r,\delta)\big)\hfill\cr
=
-\frac{1}{2(|\cos\theta| +
|\sin\theta|)^2}\,\,\xi_{\gamma}(s)\,,
}$$
where $\xi_{\gamma}(s)$ is the curvature of $\gamma$ at
$s$ and $\theta$ is the angle between the horizontal axis and the
tangent to the curve $\gamma$ at $s$.
\end{pro}
\vspace{0.5cm} The limits obtained in propositions \ref{theo1heu1}
and \ref{P1} are very different because the initial conditions
differ. Spohn's velocity is recovered in proposition \ref{P1} (cf.
(4.26) of Spohn (1993)). The choice of the measure $\mu_N$ is the
good one, since as noticed by Spohn (1993), the height differences
are governed by the zero-range process with rate function
$c(n)=\BBone_{n\geq 1}$. The product measure $\mu_N$ with
geometric distribution is invariant for the zero range process
(cf. Andjel (1982)). Motion by mean curvature for the sets $({\cal
A}_N^{\sigma_{N^2t}})$ corresponds then to the hydrodynamic limit
for the zero range process.
\\
\\
\noindent Propositions \ref{theo1heu1} and \ref{P1} are
consequences of the following theorem~\ref{T2}, which handles the
initial conditions described thereafter.
The distance between a point $a\in\BBr^2$ and a subset $B$ of $\BBr^2$
is $d(a,B)=\inf_{b\in B}|a-b|$; the Hausdorff
distance $d_H$ between two subsets $A$ and $B$ of $\BBr^2$ is
$$
d_H(A,B)=\max\left(\sup_{a\in A}d(a,B), \sup_{b\in
B}d(b,A)\right)\,.
$$
 {\bf{Initial condition.}}
 Let $\gamma$ be a Jordan curve of $\BBr^2$ of class ${\cal C}_1$. Suppose that $\gamma$ encloses a connected,
compact and bounded set $\Omega$ of $\BBr^2$.
Let $s$ be a point of
  $\gamma$. Let $r$ be a positive real number sufficiently small
  such that $\partial B(s,r)\cap \gamma$ contains exactly two
  points
  $x_0$ and $x_1$. Suppose that $x_0$, $s$ and $x_1$ are arranged counterclockwise.
Let $\theta_1\in [0,2\pi]$
(respectively $\theta_{0}\in [0,2\pi]$)
 be the oriented angle between the half horizontal axis $[0,+\infty[$
 and $T_{x_1}\gamma$
 (respectively $T_{x_0}\gamma$).
We suppose that there exists a neighborhood
$V_s$ of s and a probability measure $\nu_N$ such that
\begin{equation}\label{haus}
\forall\varepsilon>0\qquad
\lim_{N\rightarrow +\infty}\nu_N\left(d_H({\cal A}^{\sigma}_N\cap
V_s, \Omega\cap V_s)\geq \varepsilon\right)=0,
\end{equation}
and that, with probability one, the boundaries $\gamma$
and $\partial {\cal A}^{\sigma}_N$ are, in $V_s$, either
both non-increasing or either both non-decreasing.

\begin{figure}[hbt]
\vskip 9truecm \vbox{ \centerline{ \psset{unit=1.75cm}
\pspicture(2.5,2.5)(5,-0.5)
\psline[linewidth=0.5pt,linecolor=red](6,7.5)(3,4.5)(5.1,2.4)(6,1.5)
\psline[linewidth=0.5pt,linecolor=red,linestyle=dashed](3,4.5)(3.25,4.5)(3.5,4.5)(4.5,4.5)
\psarc[arcsepA=2pt,arcsepB=2pt]{->}(3,4.5){0.8}{0}{45}
\rput(3.65,4.75){$\theta_0$}
\psarc[arcsepA=2pt,arcsepB=2pt]{->}(3,4.5){0.5}{0}{315}
\rput(3,3.75){$\theta_1$}
\psline[linewidth=2pt,linecolor=red](5.5,6.5)(5,6.5)(4.75,6.5)(4.5,6.5)(4.25,6.5)(4.25,6.25)(4.25,6)(3.9,6)(3.9,5)
(3.9,4)(4.25,4)(4.25,3)(4.75,3)(4.75,2.25)(5.5,2.25)
\pscurve[linewidth=2pt,linecolor=blue,showpoints=true](6.5,6.8)(6,6.75)(5.25,6.5)(4.5,6)(4,4.5)(4.5,3)(5.25,2.6)(6.75,3)(6.5,6.8)
\rput(6,7.7){$T_{x_0}\gamma$} \rput(6.4,1.52){$T_{x_1}\gamma$}
\rput(3.7,4.3){$s$} \rput(4.5,6.2){$x_0$} \rput(4.5,2.7){$x_1$}
\psarc[arcsepA=2pt,arcsepB=2pt](4,4.5){1.6}{0}{364}
\psarc[arcsepA=2pt,arcsepB=2pt](4.5,6){0.9}{0}{364}
\psarc[arcsepA=2pt,arcsepB=2pt](4.5,3){0.88}{0}{364}
\psline[linewidth=2pt,linecolor=blue]{->}(6.6,2)(6.3,2.7)
\rput(6.8,2){$\gamma$}
\endpspicture
} } \vskip -2.5truecm {\centerline{The polygonal curve $
\partial {\cal A}^{\sigma}_N$ behaves in $V_s$ as $\gamma$. }}
\end{figure}

\vspace{1cm} \noindent Let, for $x\in \gamma\cap V_s$ and
$\delta>0$,
$$C_N(x,\delta)=\sum \BBone_{\sigma(y)=+1,\,s(\sigma,y)=2},$$
where the sum is taken over all $y\in \BBz^2$ for which
$\frac{\textstyle y}{\textstyle N}$ is a point of
$(B(x,\delta))_N\setminus B(s,|x-s|)$. The quantity
$C_N(x,\delta)$ is equal to half of the number of the corners of
the polygonal line $\partial {\cal A}^{\sigma}_N$ belonging to
$(B(x,\delta))_N \setminus B(s,|x-s|)$.
\newpage
\noindent We first suppose that $\gamma$ is a polygon and that $s$
is a corner point of $\gamma$. In this case, the following theorem
proves that, for $r$ and $\delta$ sufficiently small, the limit as
$N$ goes to infinity of
$\nu_N({\bf{A}}^{{\sigma,\gamma}}_N(s,r,\delta))$ exists under  a
suitable behavior of the expected proportions of corners
$\frac{1}{N}\nu_N(C_N(x_k,\delta))$, for $k\in\{0,1\}$.

\begin{theo}\label{T1}
Let $\gamma,s,r,\delta$ and $\nu_N$ be as described in the
previous initial condition. Suppose that $\gamma$ is a polygon and
that for $k=0,1$ and for $r,\delta$ sufficiently small, the
following limit holds:
\begin{equation}\label{c5}
\lim_{N\rightarrow
\infty}\frac{1}{N}{\nu_N}(C_N(x_k,\delta))={\delta}C(\theta_k)\,.
\end{equation}
Then, for $r$ and $\delta$ sufficiently small, one has
\begin{eqnarray}\label{th1}
{\lefteqn{
\lim_{N\rightarrow\infty}{\nu_N}({\bf{A}}^{{\sigma,\gamma}}_N(s,r,\delta))=
}}\nonumber\\
&& =-\frac{1}{2}\,\sg(\tan\theta_0)\left(\cos^2\theta_0
  +
{C(\theta_0)}\left(\left|\sin\theta_0\right|-
\left|\cos\theta_0\right|\right)\right)\,
\nonumber\\
&& +\,\frac{1}{2}\,\sg(\tan\theta_1)\left(\cos^2\theta_1
  +\,
{C(\theta_1)}\left(\left|\sin\theta_1\right|-
\left|\cos\theta_1\right|\right)\right)\nonumber \\
&& +\,
\BBone_{\sin\theta_0\sin\theta_1>0}\left(\sg(\theta_1-\theta_0)\BBone_{\cos\theta_0\cos\theta_1>0}+\,
\sg(\tan\theta_0)\BBone_{\cos\theta_0\cos\theta_1<0}\right).
  \end{eqnarray}
\end{theo}
\vspace{1cm}
 Suppose that
$C(\theta)=f(|\sin\theta|,|\cos\theta|)$, where $f$ is a positive
function defined on $[0,1]\times[0,1]$ and that $s$ is not a
corner point of the polygon $\gamma$. Theorem \ref{T1} then
implies that, for $r$ and $\delta$ sufficiently small, the limit
as $N$ goes to infinity of
$\nu_N({\bf{A}}^{{\sigma,\gamma}}_N(s,r,\delta))$ vanishes (since
in this case $\theta_1=\theta_0\pm\pi$). This constatation is not
surprising since the inverse of the curvature of a straight line
vanishes.
\\
\\
The following theorem extends theorem \ref{T1} to Jordan curves.
\begin{theo}\label{T2}
Let $\gamma,s,r$ and $\nu_N$ be as described in the previous
initial condition. Suppose that for $r$ sufficiently small and for
$k=0,1$, the following limits exist:
\begin{equation}\label{c5bis}
\lim_{\delta\rightarrow
  0}\,\,\liminf_{N\rightarrow +\infty}\, \,
\frac{1}{\delta N}{\nu_N}(C_N(x_k,\delta))=
\lim_{\delta\rightarrow
  0}\,\,\limsup_{N\rightarrow +\infty}\, \,
\frac{1}{\delta N}{\nu_N}(C_N(x_k,\delta))=
C(\theta_k)\,.
\end{equation}
Then
\begin{eqnarray*}
\,\lim_{\delta\rightarrow
  0}\,\,\liminf_{N\rightarrow +\infty}\,
\nu_N\big({\bf{A}}^{\sigma, \gamma}_N(s,r,\delta)\big) \,=\,
\,\lim_{\delta\rightarrow
  0}\, \,\limsup_{N\rightarrow +\infty}\,
\nu_N\big({\bf{A}}^{\sigma, \gamma}_N(s,r,\delta)\big),
  \end{eqnarray*}
  the common value is as in (\ref{th1}) with the function $C(.)$
  given by (\ref{c5bis}).
\end{theo}

\newpage
\section{Proofs}
We first prove theorems \ref{T1} and \ref{T2}. Next, we prove the
two propositions. For the proof of the theorems, we need the
following preliminary lemma.
\begin{lemma}\label{lem}
Let ${\cal{S}}$ be a compact set of $\BBr^2$. Let $\sigma \in
\{-1,+1\}^{\BBz^2}$ be fixed. Then
$$\lim_{t\rightarrow
0}\,\frac{1}{t}\,\left(\BBe_{\sigma}({\vol}({\cal{A}}^{\sigma_{tN^2}}_{N}\cap
{\cal{S}}_N))-{\vol}({\cal{A}}^{\sigma}_{N}\cap
{\cal{S}}_N)\right)\, = $$ $$ \sum_{x\in\BBz^2:\,
\Lambda_{\frac{x}{N}}\subset
{{\cal{S}}_N}}\kern-20pt\left(\BBone_{\sigma(x)=-1,\
s(\sigma,x)\geq 3}-\BBone_{\sigma(x)=+1,\ s(\sigma,x)\geq
3}\right) + \alpha\kern-20pt \sum_{x\in\BBz^2:\,
\Lambda_{\frac{x}{N}}\subset
{{\cal{S}}_N}}\kern-20pt\left(\BBone_{\sigma(x)=-1,\
s(\sigma,x)=2}-\BBone_{\sigma(x)=+1,\ s(\sigma,x)=2}\right). $$
\end{lemma}
\noindent {\bf{Proof of lemma \ref{lem}.}} Let
$f_N(\sigma)={\vol}({\cal{A}}^{\sigma}_{N}\cap {\cal{S}}_N)$ and
$S(t)f_N(\sigma)=\BBe_{\sigma}({\vol}({\cal{A}}^{\sigma_{t}}_{N}\cap
{\cal{S}}_N))$. We deduce from
$$
\lim_{t\rightarrow 0}\,\frac{1}{t}\left(S(t)f_N-f_N\right)=Lf_N,
$$
that
\begin{equation}\label{lim1he}
\lim_{t\rightarrow
0}\,\frac{1}{t}\left(S(t)f_N(\sigma)-f_N(\sigma)\right)=\sum_{x\in
\BBz^2}c(x,\sigma)(f_N(\sigma^{x})-f_N(\sigma)).
\end{equation}
Now,
$$
f_N(\sigma^{x})-f_N(\sigma)=\frac{1}{N^2}\BBone_{\Lambda_{\frac{x}{N}}\subset
{\cal{S}}_N}\left(\BBone_{\sigma(x)=-1}-\BBone_{\sigma(x)=1}\right),
$$
this fact together with (\ref{lim1he}) gives
$$
\lim_{t\rightarrow
0}\,\frac{1}{t}\left(S(tN^2)f_N(\sigma)-f_N(\sigma)\right)=
\sum_{x\in\BBz^2:\, \Lambda_{\frac{x}{N}}\subset
{{\cal{S}}_N}}c(x,\sigma)\left(\BBone_{\sigma(x)=-1}-\BBone_{\sigma(x)=1}\right),
$$
which proves lemma \ref{lem} since
$c(x,\sigma)=\BBone_{s(\sigma,x)\geq 3}+\alpha
\BBone_{s(\sigma,x)=2}$. \ \ \ $\Box$

\subsection{Evaluation of
$\nu_N\left(L^{\sigma,\gamma}_N(s,r,\alpha_1,\alpha_2)\right)$}
 Throughout this
step, we consider the set
\begin{equation}\label{s}
{\cal{S}}_N=\left({\cal{S}}(s,r,\alpha_1,\alpha_2)\right)_N=\left({
B}(s,r)\cup {B}(x_{0},\alpha_1) \cup { B}(x_1,\alpha_2)\right)_N,
\end{equation}
 where $\alpha_1,\alpha_2$ are positive
real numbers less than $\delta$, the positive real numbers $r$ and
$\delta$ are small enough so that
  ${\partial{B}}(s,r)\cap \gamma$ contains exactly 2 points
  $x_{0}$ and $x_1$.
\medskip

\noindent
 The boundary of ${\cal A}^{\sigma}_{N}$ which is included in ${\cal{S}}_N$ can be
described as a sequence $v_1,\ldots,v_r$ of horizontal or vertical
vectors of norm $\frac{1}{N}$, enumerated counterclockwise. We
denote by ${e}_N^1(\alpha_1)$, ${e}_N^2(\alpha_2)$ the two unit
vectors defined by
\begin{equation}\label{v}
{e}_N^1(\alpha_1)={N}{v}_1,\ \ \ \ \ \ \
{e}_N^2(\alpha_2)={N}{v}_{r},
\end{equation}

\newpage
and by ${\cal L}^{\sigma}_N$ the maximal subgraph of
$\partial{\cal A}^{\sigma}_{N}$ included in ${\cal{S}}_N$:
\begin{equation}\label{vbis}
{\cal L}^{\sigma}_N=(v_1,\ldots,v_{r}).
\end{equation}

\begin{figure}[hbt]
\vspace{1.7cm}
 \vbox{ \centerline{ \psset{unit=0.25cm} \pspicture(0,-1)(70,41)
\psline[linewidth=2.23pt](8,6)(8,8)(10,8)(10,10)(12,10)(14,10)(14,12)(16,12)
(16,14)(20,14)(22,14)(22,16)(22,18)(26,18)(28,18)(28,20)(28,22)
(38,22)(48,22)(48,24)(54,24)(54,26)(62,26)
\rput(63,28){$v_1$} \psline[linewidth=1.5pt]{->}(62,26)(60,26)
\psline[linewidth=2pt]{->}(8,8)(8,6) \rput(4.5,7){$v_r$}
\pscircle(30.9,20.8){20}
\psline[linewidth=1pt]{->}(30.9,20.8)(24.71,39.82)
\rput(29,31){$r$} \psline[linewidth=1pt]{->}(34.9,29)(37.6,22)
\rput(34,30){${\cal L}^{\sigma}_N$}
\psdots[dotscale=1.2](30.9,20.8) \rput(30.9,22.8){$s$}
\psdots[dotscale=1.2](50.72,23.46) \rput(52.2,25){$x_{0}$}
\pscircle(50.72,23.46){10}
\psline[linewidth=1pt]{->}(50.72,23.46)(58.32,16.96)
\rput(54.5,18.5){$\alpha_{1}$}
\psline[linewidth=1pt](30.9,20.8)(65.6,25.46)
\psdots[dotscale=1.2](14.36,9.56) \rput(12.7,10.8){$x_{1}$}
\pscircle(14.36,9.56){6.2}
\psline[linewidth=1pt]{->}(14.36,9.56)(17.54,3.32)
\rput(15,5){$\alpha_{2}$}
\psline[linewidth=1pt](30.9,20.8)(4.43,2.81)
\psgrid[subgriddiv=1,unit=2,griddots=4,gridlabels=0](0,-1)(35,21)
\endpspicture
} \vskip 1truecm \centerline{The polygonal line ${\cal
L}^{\sigma}_N=(v_1,\cdots,v_r)$. } \centerline{ Here ${\cal
S}_N=\big(B(s,r)\cup B(x_{0},\alpha_1) \cup
B(x_{1},\alpha_2)\big)_N$.} }
\end{figure}
\vspace{1cm} \noindent We now need the following definition and
notation.
\begin{defi}
We say that ${\cal L}_N$ is a path on $\BBz_N^2$ if ${\cal L}_N$
is a finite sequence of consecutive vectors $(v_i)_{1\leq i\leq
r}$ (this means that the endpoint of $v_i$ is the starting point
of $v_{i+1}$ for $1\leq i<r$) of norm $1/N$, drawn on the grid
$\BBz_N^2$, and such that the endpoints of these vectors (resp.
the starting points) are distinct.
\end{defi}
\newpage
\vspace{-1.75cm}
\begin{figure}[hbt]
\centerline{The following family of vectors $(v_1,\ldots,v_r)$ is
a
 path
on the grid $\BBz_N^2$.} \vskip -1truecm \vbox{ \centerline{
\psset{unit=0.7cm} \pspicture(0,-1)(10,9)
\psgrid[subgriddiv=0,griddots=10,gridlabels=0](7,7)
\psline[linewidth=1.5pt]{->}(5,1)(6,1)
\psline[linewidth=1.5pt]{->}(6,1)(6,2)
\psline[linewidth=1.5pt]{->}(6,2)(6,3)
\psline[linewidth=1.5pt]{->}(6,3)(6,4)
\psline[linewidth=1.5pt]{->}(6,4)(5,4)
\psline[linewidth=1.5pt]{->}(5,4)(4,4)
\psline[linewidth=1.5pt]{->}(4,4)(4,3)
\psline[linewidth=1.5pt]{->}(4,3)(3,3)
\psline[linewidth=1.5pt]{->}(3,3)(3,4)
\psline[linewidth=1.5pt]{->}(3,4)(3,5)
\psline[linewidth=1.5pt]{->}(3,5)(4,5)
\psline[linewidth=1.5pt]{->}(4,5)(4,6)
\psline[linewidth=1.5pt]{->}(4,6)(3,6) \rput(5.5,0.74){$v_1$}
\rput(6.3,1.5){$v_2$} \rput(3.5,6.3){$v_r$} \rput(4.5,4.2){$v_i$}
\endpspicture
} }
\end{figure}
\vspace{-0.8cm} \hspace{-0.7cm}{\bf{Notation.}} Let ${\cal
L}_N=(v_1,v_2,\ldots,v_r)$ be a path on $\BBz_N^2$. We define
\begin{equation}\label{def}
N_{+}\left({\cal L}_N\right)=\card\left\{i:\,
\widehat{(v_i,v_{i+1})}=-\frac{\pi}{2}\right\},\ \
N_{-}\left({\cal L}_N\right)=\card\left\{i:\,
\widehat{(v_i,v_{i+1})}=+\frac{\pi}{2}\right\},
\end{equation}
 where
$\widehat{(v_i,v_{i+1})}$ denotes the oriented angle between $v_i$
and $v_{i+1}$.
\medskip

\noindent
The purpose of the following proposition is to establish the
relation between $N_{-}({\cal L}^{\sigma}_N)-N_{+}({\cal
L}^{\sigma}_N)$ and $L_N^{\sigma,\gamma}(s,r,\alpha_1,\alpha_2)$,
for the path ${\cal L}^{\sigma}_N$ as defined by (\ref{vbis}).
\begin{pro}\label{a1}
Let $N$ be a fixed positive integer. Let ${\cal L}^{\sigma}_N$ be
the random path as defined by (\ref{vbis}). Then
\begin{equation}\label{egalite}
\nu_N\left(L_N^{\sigma,\gamma}(s,r,\alpha_1,\alpha_2)\right)=\frac{1}{2}\nu_N\left(N_{-}({\cal
L}^{\sigma}_N)-N_{+}({\cal L}^{\sigma}_N)\right).
\end{equation}
\end{pro}
{\bf{Proof of proposition \ref{a1}.}} Let $N\in \BBn^*$ be fixed
and ${\cal{S}}_N=\left({ B}(s,r)\cup {B}(x_{0},\alpha_1) \cup {
B}(x_1,\alpha_2)\right)_N$. Let $f$ be the function defined from
$\{0,1,\ldots,4\}$ to $\{0,1,2\}$ by
$$f(s(\sigma,x))=\left\{
\begin{array}{rl}
1\ \ & {\mbox{if}}\  s(\sigma,x)=2 \\ 2\ \ & {\mbox{if}}\
s(\sigma,x)=3\\ 
0\ \ & otherwise.
\end{array}
\right.
 $$
On the one hand, by definition of $N_{-}({\cal L}^{\sigma}_N)$ and
$N_{+}({\cal L}^{\sigma}_N)$, we have
\begin{equation}\label{h1}
\sum_{x\,\in\,\BBz^2\,\cap\, N{\cal{S}}_N}\sigma(x)f(s(\sigma,x))
=N_{+}({\cal L}^{\sigma}_N)-N_{-}({\cal L}^{\sigma}_N),
\end{equation}
 on the other hand, we deduce from the definition of the function
 $f$,
\begin{eqnarray*}
{\lefteqn{\sum_{x\,\in\,\BBz^2\,\cap\,
N{\cal{S}}_N}\sigma(x)f(s(\sigma,x))
 = - \kern -15pt\sum_{x\,\in\,\BBz^2,\,\Lambda_{x/N}\subset
{\cal{S}}_N}\kern
-20pt\left(\BBone_{\sigma(x)=-1\,,\,s(\sigma,x)=2}-\BBone_{\sigma(x)=+1,\,s(\sigma,x)=2}\right)}}\\
&& -\, 2\kern -15pt\sum_{x\,\in\,\BBz^2,\,\Lambda_{x/N}\subset
{\cal{S}}_N}\kern
-20pt\left(\BBone_{\sigma(x)=-1,\,s(\sigma,x)=3}-\BBone_{\sigma(x)=+1,\,s(\sigma,x)=3}\right).
\end{eqnarray*}
We combine the last formula, lemma \ref{lem} (with $\alpha =
\frac{\textstyle 1}{\textstyle 2}$) together with the fact that
$\BBone_{s(\sigma,x)=4}=0$, and we obtain
\begin{equation}\label{h2}
\nu_N\left(L_N^{\sigma,\gamma}(s,r,\alpha_1,\alpha_2)\right)=
-\frac{1}{2}\nu_N\left(\sum_{x\in\BBz^2\cap
{\cal{S}}_N}\sigma(x)f(s(\sigma,x))\right).
\end{equation}
The statement of proposition \ref{a1} follows from (\ref{h1}) and
(\ref{h2}) by taking the expectation with respect to~$\nu_N$.\,
$\Box$
\medskip

\noindent
In view of proposition \ref{a1}, in order to control
$\nu_N\left(L_N^{\sigma,\gamma}(s,r,\alpha_1,\alpha_2)\right)$, it
remains to evaluate $\nu_N\left(N_{+}\left({\cal
L}^{\sigma}_N\right)-N_{-}\left({\cal
L}^{\sigma}_N\right)\right)$. For this, we begin by controlling
the quantity $N_{+}\left({\cal L}_N\right)-N_{-}\left({\cal
L}_N\right)$ for monotone deterministic paths ${\cal L}_N$ defined
as follows.

\begin{defi}
A path on $\BBz_N^2$ is said to be monotone if all its horizontal as
well as all its vertical vectors are oriented in the same sense.
\end{defi}
\vspace{1.7cm}
\begin{figure}[hbt]
\vbox{ \centerline{ \psset{unit=0.45cm} \pspicture(0,-1)(8,1)
\psgrid[subgriddiv=0,griddots=10,gridlabels=0](5,5)
\psline[linewidth=1.5pt]{->}(4,4)(3,4)
\psline[linewidth=1.5pt]{->}(3,4)(2,4)
\psline[linewidth=1.5pt]{->}(2,4)(2,3)
\psline[linewidth=1.5pt]{->}(2,3)(1,3)
\psline[linewidth=1.5pt]{->}(1,3)(1,2)
\psline[linewidth=1.5pt]{->}(1,2)(1,1)
\endpspicture
} \vskip -0.5truecm \centerline{A monotone path on the grid
$\BBz_N^2$.} }
\end{figure}
The following lemma evaluates $N_{+}\left({\cal
L}_N\right)-N_{-}\left({\cal L}_N\right)$, whenever ${\cal L}_N$
is a monotone path on $\BBz_N^2$.
\begin{lemma}\label{base}
Let $(v_i)_{1\leq i\leq r}$ be a sequence of $r$ consecutive
vectors drawn on the grid $\BBz_N^2$. These vectors are enumerated
beginning from $N^{-1}u_e:=v_1$ until $N^{-1}u_s:=v_r$. We suppose
that they form a {\bf monotone} path on $\BBz_N^2$, say ${\cal
L}_N$. Let $[u_e \wedge u_s]= (u_e\cdot i)(u_s\cdot j) - (u_e\cdot
j)(u_s\cdot i)$. Then $$ N_{+}\left({\cal
L}_N\right)-N_{-}\left({\cal L}_N\right)=[u_e \wedge u_s]. $$
\end{lemma}

\begin{figure}[hbt]
\vskip -1truecm \vbox{ \centerline{ \psset{unit=0.59cm}
\pspicture(0,-1)(10,9)
\psgrid[subgriddiv=0,griddots=10,gridlabels=0](-1,0)(8,8)
\psline[linewidth=1.5pt]{->}(7,7)(6,7)
\psline[linewidth=1.5pt]{->}(6,7)(6,6)
\psline[linewidth=1.5pt]{->}(6,6)(5,6)
\psline[linewidth=1.5pt]{->}(5,6)(4,6)
\psline[linewidth=1.5pt]{->}(4,6)(4,5)
\psline[linewidth=1.5pt]{->}(4,5)(4,4)
\psline[linewidth=1.5pt]{->}(4,4)(3,4)
\psline[linewidth=1.5pt]{->}(3,4)(2,4)
\psline[linewidth=1.5pt]{->}(2,4)(2,3)
\psline[linewidth=1.5pt]{->}(2,3)(1,3)
\psline[linewidth=1.5pt]{->}(1,3)(1,2)
\psline[linewidth=1.5pt]{->}(1,2)(1,1) \rput(6.5,7.5){$u_e/N$}
\rput(0.1,1.5){$u_s/N$}
\psarc[arcsepA=2pt,arcsepB=2pt]{<-}(6,7){0.5}{270}{0}
\rput(6.5,6.5){$1$}
\psarc[arcsepA=2pt,arcsepB=2pt]{->}(6,6){0.3}{90}{180}
\rput(5.5,6.5){$-1$}
\psarc[arcsepA=2pt,arcsepB=2pt]{<-}(4,6){0.5}{270}{0}
\rput(4.5,5.5){$1$}
\psarc[arcsepA=2pt,arcsepB=2pt]{<-}(2,4){0.5}{270}{0}
\rput(2.5,3.5){$1$}
\psarc[arcsepA=2pt,arcsepB=2pt]{<-}(1,3){0.5}{270}{0}
\rput(1.5,2.5){$1$}
\psarc[arcsepA=2pt,arcsepB=2pt]{->}(4,4){0.3}{90}{180}
\rput(3.5,4.5){$-1$}
\psarc[arcsepA=2pt,arcsepB=2pt]{->}(2,3){0.3}{90}{180}
\rput(1.5,3.5){$-1$}
\endpspicture
} \centerline {For this monotone path ${{\cal L}}_N$, we have
$u_e=(-1,0)$ and $u_s=(0,-1)$,} \centerline { hence $[u_e\wedge
u_s]=1$. On the other hand $N_{+}({{\cal L}}_N)-N_{-}({{\cal
L}}_N) = 1-1+1-1+1-1+1=1.$ } }
\end{figure}

\noindent {\bf{Remark.}} Let us note that for any path ${\cal
L}_N=(v_1,\ldots,v_r)$, we have $$\widehat{(u_e,u_s)}
=\frac{\pi}{2}\left(N_{-}\left({\cal L}_N\right)-N_{+}\left({\cal
L}_N\right)\right),$$ where $u_e=Nv_1$ and $u_s=Nv_r$.
\\
 {\bf{Proof of lemma \ref{base}.}} We denote
by ${\cal L}_N(r)=(v_1,\ldots,v_r)$
   a monotone path on $\BBz_N^2$.
  The proof of lemma \ref{base}
  is done by
induction on $r$.
\\
 For $r=1$, we have $N_{-}\left({\cal
L}_N(1)\right)-N_{+}\left({\cal L}_N(1)\right)=0$
  which corresponds to $[u_e \wedge u_s]$, since in this case
  $N^{-1}u_e=N^{-1}u_s=v_1$. \\
We suppose now that the property is true at step $r\geq 1$ and we
prove it at step $r+1$. We consider the path
  ${\cal L}_N(r+1)$. Since ${\cal L}_N(r+1)$ is monotone, we can suppose
  without loss of generality that
  $$
{(\cal{H})}\ \ \ \ \ \ \ \forall\, l\in\{1,\ldots, r+1\}\qquad
(Nv_l)\cdot i\in\{0,-1\},\ \ \  (Nv_l)\cdot j\in\{0,-1\}. $$ Once
the hypothesis ${(\cal{H})}$ is assumed, we have only three cases
to discuss on the expression of $(v_r,v_{r+1})$,
\\
  $\bullet$ If $v_r=v_{r+1}$, then  $
  N_{+}\left({\cal L}_N(r+1)\right)-N_{-}\left({\cal L}_N(r+1)\right)=
  N_{+}\left({\cal
  L}_N(r)\right)-N_{-}\left({\cal
  L}_N(r)\right)$, and the inductive assumption gives
  $$
  N_{+}\left({\cal L}_N(r+1)\right)-N_{-}\left({\cal L}_N(r+1)\right) =
[Nv_1\wedge Nv_{r+1}].
  $$
  $\bullet$ If $(Nv_r)\cdot j=-1=(Nv_{r+1})\cdot i$, then
$\widehat{(v_r,v_{r+1})}=\frac{\textstyle \pi}{\textstyle 2}$ and
$
  N_{+}\left({\cal L}_N(r+1)\right)-N_{-}\left({\cal L}_N(r+1)\right)=
  N_{+}\left({\cal
  L}_N(r)\right)-N_{-}\left({\cal
  L}_N(r)\right)-1$. Together with the inductive assumption, this gives
  $$
  N_{+}\left({\cal L}_N(r+1)\right)-N_{-}\left({\cal L}_N(r+1)\right) =
-(Nv_1)\cdot i-1=(Nv_1)\cdot j=[Nv_1\wedge Nv_{r+1}].
  $$
  $\bullet$ If $(Nv_r)\cdot i=-1=(Nv_{r+1})\cdot j$, then $\widehat{(v_r,v_{r+1})}=-
  \frac{\textstyle \pi}{\textstyle 2}$, $
  N_{+}\left({\cal L}_N(r+1)\right)-N_{-}\left({\cal L}_N(r+1)\right)=
  N_{+}\left({\cal
  L}_N(r)\right)-N_{-}\left({\cal
  L}_N(r)\right)+1$
  and
  $$
  N_{+}\left({\cal L}_N(r+1)\right)-N_{-}\left({\cal L}_N(r+1)\right) =
(Nv_1)\cdot j+1=-(Nv_1)\cdot i=[Nv_1\wedge Nv_{r+1}].
  $$
The equality $
  N_{+}\left({\cal L}_N(r+1)\right)-N_{-}\left({\cal L}_N(r+1)\right)=[Nv_1\wedge Nv_{r+1}]
  $ is then always valid and lemma \ref{base} is proved.\ \  $\Box$
\medskip

\noindent
  The following lemma generalizes lemma \ref{base}. Its purpose is
  to evaluate $
  N_{+}\left({\cal L}_N\right)-N_{-}\left({\cal L}_N\right)$ for
  a path ${\cal L}_N$ constructed by concatenating two monotone paths.
  \begin{lemma}\label{basebis}
Let ${\cal L}_N=(v_1,\ldots,v_r,w_1,\ldots,w_{s})$ be a path on
$\BBz_N^2$.
Suppose that $(v_1,\ldots,v_r)$
(respectively $(w_1,\ldots,w_{s})$) forms
 a monotone path on $\BBz_N^2$ and that $v_r\cdot
w_1=0$. Let $a_1,a_2,b_1,b_2\in\{-1,+1\}$. Suppose that for each
$1\leq i\leq r$ (resp. $1\leq j\leq s$), the vector $Nv_i$ (resp.
$Nw_j$) is either $(a_1,0)$ (resp. $(b_1,0)$) or $(0,a_2)$ (resp.
$(0,b_2)$). Then,
\begin{equation}\label{abb}
  N_{-}\left({\cal
  L}_N\right)-N_{+}\left({\cal
  L}_N\right)=-a_2(Nv_1)\cdot i +
  b_2(Nw_s)\cdot i + f(a_1,a_2,b_1,b_2),
  \end{equation}
where $i$ is the unit vector $(1,0)$, $\cdot$ is the usual scalar
product in $\BBr^2$ and $$ f(a_1,a_2,b_1,b_2)= \left\{
\begin{array}{rl}
2a_1a_2\ \ & {{if}}\ \ a_2b_2=-1,\,\ \left((Nv_r)\cdot i=a_1\,{\mbox{or}}\,\, a_1b_1=1\right) \\
2b_1a_2\ \ & {{if}}\ \ a_2b_2=-1,\,\  (Nv_r)\cdot i=0\\
0\ \  & \ \ {{if}}\ \ a_2b_2=1.
\end{array}
\right. $$
  \end{lemma}
\noindent{\bf{Proof of lemma \ref{basebis}.}}
  We deduce, applying lemma \ref{base} to the monotone paths
$(v_1,\ldots,v_r)$, $(w_1,\ldots,w_{s})$ and $(v_r,w_1)$ that, for
${\cal L}_N=(v_1,\ldots,v_r,w_1,\ldots,w_{s})$,
\begin{equation}\label{somme}
  N_{+}\left({\cal L}_N\right)-N_{-}\left({\cal L}_N\right)=
  [Nv_1\wedge Nv_{r}]+[Nv_r\wedge Nw_1]
  +[Nw_1\wedge Nw_{s}].
  \end{equation}

\noindent  {In the following picture, we have $Nv_1=(1,0)$,
$Nv_r=(0,-1)$, $Nw_1=(-1,0)$, $Nw_s=(-1,0)$. Hence $[Nv_1\wedge
Nv_r]+[Nw_1\wedge Nw_s]+[Nv_r\wedge Nw_1]=-2.$ On the other hand,
we have $N_{+}({{\cal L}}_N)-N_{-}({{\cal L}}_N)=-2$.}

\begin{figure}[hbt]
\vbox{ \centerline{ \psset{unit=0.7cm} \pspicture(0,-1)(6,9)
  \psgrid[subgriddiv=0,griddots=10,gridlabels=0](-3,0)(6,9)
\psline[linewidth=1.5pt]{->}(1,8)(2,8)
\psline[linewidth=1.5pt]{->}(2,8)(2,7)
\psline[linewidth=1.5pt]{->}(2,7)(2,6)
\psline[linewidth=1.5pt]{->}(2,6)(3,6)
\psline[linewidth=1.5pt]{->}(3,6)(4,6)
\psline[linewidth=1.5pt]{->}(4,6)(4,5)
\psline[linewidth=1.5pt]{->}(4,5)(5,5)
\psline[linewidth=1.5pt]{->}(5,5)(5,4)
\psline[linewidth=1.5pt]{->}(5,4)(4,4)
\psline[linewidth=1.5pt]{->}(4,4)(3,4)
\psline[linewidth=1.5pt]{->}(3,4)(2,4)
\psline[linewidth=1.5pt]{->}(2,4)(2,3)
\psline[linewidth=1.5pt]{->}(2,3)(1,3)
\psline[linewidth=1.5pt]{->}(1,3)(1,2)
\psline[linewidth=1.5pt]{->}(1,2)(1,1)
\psline[linewidth=1.5pt]{->}(1,1)(0,1) \rput(1.5,8.5){$v_1$}
\rput(2.3,7.5){$v_2$} \rput(1.6,6.5){$v_3$} \rput(5.5,4.5){$v_r$}
\rput(4.5,3.5){$w_1$} \rput(3.5,3.5){$w_2$} \rput(0.5,0.5){$w_s$}
\psarc[arcsepA=2pt,arcsepB=2pt]{->}(2,8){0.3}{180}{270}
\rput(1.5,7.5){$-1$}
\psarc[arcsepA=2pt,arcsepB=2pt]{<-}(2,6){0.3}{360}{90}
\rput(2.5,6.5){$1$}
\psarc[arcsepA=2pt,arcsepB=2pt]{->}(4,6){0.3}{180}{270}
\rput(3.5,5.5){$-1$}
\psarc[arcsepA=2pt,arcsepB=2pt]{<-}(4,5){0.3}{360}{90}
\rput(4.5,5.5){$1$}
\psarc[arcsepA=2pt,arcsepB=2pt]{->}(5,5){0.3}{180}{270}
\rput(4.5,4.6){$-1$}
\psarc[arcsepA=2pt,arcsepB=2pt]{->}(5,4){0.3}{90}{180}
\rput(4.2,4.3){$-1$}
\psarc[arcsepA=2pt,arcsepB=2pt]{<-}(2,4){0.3}{270}{360}
\rput(2.5,3.5){$1$}
\psarc[arcsepA=2pt,arcsepB=2pt]{->}(2,3){0.3}{90}{180}
\rput(1.5,3.5){$-1$}
\psarc[arcsepA=2pt,arcsepB=2pt]{->}(1,1){0.3}{90}{180}
\rput(0.5,1.5){$-1$}
\psarc[arcsepA=2pt,arcsepB=2pt]{<-}(1,3){0.3}{270}{360}
\rput(1.5,2.5){$1$}
\endpspicture
} }
\end{figure}

\noindent We deduce from
$$ (Nv_l)\cdot i\in\{0,a_1\},\ \ \ {\mbox{and}}\ \ \ (Nv_l)\cdot j\in\{0,a_2\},$$for $1\leq l\leq r$, that
$$
 a_1 (Nv_l)\cdot i+ a_2 (Nv_l)\cdot j=1.
$$
This fact gives
\begin{equation}\label{basebisbis}
[Nv_1\wedge Nv_{r}]=a_2(Nv_1)\cdot i -a_2 (Nv_r)\cdot i.
\end{equation}
In the same way, we deduce that for any $ 1\leq l\leq s$,
\begin{equation}\label{basebisbisbis}
 b_1(Nw_l)\cdot i+
b_2(Nw_l)\cdot j=1,\ \ \ [Nw_1\wedge Nw_{s}]=b_2(Nw_1)\cdot i -b_2
(Nw_s)\cdot i.
\end{equation}
We also have, since $v_r\cdot w_1=0$,
\begin{equation}\label{spbasebisbisbis}
[Nv_r\wedge Nw_{1}]=b_2 (Nv_r)\cdot i - a_2 (Nw_1)\cdot i.
\end{equation}
We obtain, collecting (\ref{basebisbis}), (\ref{basebisbisbis}),
(\ref{spbasebisbisbis}) and (\ref{somme}),
$$
 N_{+}\left({\cal L}_N\right)-N_{-}\left({\cal L}_N\right)=a_2(Nv_1)\cdot i-b_2
(Nw_s)\cdot i + (b_2-a_2)((Nv_r)\cdot i + (Nw_1)\cdot i).
$$
From, the last equality we deduce the following,
\medskip

\noindent
$\bullet$ If $a_2=b_2$ i.e. $a_2b_2=1$, then $
 N_{-}\left({\cal L}_N\right)-N_{+}\left({\cal L}_N\right)=-a_2(Nv_1)\cdot i+b_2
(Nw_s)\cdot i.
$
\\
$\bullet$ If $a_2b_2=-1$  then since $v_r\cdot w_1=0$,
$(Nv_r)\cdot i + (Nw_1)\cdot i\in\{a_1,b_1\}$ and $$
 N_{-}\left({\cal L}_N\right)-N_{+}\left({\cal L}_N\right)+ a_2(Nv_1)\cdot i-b_2
(Nw_s)\cdot i $$ is either $2a_1a_2$ or $2b_1a_2$. \ \ $\Box$
\medskip

\noindent
The following corollary evaluates $N_{-}\left({\cal
  L}_N\right)-N_{+}\left({\cal
  L}_N\right)$ for a path ${\cal
L}_N$  behaving like a polygonal line. It will be very useful for
the control of $L_N^{\sigma,\gamma}(s,r,\alpha_1,\alpha_2)$.
\begin{coro}\label{spbb}
Let $s_0$, $s_1$ and $s_2$ be three points in $\BBr^2$. Let
$\theta_0$ (resp. $\theta_1$) be the oriented angle between the
half horizontal axis $[0,+\infty[$ and the segment $[s_1,s_0[$
(respectively $[s_1,s_2[$). Let ${\cal
L}_N=(v_1,\ldots,v_r,w_1,\ldots,w_{s})$ be a path on $\BBz_N^2$.
Suppose that the family $(v_1,\ldots,v_r)$ (respectively
$(w_1,\ldots,w_{s})$) forms
 a monotone path on $\BBz_N^2$ and that $v_r\cdot
w_1=0$. Suppose moreover that $(v_1,\ldots,v_r)$ and $[s_0,s_1]$
(respectively $(w_1,\ldots,w_{s})$ and $[s_1,s_2]$) are either
both non-increasing or either both non-decreasing. Then
   \begin{equation}\label{bbb}
  N_{-}\left({\cal
  L}_N\right)-N_{+}\left({\cal
  L}_N\right)=\sg(\sin\theta_0)(Nv_1)\cdot i +
  \sg(\sin\theta_1)(Nw_s)\cdot i + f(\theta_1,\theta_0),
  \end{equation}
where $$ f(\theta_1,\theta_0)= \left\{
\begin{array}{rl}
2\sg(\theta_1-\theta_0)\ \ & {{if}}\ \
\sin\theta_0\sin\theta_1>0,\ \ \cos\theta_0\cos\theta_1>0, \\
2\sg(\tan\theta_0)\ \ & {{if}}\ \
\sin\theta_0\sin\theta_1>0,\ \ \cos\theta_0\cos\theta_1<0, \\
0\ \  & \ \ {{otherwise.}}
\end{array}
\right. $$
\end{coro}

\noindent We illustrate the conclusion of the previous corollary
with the help of the following pictures.

\begin{figure}[hbt]
\vskip -0.5truecm \vbox{ \centerline{ \psset{unit=0.9cm}
\pspicture(0,-2)(10,8)
\psline[linewidth=1pt,linestyle=dashed](3.75,3.3)(4,3.3)(4.25,3.3)(5,3.3)
(6,3.3) (7,3.3)(8,3.3)(9,3.3)
\psline[linewidth=1pt,linestyle=dashed](3.75,3.3)(3.75,3.9)(3.75,4)
(3.75,5)(3.75,6)(3.75,7)
\rput(3.5,3.1){$s_1$} \psline[linewidth=1pt,showpoints=true]
(7.25,5.25)(3.75,3.3)(4.25,6.4) \rput(7.2,5){$s_0$}
\rput(4.2,6.6){$s_2$}
\psarc[arcsepA=2pt,arcsepB=2pt]{->}(3.75,3.3){1}{0}{26}
\psarc[arcsepA=2pt,arcsepB=2pt]{->}(3.75,3.3){0.5}{0}{79}
\rput(4.5,4){$\theta_1$} \rput(5,3.5){$\theta_0$}
\psline[linewidth=1.5pt]{->}(7.5,6)(6.75,6) \rput(7,6.25){ $v_1$}
\psline[linewidth=1.5pt]{->}(6.75,6)(6.75,5.25)
\psline[linewidth=1.5pt]{->}(6.75,5.25)(6.75,4.5)
\psline[linewidth=1.5pt]{->}(6.75,4.5)(6,4.5)
\psline[linewidth=1.5pt]{->}(6,4.5)(5.25,4.5)
\psline[linewidth=1.5pt]{->}(5.25,4.5)(4.5,4.5)
\psline[linewidth=1.5pt]{->}(4.5,4.5)(4.5,5.25)
\psline[linewidth=1.5pt]{->}(4.5,5.25)(5.25,5.25)
\psline[linewidth=1.5pt]{->}(5.25,5.25)(5.25,6)
\psline[linewidth=1.5pt]{->}(5.25,6)(5.25,6.75)
\rput(4.75,6.4){{$w_s$}} \rput(4.25,5){{$w_1$}} \rput(5,4.25){
$v_r$} \psgrid[subgriddiv=0,
unit=0.75,griddots=10,gridlabels=0](-1,4)(12,10)
\endpspicture
} \vskip -4truecm \centerline{${\cal L}_N$ is the circuit
$(v_1,\ldots,v_r,w_1,\ldots,w_s)$. } \centerline{ Here,
$f(\theta_1,\theta_0)=2\sg(\theta_1-\theta_0)=2$.} }
\end{figure}

\newpage

\begin{figure}[hbt]
\vskip-0.5truecm \vbox{ \centerline{ \psset{unit=0.9cm}
\pspicture(0,-1)(15,12) \psline[linewidth=0.5pt,showpoints=true]
(8.6,7.1)(6.3, 10.2)(4.2,7.7)
\psline[linewidth=1.5pt,showpoints=true](8.6,7.1)(6.3,10.2)(4.2,7.7)
\psline[linewidth=1.5pt]{->}(9,6.75)(9,7.5)
\psline[linewidth=1.5pt]{->}(9,7.5)(9,8.25)
\psline[linewidth=1.5pt]{->}(9,8.25)(8.25,8.25)
\psline[linewidth=1.5pt]{->}(8.25,8.25)(7.5,8.25)
\psline[linewidth=1.5pt]{->}(7.5,8.25)(6.75,8.25)
\psline[linewidth=1.5pt]{->}(6.75,8.25)(6.75,9)
\psline[linewidth=1.5pt]{->}(6.75,9)(6.75,9.75)
\psline[linewidth=1.5pt]{->}(6.75,9.75)(6.75,10.5)
\psline[linewidth=1.5pt]{->}(6.75,10.5)(6,10.5)
\psline[linewidth=1.5pt]{->}(6,10.5)(5.25,10.5)
\psline[linewidth=1.5pt]{->}(5.25,10.5)(5.25,9.75)
\psline[linewidth=1.5pt]{->}(5.25,9.75)(5.25,9)
\psline[linewidth=1.5pt]{->}(5.25,9)(5.25,8.25)
\psline[linewidth=1.5pt]{->}(5.25,8.25)(4.5,8.25)
\psline[linewidth=1.5pt]{->}(4.5,8.25)(3.75,8.25)
\psline[linewidth=1.5pt]{->}(3.75,8.25)(3.75,7.5)
\psgrid[subgriddiv=0,unit=0.75,griddots=10,gridlabels=0](3,8)(15.25,16)
\rput(9.25,7.25){$v_1$} \rput(9.25,7.75){$v_2$}
\rput(5,10){$ w_1$} \rput(5.5,10.75){$ v_r$}
\rput(3.3,7.75){$w_s$} \rput(4.2,7.4){$s_2$} \rput(8.6,6.8){$s_0$}
\rput(6.3, 9.9){$s_1$}
\psline[linewidth=0.25pt,linestyle=dashed](3.25,10.2)(4.25,10.2)(5.25,10.2)(6.25,10.2)
(7.25,10.2)(8.25,10.2)(9.25,10.2)
\psarc[arcsepA=2pt,arcsepB=2pt]{->}(6.3,10.2){1}{0}{311}
\psarc[arcsepA=2pt,arcsepB=2pt]{->}(6.3,10.2){1.75}{0}{230}
\rput(8.25,11.25){$\theta_1$} \rput(6.3,8.75){$\theta_0$}
\endpspicture
} \vskip -6truecm \centerline{${\cal L}_N$ is the circuit
$(v_1,\ldots,v_r,w_1,\ldots,w_s)$. Here  } \centerline{
$f(\theta_1,\theta_0)=2\sg(\tan\theta_0)=-2$.} }
\end{figure}

\begin{figure}[hbt]
\vskip 0.9truecm \vbox{ \centerline{ \psset{unit=0.9cm}
\pspicture(0,-1)(15,12)
\psline[linewidth=1.5pt]{->}(8.25,10.5)(7.5,10.5)
\rput(8,10.7){$v_1$}
\psline[linewidth=1.5pt]{->}(7.5,10.5)(7.5,9.75)
\rput(7.2,10){$v_2$} \psline[linewidth=1.5pt]{->}(7.5,9.75)(7.5,9)
\psline[linewidth=1.5pt]{->}(7.5,9)(6,9)
\psline[linewidth=1.5pt]{->}(6,9)(5.25,9)
\psline[linewidth=1.5pt]{->}(5.25,9)(5.25,8.25)
\psline[linewidth=1.5pt]{->}(5.25,8.25)(5.25,7.5)
\psline[linewidth=1.5pt]{->}(5.25,7.5)(5.25,6.75)
\psline[linewidth=1.5pt]{->}(5.25,6.75)(6,6.75)
\psline[linewidth=1.5pt]{->}(6,6.75)(6.75,6.75)
\psline[linewidth=1.5pt]{->}(6.75,6.75)(7.5,6.75)
\psline[linewidth=1.5pt]{->}(7.5,6.75)(7.5,6)
\psline[linewidth=1.5pt]{->}(7.5,6)(7.5,5.25)
\psline[linewidth=1.5pt]{->}(7.5,5.25)(8.25,5.25)
\psline[linewidth=1.5pt]{->}(8.25,5.25)(8.25,4.5)
\psline[linewidth=1.5pt]{->}(8.25,4.5)(9,4.5)
\psline[linewidth=1.5pt]{->}(9,4.5)(9.75,4.5)
\rput(9.4,4.2){$w_s$} \rput(9.4,5.25){$s_2$} \psgrid[subgriddiv=0,
unit=0.75, griddots=10,gridlabels=0](3.75,3.75)(17,17)
\psline[linewidth=0.5pt,showpoints=true]
(8.25,10)(4.8,7.9)(9.4,4.9) \rput(8.25,9.5){$s_0$}
\psline[linewidth=0.25pt,linestyle=dashed](4.8,3.5)(4.8,5.5)(4.8,7.2)(4.8,8.1)(4.8,9.4)(4.8,10.3)(4.8,11)
\psline[linewidth=0.25pt,linestyle=dashed](3.4,7.9)(6.3,7.9)(10.2,7.9)
\psarc[arcsepA=2pt,arcsepB=2pt]{->}(4.8,7.9){0.9}{0}{38}
\rput(5.99,8){$\theta_0$}
\psarc[arcsepA=2pt,arcsepB=2pt]{->}(4.8,7.9){0.7}{0}{335}
\rput(3.75,8.2){$\theta_1$}
\endpspicture
} \vskip -3.25truecm \centerline{ In this picture,
$f(\theta_1,\theta_0)=0$.} }
\end{figure}

\noindent{\bf{Proof of corollary \ref{spbb}.}} We first check that
${\mbox{for any}}\ \ 1\leq l\leq r,\ \ (Nv_l)\cdot i\in\{0,-
\sg(\cos\theta_0)\},$ and $(Nv_l)\cdot j\in\{0,-
\sg(\sin\theta_0)\}.$ In the same way, we have $ 1\leq l\leq s$,
$$
 (Nw_l)\cdot i\in\{0,\sg(\cos\theta_1)\},\ \ \ {\mbox{and}}\ \ \ (Nw_l)\cdot j\in\{0,
\sg(\sin\theta_1)\}.
$$
Lemma \ref{basebis} gives then
$$
  N_{-}\left({\cal
  L}_N\right)-N_{+}\left({\cal
  L}_N\right)=\sg(\sin\theta_0)(Nv_1)\cdot i +
  \sg(\sin\theta_1)(Nw_s)\cdot i + f(\theta_1,\theta_0),
$$
\newpage
\noindent where \\
$\bullet$ If $a_2b_2= -\sg(\sin\theta_0)\sg(\sin\theta_1)>0$, then
$f(\theta_1,\theta_0)=0$. \\
$\bullet$ If $\sg(\sin\theta_0)\sg(\sin\theta_1)>0$ and $\sg(\cos\theta_0)\sg(\cos\theta_1)<0$,  then $f(\theta_1,\theta_0)=2\sg(\tan\theta_0)$. \\
In fact this case corresponds to $a_2b_2=-1$,
$a_1=-\sg(\cos\theta_0)=\sg(\cos\theta_1)=b_1$.

Now we have to discuss the case
$\sg(\sin\theta_0)\sg(\sin\theta_1)>0$ and
$\sg(\cos\theta_0)\sg(\cos\theta_1)>0$ i.e. when
 $ a_2b_2=-1$ and $a_1b_1=-1$. We distinguish all the cases on the
 possible values of $(a_1,a_2)$ and we deduce the following:
 $
 (Nv_r)\cdot i=0
 $
 if and only if $(a_1a_2<0\,{\mbox{and}}\, \theta_0<\theta_1)$ or
$(a_1a_2>0\,{\mbox{and}}\, \theta_0>\theta_1)$. So $
 (Nv_r)\cdot i=0
 $
 if and only if $a_1a_2\sg(\theta_0-\theta_1)>0$. We apply again lemma
 \ref{basebis} and we deduce that in this last case
 $f(\theta_1,\theta_0)=2\sg(\theta_1-\theta_0)$.\,\, \, $\Box$
\medskip

\noindent
We have now all the ingredients in order to evaluate
$\nu_N\left(N_{-}\left({\cal
  L}^{\sigma}_N\right)-N_{+}\left({\cal
  L}^{\sigma}_N\right)\right)$ for the random path ${\cal
  L}^{\sigma}_N$ as defined by (\ref{vbis}). The curve $\gamma$ is of class ${\cal C}_1$, hence
  for $r$ small enough, the part of $\gamma$ situated between  $x_0$ and $s$
  (resp. between $s$ and $x_1$) is either nondecreasing or nonincreasing.
  We conclude from the assumptions of
  theorem \ref{T1} that, for $N$ large enough and with probability one,
  the
  random path ${\cal
  L}^{\sigma}_N$ respects  the behavior of the
  curve $\gamma$, thus ${\cal
  L}^{\sigma}_N$ is either  monotone or it is constructed by concatenating
two
  monotone paths, say ${\cal
  L}^{\sigma}_N= ({\cal
  L}^{\sigma}_{1,N}, {\cal
  L}^{\sigma}_{2,N})$. These monotone paths are such that, noting by $s'$ the point of $T_{x_0}\gamma\cap T_{x_1}\gamma$, ${\cal
  L}^{\sigma}_{1,N}$ and $[x_0,s']$ (resp. ${\cal
  L}^{\sigma}_{2,N}$ and $[s',x_1]$) are either both nondecreasing
  or both nonincreasing. Corollary \ref{spbb} applies and gives,
  for $N$ large enough,
$$
\nu_N\left( N_{-}\left({\cal
  L}^{\sigma}_N\right)-N_{+}\left({\cal
  L}^{\sigma}_N\right)\right)=\sg(\sin\theta_0)\,\nu_N(e_N^{1}(\alpha_1)\cdot i) +
  \sg(\sin\theta_1)\,\nu_N(e_N^{2}(\alpha_2)\cdot i )+ f(\theta_1,\theta_0),
$$
the function $f(\theta_1,\theta_0)$ is defined in corollary
\ref{spbb}, the angles $\theta_1,\theta_0$ are those defined by
theorem \ref{T1}, the random vectors ${e}_N^1(\alpha_1)$,
${e}_N^2(\alpha_2)$ are the two unit vectors as defined by
(\ref{v}). We then deduce from proposition \ref{a1}  that there
exists $N_0$ depending only on $\gamma$ such that, for any $N\geq
N_0$, we have,
 $$\nu_N\left(L_N^{\sigma,\gamma}(s,r,\alpha_1,\alpha_2)\right)= \frac{1}{2}\sg(\sin\theta_0)\,\nu_N\left(e_N^{1}(\alpha_1)\cdot
i\right)+\frac{1}{2}\sg(\sin\theta_1)\,\nu_N\left(e_N^{2}(\alpha_2)\cdot
i\right) + \frac{1}{2}f(\theta_1,\theta_0).$$

\subsection{Evaluation of
$\int_0^{\delta}\int_0^{\delta}\nu_N\left(L_N^{\sigma,\gamma}(s,r,\alpha_1,\alpha_2)\right)\,d\alpha_1\,d\alpha_2$
}\label{SEC} By the previous formula, in order to evaluate the
quantity $$ \int_0^{\delta}\int_0^{\delta}
\nu_N\big(L_N^{\sigma,\gamma}(s,r,\alpha_1,\alpha_2)\big)
\,d\alpha_1\,d\alpha_2,
$$ for $\delta$ and $r$ small enough, it suffices to evaluate the
terms
$$\nu_N\left(\int_{0}^{\delta}\,e_N^{1}(\alpha)\cdot
i\,d\alpha\right)\,,\qquad
\nu_N\left(\int_{0}^{\delta}\,e_N^{2}(\alpha)\cdot
i\,d\alpha\right)\,.$$
We begin by the first quantity, for this we
need some further notations.
\newpage

\noindent {\bf{Notation.}} For a vector  $v$ drawn on the grid
$\BBz^2_N$, we denote by $R(v)$ the union of the two boxes of the
family $(\Lambda_{x/N})_{x\in\BBz^2}$ having $v$ as an edge
vector.

\vspace{-2mm}
\begin{figure}[hbt]
\vbox{ \centerline{ \psset{unit=0.7cm} \pspicture(2.5,2.5)(4,3.5)
\pspolygon[fillstyle=solid,fillcolor=lightgray,
linewidth=0pt](0.75,0.75)(2.25,0.75)(2.25,1.5)(0.75,1.5)
\pspolygon[fillstyle=solid,fillcolor=lightgray,
linewidth=0pt](3,0.75)(3.75,0.75)(3.75,2.25)(3,2.25)
\psgrid[subgriddiv=0,gridlabels=0,griddots=10,unit=0.75](6,4)
\psline[linewidth=1.5pt]{->}(1.5,0.75)(1.5,1.5)
\psline[linewidth=1.5pt]{->}(3.75,1.5)(3,1.5) \rput(1.25,1){$v$}
\rput(3.5,1.25){$w$}
\endpspicture
} \vskip 2cm \centerline{The two blocks $R(v)$ and $R(w)$.} }
\end{figure}

{\hspace{-7mm}} Let ${\cal L}^{\sigma}_{N}=(v_1,\ldots,v_{r})$ be
the oriented path as defined by (\ref{vbis}). Let ${\cal
L}^{\sigma}_{1,N}=(v_1,\ldots,v_s)$ be the subgraph of ${\cal
L}^{\sigma}_{N}$ included in  $\partial {\cal A}_{\sigma}^N \cap
(B(x_0,\delta))_N$ such that the vector $v_s$ is the entering
vector in $(B(s,r))_N$.

To each vector $v_{l}$ ($1\leq l\leq s$), we associate the block
$R_{s-l+1}:=R(v_{l})$. These blocks $(R_l)_{1\leq l\leq s}$ are
enumerated according to their distances to $x_0$, $R_1$ being the
block containing $v_s$. Let $(a_l)_{1\leq l\leq s}$ be the
sequence of vertices such that $$ d_l:=d(x_0,R_l)=|a_l-x_0|,
$$ then this sequence of vertices $(a_l)_{1\leq l\leq s}$ is $L^1$
connected and the vector $a_{l}a_{l+1}$ is either vertical or
horizontal. Finally, let ${\cal H}_{N}$ be the set of indices
$l\in\{1,\ldots,s\}$ for which $v_l$ is horizontal.

\vspace{1cm}
\begin{figure}[hbt]
\vbox{ \centerline{ \psset{unit=0.9cm} \pspicture(0,-1)(9,7)
\pspolygon[fillstyle=solid,fillcolor=lightgray,
linewidth=0pt](6,3)(4.5,3)(4.5,3.75)(6,3.75)
\pspolygon[fillstyle=solid,fillcolor=gray,
linewidth=0pt](6.75,3)(7.5,3)(7.5,4.5)(6.75,4.5)
\psline[linewidth=1.5pt]{->}(10.5,4.5)(9.75,4.5)
\pspolygon[fillstyle=solid,fillcolor=gray,
linewidth=0pt](3,2.25)(3.75,2.25)(3.75,3.75)(3,3.75)
\rput(4.5,2){\bf $a_3$} \rput(4.3,3.2){\bf $a_4$}
\psline[linewidth=0.75pt, showpoints=true]{->}(4.5,2.25)(4.5,3)
\psline[linewidth=1.5pt]{->}(11.25,4.5)(10.5,4.5)
\psline[linewidth=1.5pt]{->}(9.75,4.5)(9,4.5)
\psline[linewidth=1.5pt]{->}(9,4.5)(9,3.75)
\psline[linewidth=0.75pt,showpoints=true,linestyle=dashed](3.2,2.3)(6.75,3)
\rput(6.75,2.7){$a_i$} \rput(7.5,2.7){$a_{i+1}$}
\psline[linewidth=1.5pt,showpoints=true]{->}(6.75,3)(7.5,3)
\psline[linewidth=1.5pt]{->}(9,3.75)(8.25,3.75)
\psline[linewidth=1.5pt]{->}(8.25,3.75)(7.5,3.75)
\psline[linewidth=1.5pt]{->}(7.5,3.75)(6.75,3.75)
\psline[linewidth=1.5pt]{->}(6.75,3.75)(6,3.75)
\psline[linewidth=1.5pt]{->}(6,3.75)(5.25,3.75)
\psline[linewidth=1.5pt]{->}(5.25,3.75)(5.25,3)
\psline[linewidth=1.5pt]{->}(5.25,3)(4.5,3)
\psline[linewidth=1.5pt]{->}(4.5,3)(3.75,3)
\psline[linewidth=1.5pt]{->}(3.75,3)(3,3)
\psline[linewidth=2.25pt]{->}(6.75,5.25)(6.9,4.5)
\rput(6.7,5.5){$R_i$}
\psline[linewidth=2.25pt]{->}(4.5,4.5)(4.75,3.75)
\rput(4.6,5){$R_4$}
\psline[linewidth=2.25pt]{->}(2.25,3.75)(3,3.5)
\rput(1.75,4){$R_1$}
\psgrid[subgriddiv=0,unit=0.75,griddots=10,gridlabels=0](16,9)
\psline[linewidth=1.5pt]
(10.9,2.8)(3.2,2.3)(1,2.1) \rput(3.2,2){$x_0$}
\rput(1,1.5){$T_{x_0}\gamma$}
\psline[linewidth=0.75pt]{->}(3.2,2.3)(10.8,1.5)
\rput(6.75,1.59){$\delta$}
\psarc[arcsepA=2pt,arcsepB=2pt](3.2,2.3){7.6}{-10}{30}
\endpspicture
} \centerline{$d_l=d(x_0,R_l)=|x_0-a_l|$.} \centerline{For $N$
large enough, the vector $a_la_{l+1}$ is either horizontal or
vertical, and $|a_l-a_{l+1}|=\frac{\textstyle 1}{\textstyle N}$.}
}
\end{figure}
{\hspace{-7mm}}With probability one, the path ${\cal
L}^{\sigma}_{1,N}$ is monotone and behaves, on a neighborhood of
$x_0$, as $T_{x_0}\gamma$. This fact ensures that, with probability
one, $e_N^{1}(\alpha)\cdot i\in\{0,-\sg(\cos\theta_0)\}$. Now, by
construction $e_N^{1}(\alpha)\cdot i=-\sg(\cos\theta_0)$ if and
only if there exists $l\in {\cal H}_{N}$ such that $\alpha\in
]d_l,d_{l+1}]$ (such an index is necessarily unique).
With probability one,
\begin{equation}\label{l1}
\left|\int_{0}^{\delta}e_N^{1}(\alpha)\cdot{i}\,d\alpha
+\sg(\cos\theta_0)\sum_{l\in {\cal H}_N(x_0,\delta)}
\left(d_{l+1}-d_l\right)\right|\leq \frac{2}{N}\,,
\end{equation}
where
${\cal H}_N(x_{0},\delta)$
is the set
 of all the horizontal edges of $\partial
 {\cal A}^{\sigma}_N$ included in $(B(x_0,\delta))_N
 \setminus B(s,r)$.
 In order to evaluate $d_{l+1}-d_l$, we need the following
lemma.
\begin{lemma}\label{lem2}
Let ${u}$ and ${v}$ be two vectors such that $\|{u}\|\leq
\|{v}\|$. Then $$\|{u}+{v}\|- \|{v}\|=\frac{({u}+{v})\cdot{u}}
{\|{u}+{v}\|} -
\frac{\|{u}\|^2}{\|{v}\|}\frac{\sin^2\theta}{1+\sqrt{1-\frac{\|{u}\|^2}
{\|{v}\|^2}\sin^2\theta}},$$ where $\theta$ is the angle between
${u}$ and ${u}+{v}$.
\end{lemma}
{\bf{Proof of lemma \ref{lem2}.}} Let ${u}$, ${v}$ and $\theta$ be
as defined in lemma \ref{lem2}.

\begin{figure}[hbt]
\vbox{ \centerline{ \psset{unit=0.9cm} \pspicture(0,-1)(9,4)
\psline[linewidth=1.5pt]{->}(2,0) \rput(1,-0.4){u}
\psline[linewidth=1.5pt]{->}(5,4) \rput(2.5,3){u+v}
\psline[linewidth=1.5pt]{->}(2.1,0)(5,3.8) \rput(3,1.6){v}
\psline[linestyle=dotted](2,0)(5,0)
\psline[linestyle=dotted](5,0)(5,4)
\psline[linewidth=1.5pt]{<->}(5.8,0)(5.8,4) \rput(6,2){H}
\psline[linewidth=1.5pt]{<->}(0,-0.8)(5,-0.8) \rput(2.5,-1.1){L}
\psarc[arcsepA=2pt,arcsepB=2pt]{->}(0,0){1}{0}{40}
\rput(1.2,0.5){$\theta$}
\endpspicture
} }
\end{figure}
\noindent We have
\begin{eqnarray*}
\|{u}+{v}\|^2 & = & L^2 + H^2
\\ & = &
\cos^2\theta\|{u}+{v}\|^2 +\|{v}\|^2- \left(\cos\theta\|{u}+{v}\|-
\|{u}\|\right)^2 \\ & = & \|{v}\|^2+ 2\cos\theta \|{u}\|\times
\|{u}+{v}\| -\|{u}\|^2.
\end{eqnarray*}
The quantity $\|{u}+{v}\|$ is then a positive solution of an
algebraic equation of degree two. We deduce from $\|{u}\|\leq
\|{v}\|$, that
$$\|{u}+{v}\|= \|{u}\|\cos\theta+
\sqrt{\|{v}\|^2-\sin^2\theta\|{u}\|^2}. $$ Hence
\begin{eqnarray*}
\|{u}+{v}\|-\|{v}\| & = & \|{u}\|\cos\theta +\|{v}\|\left(\sqrt{1-
\sin^2\theta\frac{\|{u}\|^2}{\|{v}\|^2}}-1\right)
\\
& = & \|{u}\|\cos\theta -
\frac{\|{u}\|^2}{\|{v}\|}\frac{\sin^2\theta}{1+\sqrt{1-\frac{\|{u}\|^2}
{\|{v}\|^2}\sin^2\theta}}.
\end{eqnarray*}
The last equality together with the fact that $\|{u}\|\cos\theta =
\frac{\textstyle({u}+{v})\cdot{u}} {\textstyle\|{u}+{v}\|}$ proves
lemma \ref{lem2}. $\Box$
\medskip

\noindent We continue the proofs of theorems \ref{T1} and
\ref{T2}. We apply lemma \ref{lem2} with ${u}={a_la_{l+1}}$,
${v}={x_{0}a_l}$ and we get
$$ \frac{({u}+{v})\cdot{u}} {\|{u}+{v}\|} =
\frac{(x_{0}a_{l+1})\cdot(a_la_{l+1})} {|{x_{0}-a_{l+1}}|}. $$
Moreover, we deduce from lemma \ref{lem2},
\begin{equation}\label{l2}
\left|d_{l+1}-d_l- \frac{(x_{0}a_{l+1})\cdot (a_la_{l+1})}
{|{x_{0}-a_{l+1}}|}\right|\leq
\min\Big(\frac{1}{N^2|{x_{0}-a_l}|},{2\over N}\Big).
\end{equation}
\marginpar{j'ai change ici parce que je ne comprenais pas ce que tu
voulais dire avec le $\epsilon_N$?}
We first evaluate the sum over $l\in {\cal H}_{N}(x_0,\delta)$ of
the right hand side of the last inequality.
Let $\phi(l)$ be the cardinality of the set ${\cal
  H}_N(x_0,\delta)\cap\{1,\ldots,l\}$.
For $l\in {\cal H}_{N}(x_0,\delta)$, we have $|(x_{0}a_l)\cdot
i|\geq N^{-1}(\phi(l)-1)$ , whence
  \begin{eqnarray}\label{l4}
\sum_{l\in{\cal H}_{N}(x_0,\delta)}
\min\Big(\frac{1}{N^2|{x_{0}-a_l}|},{2\over N}\Big)
 & \leq & \frac{1}{N}\sum_{l\in {\cal
H}_{N}}
\min\Big(\frac{1}{\phi(l)-1},2\Big) \, \leq \, \frac{2+
\ln|{\cal H}_{N}|}{N}\,.
  \end{eqnarray}
With probability one, we have
$$\sg(\cos\theta_0)(a_la_{l+1})\cdot i +
\sg(\sin\theta_0)(a_la_{l+1})\cdot j=\frac {1}{N}\,,$$
whence
\begin{eqnarray}\label{l31}
{\lefteqn{\sum_{l\in {\cal
H}_{N}(x_0,\delta)}\frac{(x_{0}a_{l+1})\cdot (a_la_{l+1})}
{|{x_{0}-a_{l+1}}|}}}
\\
&& = \frac{\sg(\cos\theta_0)}{N}\sum_{l\in{\cal
H}_{1,N}(x_0,\delta)} \frac{{(x_0a_{l+1})}\cdot{{i}}}
{|{x_0-a_{l+1}}|}+ \frac{\sg(\sin\theta_0)}{N}\sum_{l\in{\cal
H}_{2,N}(x_0,\delta)} \frac{({x_0a_{l+1})}\cdot{{j}}}
{|{x_0-a_{l+1}}|}\nonumber,
  \end{eqnarray}
where
\begin{eqnarray}\label{l41}
{\cal H}_{1,N}(x_0,\delta)=\big\{\,l\in{\cal H}_{N}(x_0,\delta):
  ({a_{l}a_{l+1}})\cdot{{j}}=0\,\big\}\,,\nonumber\\
{\cal H}_{2,N}(x_0,\delta)=
  \big\{\,l\in{\cal H}_{N}(x_0,\delta):
  ({a_{l}a_{l+1}})\cdot{{i}}=0\,\big\}.
  \end{eqnarray}
  \vskip 0.5truecm
\noindent We now distinguish the case of the polygons and the case
of the Jordan curves.
\subsection{End of the proof for polygons (theorem~\ref{T1}).}

\begin{lemma}\label{lem5sp}
For $\delta$ small enough,
we have
$$\displaylines{
\lim_{N\rightarrow \infty} \bigg|\nu_N\bigg( \frac{1}{N}
\sum_{l\in {\cal H}_{1,N}(x_0,\delta)}
\frac{{(x_0a_{l+1})}\cdot{{i}}} {|{x_0-a_{l+1}}|}\bigg) \,-\,
\nu_N\bigg(\frac{|{\cal H}_{1,N}(x_0,\delta)|}{N}\bigg)
\cos\theta_0\bigg|\,=\,0\,,\cr \lim_{N\rightarrow \infty}
\bigg|\nu_N\bigg( \frac{1}{N} \sum_{l\in {\cal
H}_{2,N}(x_0,\delta)} \frac{{(x_0a_{l+1})}\cdot{{j}}}
{|{x_0-a_{l+1}}|}\bigg) \,-\, \nu_N\bigg(\frac{|{\cal
H}_{2,N}(x_0,\delta)|}{N}\bigg) \sin\theta_0\bigg|\,=\,0\,. }$$
\end{lemma}
\newpage
\noindent{\bf{Proof of lemma \ref{lem5sp}.}} We only prove the
first limit since the argument for the second limit is similar.
Let $u$ be a unit vector tangent to~$\gamma$ at $x_0$ and
let $v$ be such that $(u,v)$ is a direct basis.
For $\varepsilon>0$, let ${\cal R}({\varepsilon})$ be the strip
of width $2\varepsilon$ centered on the tangent line $T_{x_0}\gamma$,
i.e.,
$${\cal R}(\varepsilon)\,=\,
\big\{\,x\in\BBr^2:|x_0x\cdot v|\leq\varepsilon\,\}\,.$$
The condition (\ref{haus}) implies
that for $\delta$ small enough,
\marginpar{ceci n'est vrai que pour une ligne: pour une courbe
quelconque, on est dans un tube autour de la courbe, pas de la tangente}
$$
\forall\varepsilon>0\qquad \lim_{N\rightarrow
+\infty}\nu_N\left(\partial {\cal A}_N^{\sigma}\cap
B(x_0,\delta)\subset {\cal R}({\varepsilon})\right)=1.
$$

\begin{figure}[hbt]
\vskip 8truecm \vbox{ \centerline{ \psset{unit=1.5cm}
\pspicture(2.5,2.5)(5,-0.5)
\psline[linewidth=0.5pt,linecolor=red](6,7.5)(3,4.5)
\rput(3,4.5){$\bullet$} \rput(4.5,6){$\bullet$}
\psline[linewidth=2pt,linecolor=red](5.5,6.5)(5,6.5)(4.75,6.5)(4.5,6.5)(4.25,6.5)(4.25,6.25)(4.25,6)(3.9,6)(3.9,5.27)
\psline[linewidth=1pt,linecolor=blue](6,8.25)(5.25,7.5)(3.3,5.6)
\psline[linewidth=1pt,linecolor=blue](6.75,7.5)(4.5,5.25)(3.75,
4.5) \rput(6,7.7){$T_{x_0}\gamma$} \rput(2.7,4.65){$s$}
\rput(4.5,6.2){$x_0$}
\psline[linewidth=0.85pt]{->}(4.5,6)(5.25,5.5)
\rput(5,6){$\delta$}
\psarc[arcsepA=2pt,arcsepB=2pt](4.5,6){0.9}{0}{360}
\psline[linewidth=2pt,linecolor=blue]{<->}(5.25,7.5)(6,6.75)
\rput(6.12,7.2){$2\epsilon$}
\psline[linewidth=0.5pt,linestyle=dashed,
linecolor=red](3,4.5)(3.75, 4.5)(4,4.5)(4.25,4.5)
\psarc[arcsepA=2pt,arcsepB=2pt]{->}(3,4.5){0.7}{0}{45}
\rput(3.5,4.75){$\theta_0$}
\endpspicture
} } \vskip -7truecm {\centerline{For $\delta$ small enough,
$\lim_{N\rightarrow +\infty}\nu_N\left(\partial {\cal
A}_N^{\sigma}\cap B(x_0,\delta)\subset {\cal
R}({\varepsilon})\right)=1$. }}
\end{figure}

\noindent Let $\delta>0$ be small enough so that the above limit
holds. Let $\delta_0$, $\varepsilon$ such that
$0<\varepsilon<\delta_0<\delta$ and let $x\in{\cal
R}(\varepsilon)\setminus B(x_0,\delta_0)$. We have
$$\displaylines{x_0x\cdot i\,=\,
(x_0x\cdot u)\cos\theta_0-
(x_0x\cdot v)\sin\theta_0\,,\cr
|x_0-x|^2\,=\,
(x_0x\cdot u)^2+
(x_0x\cdot v)^2\,,}$$
whence
$$
\frac{x_0x\cdot i}{|x_0-x|}\,=\,
\bigg(1-
\frac{(x_0x\cdot v)^2}{|x_0-x|^2}
\bigg)^{1/2}
\cos\theta_0-
\frac{(x_0x\cdot v)}{|x_0-x|}\sin\theta_0\,$$
and
$$\bigg|
\frac{x_0x\cdot i}{|x_0-x|}-\cos\theta_0\bigg|\,\leq\,
1-\sqrt{1-\varepsilon^2/\delta_0^2}+
\frac{\varepsilon}{\delta_0}\,\leq\,2
\frac{\varepsilon}{\delta_0}\,.$$ If the event $\{\,\partial
{{\cal A}}_N^{\sigma}\cap B(x_0,\delta)\subset {\cal
R}({\varepsilon})\,\}$ occurs, then for $l\in {\cal
H}_{N}(x_0,\delta)\setminus {\cal H}_{N}(x_0,\delta_0)$, we have
$a_{l+1}\in {\cal R}({\varepsilon})\setminus B(x_0,\delta_0)$, and
thus
$$
\limsup_{N\rightarrow \infty}\,\,\nu_N\left(\sup_{l\in {\cal
H}_{N}(x_0,\delta)\setminus {\cal H}_{N}(x_0,\delta_0)}
\left|\frac{({x_0a_{l+1})}\cdot{{i}}} {|{x_0-a_{l+1}}|}
-\cos\theta_0\right|\right)\,\leq\,
2\frac{\varepsilon}{\delta_0}\,.$$ Moreover, we have $\big|{\cal
H}_{N}(x_0,\delta_0)\big|\,\leq\,2N\delta_0$, whence, by splitting
the sum over ${\cal H}_{N}(x_0,\delta_0)$ and ${\cal
H}_{N}(x_0,\delta)\setminus {\cal H}_{N}(x_0,\delta_0)$, we obtain
$$\limsup_{N\rightarrow \infty}
\bigg|\nu_N\bigg( \frac{1}{N} \sum_{l\in {\cal
H}_{1,N}(x_0,\delta)} \frac{{(x_0a_{l+1})}\cdot{{i}}}
{|{x_0-a_{l+1}}|}\bigg) \,-\, \nu_N\bigg(\frac{|{\cal
H}_{1,N}(x_0,\delta)|}{N}\bigg) \cos\theta_0\bigg|\,\leq\,
4\frac{\delta\varepsilon}{\delta_0}+4\delta_0\,.$$
\newpage
\noindent We conclude by sending successively $\varepsilon$ to~$0$
and $\delta_0$ to~$0$.
 \ \ $\Box$
\medskip

\noindent We obtain, combining (\ref{l31}) and lemma \ref{lem5sp},
that for $\delta$ small enough,
$${\lim_{N\rightarrow \infty}
\bigg|\nu_N\left(\sum_{l\in {\cal
H}_{N}(x_0,\delta)}\kern-4pt\frac{(x_{0}a_{l+1})\cdot (a_la_{l+1})}
{|{x_{0}-a_{l+1}}|}\right)\kern150pt}
$$
\vspace{-0.6cm}
\begin{eqnarray}\label{e6} &&
-\left(\nu_N\left(\frac{|{\cal
H}_{1,N}(x_0,\delta)|}{N}\right)|\cos\theta_0|+
\nu_N\left(\frac{|{\cal
H}_{2,N}(x_0,\delta)|}{N}\right)|\sin\theta_0|\right)\bigg|
\,=\,0\,.
\end{eqnarray}
Our purpose now is to evaluate, for $N$ large enough, the
expectations over $\nu_N$ of $\frac{\textstyle |{\cal
H}_{N}(x_0,\delta)|}{\textstyle N}$, $\frac{\textstyle |{\cal
H}_{1,N}(x_0,\delta)|}{\textstyle N}$ and $\frac{\textstyle |{\cal
H}_{2,N}(x_0,\delta)|}{\textstyle N}$.  For this, we prove the
following lemma.
\begin{lemma}\label{lembis5sp}
 For
$\delta$ small enough, one has
\begin{equation}\label{c1bis}
\lim_{N\rightarrow +\infty}\nu_N\left(\frac{|{\cal
H}_N(x_{0},\delta)|}{N}\right) =\delta|\cos\theta_0| .
\end{equation}
\end{lemma}
{\bf{Proof of lemma \ref{lembis5sp}.}} We denote by $a$ and $x'$
the points of $\partial B(x_0,\delta)\setminus B(s,r)$ belonging
respectively to $\partial {\cal A}_N^{\sigma}$ and to
$T_{x_0}\gamma$. Let $b$ be the point of $\partial {\cal
A}_N^{\sigma}\cap \partial B(s,r)\setminus B(x_1,\delta)$. We
suppose without loss of generality that $(ba)\cdot i\geq 0$.

\begin{figure}[hbt]
\vskip 8truecm \vbox{ \centerline{ \psset{unit=1.5cm}
\pspicture(2.5,2.5)(5,-0.5) \rput(5.1,6.9){$a$}
\psline[linewidth=0.75pt,linecolor=red,showpoints=true](3,4.5)(4.99,3.81)(6,3.5)(6.75,4.5)(6,7.5)
\rput(4.75,3.75){$x_1$}
\psline[linewidth=1.2pt,linecolor=red](5.2,6.8)(5,6.8)(4.8,6.8)(4.8,6.2)(4.7,6.2)(4.7,6)(4.7,5.8)(4.4,5.8)(4.2,5.8)(4.2,5.6)(4,5.6)
\rput(5.4,6.7){$x'$} \rput(2.9,4.4){$s$} \rput(4.21,6){$x_0$}
\rput(4.66,5.8){$\bullet$} \rput(4.8,5.8){$b$}
\psline[linewidth=0.85pt]{->}(4.5,6)(4.5,6.9)
\rput(4.3,6.4){$\delta$}
\psline[linewidth=0.85pt]{->}(3,4.5)(1.5,6) \rput(2.5,5.3){$r$}
\psarc[arcsepA=2pt,arcsepB=2pt](4.5,6){0.9}{0}{362}
\psarc[arcsepA=2pt,arcsepB=2pt](3,4.5){2.1}{0}{362}
\psline[linewidth=0.5pt,linestyle=dashed,
linecolor=red](3,4.5)(3.75, 4.5)(4,4.5)(4.25,4.5)
\psarc[arcsepA=2pt,arcsepB=2pt]{->}(3,4.5){0.7}{0}{45}
\rput(3.5,4.75){$\theta_0$} \rput(4.9,6.8){$\bullet$}
\psline[linewidth=0.75pt,linecolor=red,showpoints=true](6,7.5)(5.12,6.6)(4.5,6)(3,4.5)
\endpspicture
} } \vskip -3truecm {\centerline{In this case, $\gamma$ is a
polygon. The proportion of the horizontal edges of $\partial {\cal
A}_N^{\sigma}$ which are}} \centerline{{ in
$B(x_0,\delta)\setminus B(s,r)$ is controlled by
$\delta|\cos\theta_0|$. }}
\end{figure}

\newpage
\noindent We have, by definition of ${\cal H}_N(x_{0},\delta)$,
$$
\left|\frac{|{\cal H}_N(x_{0},\delta)|}{N} -|ba\cdot i|\right|\leq
\frac{2}{N}.
$$
 We use the same notation as in the proof of lemma
\ref{lem5sp}. We have
$$\displaylines{\bigcap_{\varepsilon>0}
{\cal R}(\varepsilon)\cap\partial B(x_0,\delta)\setminus B(s,r)
\,=\, T_{x_0}\gamma\cap\partial B(x_0,\delta)\setminus B(s,r)\,
=\,\{x'\}\,,\cr \bigcap_{\varepsilon>0} {\cal R}(\varepsilon)\cap
\partial B(s,r)\setminus B(x_1,\delta)
\,=\,
T_{x_0}\gamma\cap
\partial B(s,r)\setminus B(x_1,\delta)\,
=\,\{x_0\}\,.}
$$
Let $\alpha>0$. By the above identities,
there exists $\varepsilon>0$ such that,
if $\{\,\partial {\cal A}_N^{\sigma}\cap B(x_0,\delta)\subset
{\cal R}({\varepsilon})\,\}$, then $|b-x_0|<\alpha$ and
$|a-x'|<\alpha$. Now, the condition (\ref{haus}) implies that for
$\delta$ small enough,
$$
\qquad \lim_{N\rightarrow +\infty}\nu_N\left(\partial {\cal
A}_N^{\sigma}\cap B(x_0,\delta)\subset {\cal
R}({\varepsilon})\right)=1.
$$
Putting together the previous facts, we obtain that
$$\limsup_{N\rightarrow +\infty}\,\,
\bigg|\nu_N\left(\frac{|{\cal H}_N(x_{0},\delta)|}{N}\right)-
|x_0x'\cdot i| \bigg|\,\leq\,2\alpha\,.$$ Remarking that
$|x_0x'\cdot i|= \delta|\cos\theta_0|$, we conclude the proof by
sending~$\alpha$ to~$0$. \,\, $\Box$
\medskip

\noindent
Now, let ${\cal L}^{\sigma}_{1,N}$ be the monotone path
$(v_1,\ldots,v_s)$ as defined in the subsection \ref{SEC}. We
obtain using the definition of ${\cal H}_{2,N}(x_0,\delta)$,
that $|{\cal H}_{2,N}(x_0,\delta)|$ is either $N_+({\cal
L}^{\sigma}_{1,N})$ or $N_-({\cal L}^{\sigma}_{1,N})$. This fact
together with the constatation that $\left|N_+({\cal
L}_{1,N})-N_-({\cal L}_{1,N})\right|\leq 1$, gives
$$
\left|\frac{|{\cal H}_{2,N}(x_0,\delta)|}{N}-\frac{ N_+({\cal
L}^{\sigma}_{1,N})}{N}\right|\leq \frac{1}{N}\,.
$$
Condition (\ref{c5}) together with the last inequality ensures,
since $\left||{\cal H}_{2,N}(x_0,\delta)|-C_N(x_0,\delta)\right|\leq
{1}$,
\begin{equation}\label{ef}
\lim_{N\rightarrow+\infty}{\nu_N}\kern-4pt\left(\frac{|{\cal
H}_{2,N}(x_0,\delta)|}{N}\right)={\delta}C(\theta_0).
\end{equation}
The two sets of indices ${\cal H}_{1,N}(x_0,\delta)$ and ${\cal
H}_{2,N}(x_0,\delta)$ form a partition of ${\cal H}_{N}(x_0,\delta)$,
 hence
\begin{equation}\label{ef1}
\lim_{N\rightarrow+\infty}\left|{\nu_N}\kern-4pt\left(\frac{|{\cal
H}_{1,N}(x_0,\delta)|}{N\delta}\right)-\left(\left|\cos\theta_0\right|
-{C(\theta_0)}\right)\right|=0.
\end{equation}
We obtain, collecting (\ref{l1}), (\ref{l2}), (\ref{l4}),
(\ref{e6}), (\ref{ef}), (\ref{ef1}) that
\begin{eqnarray*}
\lim_{N\rightarrow\infty}\left|
\frac{1}{\delta}{\nu_N}\kern-4pt\left(\int_{0}^{\delta}\,e_N^{1}(\alpha)\cdot{i}\,d\alpha\right)
+\sg(\cos\theta_0)\left(\cos^2\theta_0 +
{C(\theta_0)}\left(\left|\sin\theta_0\right|-
\left|\cos\theta_0\right|\right)\right)\right|=0.
\end{eqnarray*}
Using the same method, we prove that
\begin{eqnarray*}
\lim_{N\rightarrow\infty}\left|\frac{1}{\delta}{\nu_N}\kern-4pt\left(\int_{0}^{\delta}\,e_N^{2}(\alpha)\cdot{i}\,d\,\alpha\right)
- \sg(\cos\theta_1) \left(\cos^2\theta_1+
{C(\theta_1)}\left(\left|\sin\theta_1\right|-
\left|\cos\theta_1\right|\right)\right)\right|=0.
\end{eqnarray*}
We finish the proof of theorem \ref{T1} by combining proposition
\ref{a1} together with the two last limits. \, $\Box$
\newpage
\subsection{ End of the proof for Jordan curves
(theorem~\ref{T2}).}  To extend the proofs to Jordan curves, we
have to generalize lemmas \ref{lem5sp} and \ref{lembis5sp} as
follows.
\begin{lemma}\label{jorlem5sp}
We have
$$\displaylines{
\lim_{\delta\to 0} \limsup_{N\rightarrow \infty} \bigg|
\nu_N\bigg( \frac{1}{\delta N} \sum_{l\in {\cal
H}_{1,N}(x_0,\delta)} \frac{{(x_0a_{l+1})}\cdot{{i}}}
{|{x_0-a_{l+1}}|}\bigg) \,-\, \nu_N\bigg(\frac{|{\cal
H}_{1,N}(x_0,\delta)|}{N}\bigg) \cos\theta_0\bigg|\,=\,0\,,\cr
\lim_{\delta\to 0}\limsup_{N\rightarrow \infty} \bigg| \nu_N\bigg(
\frac{1}{\delta N} \sum_{l\in {\cal H}_{2,N}(x_0,\delta)}
\frac{{(x_0a_{l+1})}\cdot{{j}}} {|{x_0-a_{l+1}}|}\bigg) \,-\,
\nu_N\bigg(\frac{|{\cal H}_{2,N}(x_0,\delta)|}{N}\bigg)
\sin\theta_0\bigg|\,=\,0\,. }$$
\end{lemma}
\noindent{\bf{Proof of lemma \ref{jorlem5sp}.}} We only prove the
first limit since the argument for the second limit is similar.
Let $u$ be a unit vector tangent to~$\gamma$ at $x_0$ and let $v$
be such that $(u,v)$ is a direct basis. For $\varepsilon>0$, let
${\cal R}({\varepsilon})$ be the strip of width $2\varepsilon$
centered on the tangent line $T_{x_0}\gamma$, i.e.,
$${\cal R}(\varepsilon)\,=\,
\big\{\,x\in\BBr^2:|x_0x\cdot v|\leq\varepsilon\,\}\,.$$ Since
$T_{x_0}\gamma$ is the tangent to $\gamma$ at $x_0$, we have
\begin{equation}\label{hausbis}
\lim_{\delta\to 0} \,{1\over\delta}d_H\Big(\gamma\cap
B(x_0,\delta), T_{x_0}\gamma\cap B(x_0,\delta)\Big)\,=\,0\,.
\end{equation}
 Let
$0<\varepsilon<1$, there exists $\delta_0>0$ such that, for
$\delta<\delta_0$,
$$d_H\Big(\gamma\cap B(x_0,\delta),
T_{x_0}\gamma\cap
B(x_0,\delta)\Big)\,\leq\,\varepsilon\delta/4\,.$$ This fact
together with condition (\ref{haus}) implies that there exists
$\delta_1>0$ such that
$$
\forall\delta<\delta_1\qquad \lim_{N\rightarrow
+\infty}\nu_N\left(\partial {\cal A}_N^{\sigma}\cap
B(x_0,\delta)\subset {\cal R}({\varepsilon}\delta/2)\right)=1.
$$
Let $\delta>0$ be such that $\delta<\min(\delta_0,\delta_1)$. Let
$x\in{\cal R}(\varepsilon\delta)\setminus
B(x_0,\sqrt{\varepsilon}\delta)$. We have
$$\displaylines{x_0x\cdot i\,=\,
(x_0x\cdot u)\cos\theta_0- (x_0x\cdot v)\sin\theta_0\,,\cr
|x_0-x|^2\,=\, (x_0x\cdot u)^2+ (x_0x\cdot v)^2\,,}$$ whence
$$
\frac{x_0x\cdot i}{|x_0-x|}\,=\, \bigg(1- \frac{(x_0x\cdot
v)^2}{|x_0-x|^2} \bigg)^{1/2} \cos\theta_0- \frac{(x_0x\cdot
v)}{|x_0-x|}\sin\theta_0\,$$ and
$$\bigg|
\frac{x_0x\cdot i}{|x_0-x|}-\cos\theta_0\bigg|\,\leq\,
1-\sqrt{1-\varepsilon}+ \sqrt{\varepsilon}\,\leq\,2
\sqrt{\varepsilon}\,.$$ If the event $\{\,\partial {{\cal
A}}_N^{\sigma}\cap B(x_0,\delta)\subset {\cal
R}({\varepsilon\delta}/2)\,\}$ occurs, then $a_{l+1}\in {\cal
R}({\varepsilon}\delta)\setminus B(x_0,\sqrt{\varepsilon}\delta)$
for $l\in {\cal H}_{N}(x_0,\delta)\setminus {\cal
H}_{N}(x_0,\sqrt{\varepsilon}\delta)$, and thus
$$
\limsup_{N\rightarrow \infty}\,\,\nu_N\left(\sup_{l\in {\cal
H}_{N}(x_0,\delta)\setminus {\cal
H}_{N}(x_0,\sqrt{\varepsilon}\delta)}
\left|\frac{({x_0a_{l+1})}\cdot{{i}}} {|{x_0-a_{l+1}}|}
-\cos\theta_0\right|\right)\,\leq\, 2\sqrt{\varepsilon}\,.$$
Moreover, we have $\big|{\cal H}_{N}(x_0,\sqrt{\varepsilon}\delta)
\big|\,\leq\,2N\delta\sqrt{\varepsilon}$, whence, by splitting the
sum over ${\cal H}_{N}(x_0,\sqrt{\varepsilon}\delta)$ and ${\cal
H}_{N}(x_0,\delta)\setminus {\cal
H}_{N}(x_0,\sqrt{\varepsilon}\delta)$, we obtain
$$\limsup_{N\rightarrow \infty}
\bigg|\nu_N\bigg( \frac{1}{\delta N} \sum_{l\in {\cal
H}_{1,N}(x_0,\delta)} \frac{{(x_0a_{l+1})}\cdot{{i}}}
{|{x_0-a_{l+1}}|}\bigg) \,-\, \nu_N\bigg(\frac{|{\cal
H}_{1,N}(x_0,\delta)|}{N\delta}\bigg) \cos\theta_0\bigg|\,\leq\,
8\sqrt{\varepsilon}\,.$$ This inequality being valid for all
$\delta$ small enough, the proof is completed.
 \ \ $\Box$
\begin{lemma}\label{jorlembis5sp}
We have
\begin{equation}\label{jorc1bis}
\lim_{\delta\to 0}\limsup_{N\rightarrow +\infty}
\bigg|\nu_N\left(\frac{|{\cal
H}_N(x_{0},\delta)|}{\delta N}\right) -|\cos\theta_0|\bigg|\,=\,0\, .
\end{equation}
\end{lemma}
{\bf{Proof of lemma \ref{jorlembis5sp}.}} We denote, as in the
proof of lemma \ref{lembis5sp}, by $a$ and $x'$ the points of
$\partial B(x_0,\delta)\setminus B(s,r)$ belonging respectively to
$\partial {\cal A}_N^{\sigma}$ and to $\gamma$. Let $b$ be the
point of $\partial {\cal A}_N^{\sigma}\cap \partial
B(s,r)\setminus B(x_1,\delta)$.

\begin{figure}[hbt]
\vskip 8truecm \vbox{ \centerline{ \psset{unit=1.5cm}
\pspicture(2.5,2.5)(5,-0.5)
\psline[linewidth=1.3pt,linecolor=red](5.7,7.5)(3.21,4.5)
\psline[linewidth=1.75pt,linecolor=red](5.75,6.85)(5.5,6.85)(5.25,6.85)(5.25,6.65)(5.05,6.65)(5.05,6.45)(4.65,6.45)(4.65,6.05)
(4.65,5.6) \rput(4.65,5.95){$\bullet$} \rput(4.84,5.99){$b$}
\psline[linewidth=0.95pt,linecolor=red,linestyle=dashed](3.21,4.5)(3.4,4.5)(3.6,4.5)(4,4.5)(5,4.5)
\psarc[arcsepA=2pt,arcsepB=2pt]{->}(3.21,4.5){0.4}{0}{59}
\rput(3.8,4.7){$\theta_0$}
\pscurve[linewidth=2pt,linecolor=blue](6.5,6.8)(6,6.75)(5.25,6.5)(4.5,6)(4,4.5)(4.5,3)(5.25,2.6)(6.75,3)(6.5,6.8)
\rput(5.71,6.7){$\bullet$} \rput(5.6,6.5){$x'$}
\rput(5.59,6.85){$\bullet$} \rput(5.8,6.96){$a$}
\rput(6,7.7){$T_{x_0}\gamma$} \rput(3.8,4.3){$s$}
\rput(4.5,3){$\bullet$} \rput(4.5,6){$\bullet$}
\rput(4,4.5){$\bullet$}
 \rput(4.47,6.27){$x_0$}
\rput(4.5,2.7){$x_1$}
\psarc[arcsepA=2pt,arcsepB=2pt](4,4.5){1.6}{0}{362}
\psarc[arcsepA=2pt,arcsepB=2pt](4.5,6){1.38}{0}{362}
\psline[linewidth=2pt,linecolor=blue]{->}(6.6,2.3)(6.3,2.7)
\rput(6.8,2.4){$\gamma$}
\psline[linewidth=1pt,linecolor=blue]{->}(4,4.5)(3.25,3.1)
\rput(3.4,3.78){$r$}
\psline[linewidth=1pt,linecolor=blue]{->}(4.5,6)(3.25,6.67)
\rput(4,6.5){$\delta$}
\endpspicture
} } \vskip -4truecm {\centerline{The random points $a$ and $b$ are
approximated, for $N$ large enough,}} {\centerline{ respectively
by $x'$ and $x_0$. }}
\end{figure}

\noindent We
 have, by definition of ${\cal H}_N(x_{0},\delta)$,
$$
\left|\frac{|{\cal H}_N(x_{0},\delta)|}{N} -|ba\cdot i|\right|\leq
\frac{2}{N}.
$$
We suppose that $r$ is small enough so that
$T_{x_0}\gamma$ is not tangent to the circle $\partial B(s,r)$.
Let $\alpha>0$. There exists $\varepsilon>0$ depending on~$\alpha$
and the angle of the tangent $T_{x_0}\gamma$ with
$\partial B(s,r)$
such that
$$\forall\delta >0
\qquad\{\,\partial {\cal A}_N^{\sigma}\cap B(x_0,\delta)\subset
{\cal R}({\varepsilon}\delta)\,\}\quad\Rightarrow \quad
|b-x_0|<\alpha\delta\,,\quad |a-x'|<\alpha\delta\,.$$ Now as in
the proof of lemma \ref{jorlem5sp}, the condition (\ref{haus})
together with (\ref{hausbis}) implies that for $\delta$ small
enough,
$$
\lim_{N\rightarrow +\infty}\nu_N\left(\partial {\cal
A}_N^{\sigma}\cap B(x_0,\delta)\subset {\cal
R}({\varepsilon}\delta)\right)=1.
$$
Putting together the previous facts, we obtain that
$$\limsup_{N\rightarrow +\infty}\,\,
\bigg|\nu_N\left(\frac{|{\cal H}_N(x_{0},\delta)|}{N}\right)-
|x_0x'\cdot i| \bigg|\,\leq\,2\alpha\delta\,.$$ Remarking that
$|x_0-x'|= \delta$, and that
$$
\lim_{\delta\rightarrow 0}\frac{|x_0x'\cdot
i|}{|x_0-x'|}=|\cos\theta_0|,
$$
we conclude the proof by sending~$\alpha$ to~$0$. \,\, $\Box$
\begin{coro}\label{jjorlembis5sp}
We have
\begin{equation}\label{jjorc1bis}
\lim_{\delta\to 0}\limsup_{N\rightarrow +\infty}
\bigg|\nu_N\left(\frac{|{\cal H}_{1,N}(x_{0},\delta)|}{\delta
N}\right) -\left(|\cos\theta_0|-C(\theta_0)\right)\bigg|\,=\,0\, .
\end{equation}
\begin{equation}\label{jjorc2bis}
\lim_{\delta\to 0}\limsup_{N\rightarrow +\infty}
\bigg|\nu_N\left(\frac{|{\cal H}_{2,N}(x_{0},\delta)|}{\delta
N}\right) -C(\theta_0)\bigg|\,=\,0\, .
\end{equation}
\end{coro}
\noindent{\bf{Proof of corollary \ref{jjorlembis5sp}.}} The limit
in (\ref{jjorc2bis}) is deduced from the condition (\ref{c5bis})
since by definition
$$
\left||{\cal H}_{2,N}(x_0,\delta)|-C_N(x_0,\delta)\right|\leq
\frac{1}{N}.
$$
The first limit is deduced by combining (\ref{jjorc2bis}) and the
result of lemma \ref{jorlembis5sp}, since ${\cal
H}_{2,N}(x_0,\delta)$ and ${\cal H}_{1,N}(x_0,\delta)$ form a
partition of ${\cal H}_{N}(x_0,\delta)$.
\\
$\Box$
\\
We obtain, collecting (\ref{l1}), (\ref{l2}), (\ref{l4}),
(\ref{l31}), lemma \ref{jorlem5sp}, (\ref{jjorc1bis}) and
(\ref{jjorc2bis}) that
\begin{eqnarray*}
\lim_{\delta\rightarrow 0}\limsup_{N\rightarrow\infty}\left|
\frac{1}{\delta}{\nu_N}\kern-4pt\left(\int_{0}^{\delta}\,e_N^{1}(\alpha)\cdot{i}\,d\alpha\right)
+\sg(\cos\theta_0)\left(\cos^2\theta_0 +
{C(\theta_0)}\left(\left|\sin\theta_0\right|-
\left|\cos\theta_0\right|\right)\right)\right|=0.
\end{eqnarray*}
Using the same method, we prove that
\begin{eqnarray*}
\lim_{\delta\rightarrow 0}
\limsup_{N\rightarrow\infty}\left|\frac{1}{\delta}{\nu_N}\kern-4pt\left(\int_{0}^{\delta}\,e_N^{2}(\alpha)\cdot{i}\,d\,\alpha\right)
- \sg(\cos\theta_1) \left(\cos^2\theta_1+
{C(\theta_1)}\left(\left|\sin\theta_1\right|-
\left|\cos\theta_1\right|\right)\right)\right|=0.
\end{eqnarray*}
We get the expression of $\lim_{\delta\rightarrow
0}\limsup_{N\rightarrow \infty}\, \nu_N\big({\bf{A}}^{\sigma,
\gamma}_N(s,r,\delta)\big)$ of theorem \ref{T2} by combining
proposition \ref{a1} together with the two last limits. The
$\liminf$ can be handled similarly.\, $\Box$
\newpage
\subsection{Proof of proposition \ref{theo1heu1}}
In this case the condition \ref{haus} of theorem \ref{T1} is
satisfied and we have only to check the limit in (\ref{c5bis}) and
to precise the value of the function $C$ defined there. For this,
we need the following lemma.
\begin{lemma}\label{lem5bis}
Let $A$ and $B$ be two points of $\BBr^2$. For each fixed integer
$N$, let ${\cal L}_N$ denote one of the two maximal subpaths of
$\partial([AB]_N)$ not crossing the line $(AB)$.
 Let $N_+({\cal L}_N)$ be as
defined in (\ref{def}). Then $$\lim_{N\rightarrow+\infty}\frac{
N_+({\cal L}_N)}{N}= \left|(AB)\cdot{i}\right|\wedge
\left|(AB)\cdot{j}\right|, $$ where $a\wedge b=\min(a,b)$.
\end{lemma}

\begin{figure}[hbt]
\vbox{ \centerline{ \psset{unit=0.9cm} \pspicture(0,-1)(9,7)
\psgrid[subgriddiv=0,griddots=10,gridlabels=0](7,6)
\rput(0.3,0.2){A} \rput(6.5,4.6){B} \rput(7.4,5){$B_N$}
\rput(-0.2,-0.1){$A_N$}
\psline[linewidth=1.5pt,showpoints=true](0.6,0.3)(6.6,4.3)
\psline[linewidth=1.5pt,showpoints=true]{->}(7,5)(6,5)(6,4)(5,4)(4,4)(4,3)(3,3)(3,2)(2,2)(1,2)(1,1)(0,1)(0,0)
\psline(0.6,0.2)(6.5,0.2)
\psarc[arcsepA=2pt,arcsepB=2pt]{->}(0.3,0.2){1.5}{0}{30}
\rput(2.25,0.7){$\theta$}
\endpspicture
} }
\end{figure}

{\hspace{-7mm}}{\bf{{Proof of lemma \ref{lem5bis}.}}} We suppose
without loss of generality that ${AB}\cdot{i}$ and ${AB}\cdot{j}$
are positive. Let $\theta$ denote the angle between ${AB}$ and
${i}$. We consider only the case
 $0\leq \tan\theta< 1$, since the proofs for the cases $\tan\theta> 1$ and $\tan\theta=
1$ are similar. We denote by $A_N$ and $B_N$ the extreme points
of ${\cal L}_N$. Our task is to prove that
\begin{equation}\label{limlem5}
(A_NB_N)\cdot{j}=\frac{ N_{+}({\cal L}_N)}{N}.
\end{equation}
The identity (\ref{limlem5}) will prove lemma \ref{lem5bis} since
$\lim_{N\rightarrow +\infty} (A_NB_N)\cdot{j} =(AB)\cdot{j}$ and
$0\leq \tan\theta< 1$.
\\
We first prove the equality (\ref{limlem5}) for
 $N_{+}({\cal L}_N)=1$.
 When $N_{+}({\cal L}_N)=1$,
  the path ${\cal L}_N$
contains a unique monotone path ${\cal L}'_N=(v_1,w_1,\ldots,w_r)$
such that $v_1\cdot j=0$, $v_1\cdot w_1=0$ and $w_1\cdot
i=\ldots=w_r\cdot i=0$. These vectors are drawn on the lattice
$\BBz_N^2$ and arranged according to the direct sense.

\noindent Let $C_1$, $C_2$  be the two points of $(AB)$ such that
$(C_1C_2)\cdot{i}={1}/{N}$ and that the path ${\cal L}'_N$ cover
the segment $[C_1, C_2]$. By construction $$
(C_1C_2)\cdot{j}>\frac{r-1}{N}, $$ hence $$
\tan\theta=\frac{(C_1C_2)\cdot{j}}{(C_1C_2)\cdot{i}}> r-1. $$

\begin{figure}[hbt]
\vbox{ \centerline{ \psset{unit=0.9cm} \pspicture(0,-1)(9,6)
\psgrid[subgriddiv=0,griddots=10,unit=0.75,
gridlabels=0](3,0)(8,8)
\psline[linewidth=1.5pt](4.5,4.5)(3.75,4.5)(3.75,1.5)
\psline[linewidth=0.5pt,showpoints=true](4.5,4)(3.75,1)
\rput(3.75,0.6){$C_1$} \rput(4.8,4){$C_2$} \rput(4,4.75){$v_1$}
\rput(3.5,4.2){$w_1$} \rput(3.5,1.9){$w_r$}
\endpspicture
} }
\end{figure}

\noindent Since $0\leq \tan\theta< 1$, we deduce that $r=1$. Now
let  $B_1$, $B_2$ be the two points of $(AB)$ belonging to the
boundary of the box of $[AB]_N$ that contains the point $B$. Since
$0\leq \tan\theta<1$, we have $$
\left|{(B_1B_2)}\cdot{j}\right|<\frac{1}{N}. $$

\begin{figure}[hbt]
\vskip-2.5truecm \vbox{ \centerline{ \psset{unit=0.9cm}
\pspicture(0,-1)(9,7) \psgrid[subgriddiv=0,griddots=10,unit=0.75,
gridlabels=0](0.75,1.5)(7,6)
\psline[linewidth=1.5pt]{->}(4.5,3.75)(3.75,3.75)
\rput(4.75,4){$B_N$}
\psline[linewidth=1.5pt]{->}(3.75,3.75)(3,3.75)
\psline[linewidth=1.5pt]{->}(3,3.75)(3,3)
\psline[linewidth=1.5pt]{->}(3,3)(2.25,3)
\psline[linewidth=1.5pt]{->}(2.25,3)(1.5,3) \rput(1.2,3){$A_N$}
\psline[linewidth=1.5pt,showpoints=true](4.5,3.3)(4,3.15)(3.75,3.09)(1.5,2.4)
\rput(4.9,3.3){$B_2$} \rput(3.6,2.6){$B_1$} \rput(4.3,2.7){$B$}
\rput(4.25,4){$e_1$} \rput(3.5,4){$e_2$} \rput(2.6,3.6){$w_1$}
\rput(2.5,3.3){$f_1$}
\endpspicture
} }
\end{figure}
\vspace{-2cm} \noindent This fact together with $r=1$ proves that
 the path ${\cal L}_N$ is equal to $(e_1,\ldots,e_m,w_1,f_1,\ldots,f_{n})$,
 where the vectors $(e_i)$ and $(f_i)$ are  copies of
 the vector $v_1$, so that they are all horizontal.
 Hence
 $
 (A_NB_N)\cdot{j}={1}/{N}.
 $
The general case when $N_{+}({\cal
L}_N)>1$ is proved by induction on $N_{+}({\cal L}_N)$.
 $\Box$
 \\
 \\
\noindent {\bf{Proofs for polygons.}} Lemma \ref{lem5bis} together
with theorem \ref{T1} yield the control of
${\bf{A}}_N^{\sigma,\Gamma}(s,r,\delta)$ for a class of regular
polygons $\Gamma$ defined as follows.
\\
\\
{\bf{m-smooth polygons.}} Let $s_1,\ldots,s_m$ be $m$ points of
$\BBr^2$. We denote by $\Gamma(s_1,\ldots,s_m)$ or by $\Gamma$, if
there is no ambiguities, the polygon in $\BBr^2$ linking the
points $[s_1,s_2,\ldots,s_m,s_1]$; the points $s_1,s_2,\ldots,s_m$
are then the corner points of $\Gamma$. We suppose that the points
$s_1,s_2,\ldots,s_m$ are arranged counterclockwise. By convention,
we set $s_0=s_m$. To each site $s_i$, we associate two oriented
angles $\theta_i(s_i)$ and $\theta_{i-1}(s_i)$ such that
$\theta_{i-1}(s_i)$ (respectively $\theta_i(s_i)$) is the oriented
angle between the half horizontal axis $[0,+\infty[$ and the
segment $[s_i,s_{i-1}[$ (respectively $[s_i,s_{i+1}[$).
\newpage

\noindent Finally, we suppose that
 $\Gamma$ encloses a connected, compact, bounded set $U$ of
  $\BBr^2$ i.e. $\Gamma={\partial}{U}$ and that $\Gamma\cap \BBz^2/N=\emptyset$ for all $N\geq 1$.
  \\
  \\
{\bf{Initial condition.}} We will consider $\sigma$ the spin
configuration associated to the polygon $\Gamma$ at step $N$.

\begin{figure}[hbt]
\vspace{2cm} \vbox{ \centerline{ \psset{unit=0.75cm}
\pspicture(0,-1)(15,15)
\pspolygon[fillstyle=solid,fillcolor=lightgray,
linewidth=1pt](2,1)(2,2)(3,2)(3,4)(2,4)(2,5)(3,5)
(3,6)(4,6)(4,7)(5,7)(5,8)(6,8)(6,9)(7,9)(7,6)
(8,6)(8,10)(2,10)(2,11)(3,11)(3,13)(6,13)
(6,14)(8,14)(8,13)(13,13)(13,9)(14,9)(14,4)
(13,4)(13,3)(11,3)(11,2)(8,2)(8,1)
\pspolygon[fillstyle=solid,fillcolor=gray,
linewidth=1pt](2.8,1.6)(5.3,3.1)(2.5,4.5)
(6.6,8.2)(6.8,6.1)(7.3,4)(8.8,6.2)(9.1,8.6)(7.3,11.2)(5.1,10.3)
(2.4,10.1)(4.5,12.6)(7.5,13.1)(10.5,12.3)(12.5,12.7)
(12.2,10.7)(13.1,8.3)(13.5,4.2)(10.5,3.3)(9.2,2.2)
\psdots[dotscale=1.2](2.8,1.6)(5.3,3.1)(2.5,4.5)
(6.6,8.2)(6.8,6.1)(7.3,4)(8.8,6.2)(9.1,8.6)(7.3,11.2)(5.1,10.3)
(2.4,10.1)(4.5,12.6)(7.5,13.1)(10.5,12.3)(12.5,12.7)
(12.2,10.7)(13.1,8.3)(13.5,4.2)(10.5,3.3)(9.2,2.2)
\psgrid[subgriddiv=1,griddots=4,gridlabels=0](0,0)(0,0)(15,15)
\rput(10.5,12.7){$s_1$} \rput(7.5,13.6){$s_2$}
\rput(4.5,13.4){$s_3$}
\psline[linestyle=dashed,linewidth=1pt](0.5,4.5)(6,4.5)
\psarc[arcsepA=2pt,arcsepB=2pt]{->}(2.5,4.5){3}{0}{42.064}
\rput(6.1,5.7){$\theta_{i-1}$}
\psarc[arcsepA=2pt,arcsepB=2pt]{->}(2.5,4.5){1.75}{0}{333.435}
\rput(1.5,2.7){$\theta_{i}$} \rput(2.5,5.5){$s_i$}
\rput(6.8,8.7){$s_{i-1}$} \rput(5.5,3.5){$s_{i+1}$}
\pspolygon[fillstyle=solid,fillcolor=lightgray,
linewidth=0pt](-3,8)(-2,8)(-2,9)(-3,9) \rput(-1,8.5){$=\sigma$}
\pspolygon[fillstyle=solid,fillcolor=gray,
linewidth=0pt](-3,6)(-2,6)(-2,7)(-3,7) \rput(-1,6.5){$=U_0$}
\rput(2.5,8.5){$\Gamma$}
\psline[linewidth=1pt]{->}(2.8,8.7)(4,10.2)
\endpspicture
} \centerline{A polygon~$\Gamma$ and the configuration~$\sigma$} }
\end{figure}
\newpage
\noindent Lemma \ref{lem5bis} allows to apply theorem \ref{T1}
with $C(\theta_k)=|\sin\theta_k|\wedge |\cos\theta_k|$. Doing so,
we get the following proposition.
\begin{pro}\label{pro1}
Let $\Gamma$ be an m-smooth polygon in $\BBr^2$ associated to the
$m$ points $s_1,\ldots,s_m$ and let $\sigma_0:=\sigma_{0,N}$ be
 the associated initial configuration
at step $N$. Let $\theta_i\in [0,2\pi]$ (respectively
$\theta_{i-1}\in [0,2\pi]$)
 be the oriented angle between the half horizontal axis $[0,+\infty[$
 and the segment
 $[s_i,s_{i+1}[$ (respectively $[s_{i},s_{i-1}[$)
 with the convention that $s_0=s_m$. Then, for each $i=1,\ldots,m$, and for
 any
positive real numbers $r$, $\delta$ small enough, one has
$$\lim_{N\rightarrow+\infty}{\bf{A}}^{\sigma,\Gamma}_N(s_i,r,\delta)
=\frac{1}{4}\sin2\theta_{i-1}\left(2\BBone_{|\sin\theta_{i-1}|<|\cos\theta_{i-1}|}-1\right)
-\frac{1}{4}\sin2\theta_{i}\left(2\BBone_{|\sin\theta_{i}|<|\cos\theta_{i}|}-1\right)
$$ $$
+\frac{1}{2}\left(\sg(\tan\theta_i)\BBone_{|\sin\theta_{i}|<|\cos\theta_{i}|}-
\sg(\tan\theta_{i-1})\BBone_{|\sin\theta_{i-1}|<|\cos\theta_{i-1}|}\right)$$
$$+
\BBone_{\sin\theta_{i-1}\sin\theta_{i}>0}\left(\sg(\theta_{i}-\theta_{i-1})\BBone_{\cos\theta_{i-1}\cos\theta_{i}>0}+
\sg(\tan\theta_{i-1})
\BBone_{\cos\theta_{i-1}\cos\theta_{i}<0}\right).$$
\\
Hence,
\begin{itemize}
  \item if\ \
$(\theta_{i-1},\theta_i)\in[(2k+1)\frac{\textstyle \pi}{\textstyle
4},(2k+3)\frac{\textstyle \pi}{\textstyle 4}]\times
[(2k+5)\frac{\textstyle \pi}{\textstyle 4},(2k+7)\frac{\textstyle
\pi}{\textstyle 4}],$ with {\bf{$k\in\{0,2\}$}}, then $$
\lim_{N\rightarrow+\infty}{\bf{A}}^{\sigma,\Gamma}_N(s_i,r,\delta)
=
 \frac{1}{4}\left({\sin
2\theta_{i}}-{\sin 2\theta_{i-1}}\right). $$
\item if\ \
$(\theta_{i-1},\theta_i)\in[(2k+1)\frac{\textstyle \pi}{\textstyle
4},(2k+3)\frac{\textstyle \pi}{\textstyle 4}]\times
[(2k+5)\frac{\textstyle \pi}{\textstyle 4},(2k+7)\frac{\textstyle
\pi}{\textstyle 4}],$ with {\bf{$k\in\{1,3\}$}}, then $$
\lim_{N\rightarrow+\infty}{\bf{A}}^{\sigma,\Gamma}_N(s_i,r,\delta)
=
 \frac{1}{4}\left({\sin 2\theta_{i-1}}-{\sin 2\theta_i}\right). $$
\item if\ \
$(\theta_{i-1},\theta_i)\in[(2k+1)\frac{\textstyle\pi}{\textstyle
4},(2k+3)\frac{\textstyle\pi}{\textstyle 4}]^2$, with
{\bf{$k\in\{0,2\}$}}, then $$
\lim_{N\rightarrow+\infty}{\bf{A}}^{\sigma,\Gamma}_N(s_i,r,\delta)
=
 \frac{1}{4}\left(4\,\sg(\theta_{i}-\theta_{i-1})+{\sin
2\theta_{i}}-{\sin 2\theta_{i-1}}\right). $$
\item if\ \
$(\theta_{i-1},\theta_i)\in[(2k+1)\frac{\textstyle \pi}{\textstyle
4},(2k+3)\frac{\textstyle \pi}{\textstyle 4}]^2$, with
{\bf{$k\in\{1,3\}$}}, then $$
\lim_{N\rightarrow+\infty}{\bf{A}}^{\sigma,\Gamma}_N(s_i,r,\delta)
= \frac{1}{4}\left(4\,\sg(\sin(\theta_{i}-\theta_{i-1}))+{\sin
2\theta_{i-1}}-{\sin 2\theta_i}\right). $$
  \item if\ \
$(\theta_{i-1},\theta_i)\in[(2k+1)\frac{\textstyle \pi}{\textstyle
4},(2k+3)\frac{\textstyle\pi}{\textstyle 4}]\times
[(2k+3)\frac{\textstyle \pi}{\textstyle 4},(2k+5)\frac{\textstyle
\pi}{\textstyle 4}]\cup [(2k+3)\frac{\textstyle \pi}{\textstyle
4},(2k+5)\frac{\textstyle \pi}{\textstyle 4}]\times
[(2k+1)\frac{\textstyle \pi}{\textstyle 4},(2k+3)\frac{\textstyle
\pi}{\textstyle 4}] $, with $k\in\{0,2\}$, then $$
\lim_{N\rightarrow+\infty}{\bf{A}}^{\sigma,\Gamma}_N(s_i,r,\delta)
=\left\{
\begin{array}{rl}
\frac{\textstyle 1}{\textstyle 4}\left(2-{\sin 2\theta_{i}}-{\sin
2\theta_{i-1}}\right)
& if \ |\tan\theta_i|\leq 1,\ \ |\tan\theta_{i-1}|\geq 1\\
\\
 \frac{\textstyle 1}{\textstyle 4}\left(-2+{\sin 2\theta_{i-1}}+{\sin 2\theta_i}\right) & if \
|\tan\theta_i|\geq 1,\ \ |\tan\theta_{i-1}|\leq 1.
\end{array}
\right. $$
\item if\ \
$(\theta_{i-1},\theta_i)\in[(2k+1)\frac{\textstyle \pi}{\textstyle
4},(2k+3)\frac{\textstyle \pi}{\textstyle 4}]\times
[(2k+3)\frac{\textstyle \pi}{\textstyle 4},(2k+5)\frac{\textstyle
\pi}{\textstyle 4}]\cup [(2k+3)\frac{\textstyle \pi}{\textstyle
4},(2k+5)\frac{\textstyle \pi}{\textstyle 4}]\times
[(2k+1)\frac{\textstyle \pi}{\textstyle 4},(2k+3)\frac{\textstyle
\pi}{\textstyle 4}] $, with $k\in\{1,3\}$, then $$
\lim_{N\rightarrow+\infty}{\bf{A}}^{\sigma,\Gamma}_N(s_i,r,\delta)
=\left\{
\begin{array}{rl}
\frac{\textstyle 1}{\textstyle 4}\left(-2-{\sin 2\theta_{i}}-{\sin
2\theta_{i-1}}\right)
& if \ |\tan\theta_i|\leq 1,\ \ |\tan\theta_{i-1}|\geq 1\\
\\
 \frac{\textstyle 1}{\textstyle 4}\left(2+{\sin 2\theta_{i-1}}+{\sin 2\theta_i}\right) & if \
|\tan\theta_i|\geq 1,\ \ |\tan\theta_{i-1}|\leq 1.
\end{array}
\right. $$
\end{itemize}
\end{pro}
\vspace{1cm} {\bf{Remark.}} We denote by
$L_{\Gamma}(s_i)=\lim_{N\rightarrow+\infty}{\bf{A}}^{\sigma,\Gamma}_N(s_i,r,\delta),$
where $\Gamma$ is a polygon as described by proposition
\ref{pro1}. Then we can check the following comparison criterion.

\begin{figure}[hbt]
\vspace{10cm} \vbox{ \centerline{ \psset{unit=1.75cm}
\pspicture(2.5,2.5)(5,1.5)
\pspolygon[fillstyle=solid,linecolor=red,showpoints=true,
linewidth=0.5pt](6,7.5)(4.5,4.5)(2,5.25)(0.75,3)(4.5,1.5)(6,3)
\pspolygon[linewidth=1pt,
showpoints=true](4.5,7)(4.5,4.5)(3,6)(3,3)(6,1.5)(6.5,3.3)
\rput(6.55,7.52){$s_{i-1}$} \rput(4.99,4.65){$s_i$}
\rput(2.1,5.5){$s_{i+1}$}
\psarc[arcsepA=2pt,arcsepB=2pt](4.5,4.5){1}{0}{360}
\psline[linewidth=2pt]{->}(6.75,3)(5.99,3.3)
\psline[linewidth=2pt]{->}(6.75,1.5)(5.99,1.72)
\psline[linewidth=1pt,linecolor=green]{<->}(4.5,4.5)(3.8,3.72)
\rput(4,4.2){$r$}
 \rput(7.86,3){$\Gamma=\partial U$}
\rput(7.86,1.5){$\Gamma'=\partial U'$}
\endpspicture
} }
\end{figure}
\vspace{1cm} \centerline{ If $U\cap B(s_i,r)\subset U'\cap
B(s_i,r)$ for some $r>0$ and $s_i\in \Gamma\cap \Gamma'$,}
\centerline{then $L_{\Gamma}(s_i)\leq L_{\Gamma'}(s_i)$.}
\newpage
We illustrate the results of proposition \ref{pro1} with the help
 of the following pictures.

\begin{figure}[hbt]
\vspace{6cm} \vbox{ \centerline{ \psset{unit=1cm}
\pspicture(2.5,2.5)(5,-0.5) \psline[linewidth=1.5pt,
showpoints=true](6,7.5)(4.5,4.5)(6,1.5)
\pspolygon[fillstyle=solid,fillcolor=lightgray,
linewidth=0pt](6,7.5)(4.5,4.5)(6,1.5)(7.5,3)(8.25,3)
\psline[linestyle=dashed,linewidth=2pt](1.5,1.5)(3,3)(4.5,4.5)(6,6)(7.5,7.5)
\psline[linestyle=dashed,linewidth=2pt](1.5,7.5)(2.25,6.75)(3,6)(3.75,5.25)(4.5,4.5)(5.25,3.75)(6,3)(6.75,2.25)(7.5,1.5)
\psline[linestyle=dashed,linewidth=1pt](2.25,4.5)(3,4.5)(5.25,4.5)(6,4.5)(6.75,4.5)
\psgrid[subgriddiv=0,gridlabels=0,griddots=10,unit=0.75](1,1)(11,11)
\rput(4.5,2.3){$\theta_{i}$} \rput(5.75,5){$\theta_{i-1}$}
\rput(6.55,7.52){$s_{i-1}$} \rput(4.99,4.65){$s_i$}
\rput(6.4,1.52){$s_{i+1}$}
\psarc[arcsepA=2pt,arcsepB=2pt]{->}(4.5,4.5){1}{0}{69}
\psarc[arcsepA=2pt,arcsepB=2pt]{->}(4.5,4.5){1.7}{0}{297}
\psline[linewidth=2pt]{->}(7.5,6)(7,5.25)
\rput(7.75,6.5){$\Gamma$}
\endpspicture
} \vskip -1truecm \centerline{Here $
\lim_{N\rightarrow+\infty}{\bf{A}}^{\sigma,\Gamma}_N(s_i,r,\delta)
=
 \frac{1}{4}\left({\sin
2\theta_{i}}-{\sin 2\theta_{i-1}}\right). $} \centerline{ In the
first picture this limit is negative, while for the second one it
is positive. } }
\end{figure}

\vspace{7cm}

\begin{figure}[hbt]
\vbox{ \centerline{ \psset{unit=1cm} \pspicture(2.5,2.5)(5,1.5)
\psline[linewidth=1.5pt, showpoints=true](6,7.5)(4.5,4.5)(6,1.5)
\pspolygon[fillstyle=solid,fillcolor=lightgray,
linewidth=0pt](6,7.5)(4.5,4.5)(6,1.5)(3,2.25)(2.25,6)
\psline[linestyle=dashed,linewidth=2pt](1.5,1.5)(3,3)(4.5,4.5)(6,6)(7.5,7.5)
\psline[linestyle=dashed,linewidth=2pt](1.5,7.5)(2.25,6.75)(3,6)(3.75,5.25)(4.5,4.5)(5.25,3.75)(6,3)(6.75,2.25)(7.5,1.5)
\psline[linestyle=dashed,linewidth=1pt](2.25,4.5)(3,4.5)(5.25,4.5)(6,4.5)(6.75,4.5)
\psgrid[subgriddiv=0,gridlabels=0,griddots=10,unit=0.75](1,1)(11,11)
\rput(4.5,2.3){$\theta_{i-1}$} \rput(5.75,5){$\theta_{i}$}
\rput(6.55,7.52){$s_{i+1}$} \rput(4.99,4.65){$s_i$}
\rput(6.4,1.52){$s_{i-1}$}
\psarc[arcsepA=2pt,arcsepB=2pt]{->}(4.5,4.5){1}{0}{69}
\psarc[arcsepA=2pt,arcsepB=2pt]{->}(4.5,4.5){1.7}{0}{297}
\psline[linewidth=2pt]{->}(3.75,7.5)(3.75,6.7)
\rput(3.75,7.9){$\Gamma$}
\endpspicture
} }
\end{figure}
\newpage
\centerline{In the following picture, we have $
\lim_{N\rightarrow+\infty}{\bf{A}}^{\sigma,\Gamma}_N(s_i,r,\delta)
=
 \frac{1}{4}\left({\sin
2\theta_{i-1}}-{\sin 2\theta_{i}}\right). $} \centerline{ This
limit is negative.} \vspace{6cm}
\begin{figure}[hbt]
\vbox{ \centerline{ \psset{unit=1cm} \pspicture(2.5,2.5)(5,2)
\psline[linewidth=1.5pt, showpoints=true](7.5,6)(4.5,4.5)(1.5,6)
\pspolygon[fillstyle=solid,fillcolor=lightgray,
linewidth=0pt](7.5,6)(4.5,4.5)(1.5,6)(3,7.5)(6,7.5)
\psline[linestyle=dashed,linewidth=2pt](3.75,3.75)(4.5,4.5)(6,6)(7.5,7.5)
\psline[linestyle=dashed,linewidth=2pt](1.5,7.5)(2.25,6.75)(3,6)(3.75,5.25)(4.5,4.5)(5.25,3.75)
\psline[linestyle=dashed,linewidth=1pt](2.25,4.5)(3,4.5)(5.25,4.5)(6,4.5)(6.75,4.5)
\psgrid[subgriddiv=0,gridlabels=0,griddots=10,unit=0.75](0,2)(11,11)
\rput(6,4.75){$\theta_{i}$} \rput(4.5,6.75){$\theta_{i-1}$}
\rput(8,6){$s_{i+1}$} \rput(5.25,4.65){$s_i$}
\rput(1.3,5.7){$s_{i-1}$}
\psarc[arcsepA=2pt,arcsepB=2pt]{->}(4.5,4.5){1}{0}{33}
\psarc[arcsepA=2pt,arcsepB=2pt]{->}(4.5,4.5){1.7}{0}{155}
\psline[linewidth=2pt]{->}(1.75,4.75)(2.65,5.4)
\rput(1.75,4.25){$\Gamma$}
\endpspicture
} }
\end{figure}

\vspace{7.5cm}
\begin{figure}[hbt]
\vbox{ \centerline{ \psset{unit=1cm} \pspicture(2.5,2.5)(5,1.5)
\psline[linewidth=1.5pt, showpoints=true](6,7.5)(4.5,4.5)(3,6.75)
\pspolygon[fillstyle=solid,fillcolor=lightgray,
linewidth=0pt](6,7.5)(4.5,4.5)(3,6.75)(0.75,3)(4.5,1.5)(6,3)
\psline[linestyle=dashed,linewidth=2pt](1.5,1.5)(3,3)(4.5,4.5)(6,6)(7.5,7.5)
\psline[linestyle=dashed,linewidth=2pt](1.5,7.5)(2.25,6.75)(3,6)(3.75,5.25)(4.5,4.5)(5.25,3.75)(6,3)(6.75,2.25)(7.5,1.5)
\psline[linestyle=dashed,linewidth=1pt](2.25,4.5)(3,4.5)(5.25,4.5)(6,4.5)(6.75,4.5)
\psgrid[subgriddiv=0,gridlabels=0,griddots=10,unit=0.75](0,1)(11,11)
\rput(4.5,6.8){$\theta_{i}$} \rput(5.75,5){$\theta_{i-1}$}
\rput(6.55,7.52){$s_{i-1}$} \rput(4.99,4.75){$s_i$}
\rput(2.5,6.899){$s_{i+1}$}
\psarc[arcsepA=2pt,arcsepB=2pt]{->}(4.5,4.5){1}{0}{68}
\psarc[arcsepA=2pt,arcsepB=2pt]{->}(4.5,4.5){1.7}{0}{126}
\psline[linewidth=2pt]{->}(6.75,3)(5.99,3.3)
\rput(7.5,3){$\Gamma$}
\endpspicture
} \vskip 1truecm \centerline{Here $
\lim_{N\rightarrow+\infty}{\bf{A}}^{\sigma,\Gamma}_N(s_i,r,\delta)
=
 \frac{1}{4}\left(4+{\sin
2\theta_{i}}-{\sin 2\theta_{i-1}}\right). $} \centerline{ This
limit is positive.} }
\end{figure}
\newpage
\centerline{In the following picture, we have $
\lim_{N\rightarrow+\infty}{\bf{A}}^{\sigma,\Gamma}_N(s_i,r,\delta)
=
 \frac{1}{4}\left(2-{\sin
2\theta_{i-1}}-{\sin 2\theta_{i}}\right). $} \centerline{ This
limit is positive.} \vspace{6cm}
\begin{figure}[hbt]
\vbox{ \centerline{ \psset{unit=1cm} \pspicture(2.5,2.5)(5,1.5)
\psline[linewidth=1.5pt, showpoints=true](6,7.5)(4.5,4.5)(2,5.25)
\pspolygon[fillstyle=solid,fillcolor=lightgray,
linewidth=0pt](6,7.5)(4.5,4.5)(2,5.25)(0.75,3)(4.5,1.5)(6,3)
\psline[linestyle=dashed,linewidth=2pt](1.5,1.5)(3,3)(4.5,4.5)(6,6)(7.5,7.5)
\psline[linestyle=dashed,linewidth=2pt](1.5,7.5)(2.25,6.75)(3,6)(3.75,5.25)(4.5,4.5)(5.25,3.75)(6,3)(6.75,2.25)(7.5,1.5)
\psline[linestyle=dashed,linewidth=1pt](2.25,4.5)(3,4.5)(5.25,4.5)(6,4.5)(6.75,4.5)
\psgrid[subgriddiv=0,gridlabels=0,griddots=10,unit=0.75](0,1)(11,11)
\rput(4.5,6.8){$\theta_{i}$} \rput(5.75,5){$\theta_{i-1}$}
\rput(6.55,7.52){$s_{i-1}$} \rput(4.99,4.65){$s_i$}
\rput(2.1,5.5){$s_{i+1}$}
\psarc[arcsepA=2pt,arcsepB=2pt]{->}(4.5,4.5){1}{0}{68}
\psarc[arcsepA=2pt,arcsepB=2pt]{->}(4.5,4.5){1.7}{0}{165}
\psline[linewidth=2pt]{->}(6.75,3)(5.99,3.3)
\rput(7.5,3){$\Gamma$}
\endpspicture
} }
\end{figure}

\vspace{7.5cm}
\begin{figure}[hbt]
\vbox{ \centerline{ \psset{unit=1cm} \pspicture(2.5,2.5)(5,1.5)
\psline[linewidth=1.5pt,
showpoints=true](4.5,1.5)(4.5,4.5)(2.25,5.25)
\pspolygon[fillstyle=solid,fillcolor=lightgray,
linewidth=0pt](4.5,4.5)(2.25,5.25)(0.75,3)(4.5,1.5)
\psline[linestyle=dashed,linewidth=2pt](1.5,1.5)(3,3)(4.5,4.5)(6,6)(7.5,7.5)
\psline[linestyle=dashed,linewidth=2pt](1.5,7.5)(2.25,6.75)(3,6)(3.75,5.25)(4.5,4.5)(5.25,3.75)(6,3)(6.75,2.25)(7.5,1.5)
\psline[linestyle=dashed,linewidth=1pt](2.25,4.5)(3,4.5)(5.25,4.5)(6,4.5)(6.75,4.5)
\psgrid[subgriddiv=0,gridlabels=0,griddots=10,unit=0.75](0,1)(11,11)
\rput(4.5,6.8){$\theta_{i}$} \rput(5.75,5){$\theta_{i-1}$}
\rput(2.25,5.5){$s_{i+1}$} \rput(4.99,4.65){$s_i$}
\rput(4.5,1){$s_{i-1}$}
\psarc[arcsepA=2pt,arcsepB=2pt]{->}(4.5,4.5){1}{0}{270}
\psarc[arcsepA=2pt,arcsepB=2pt]{->}(4.5,4.5){1.7}{0}{162}
\psline[linewidth=2pt]{->}(5.25,1.75)(4.5,1.75)
\rput(5.25,2.25){$\Gamma$}
\endpspicture
} \vskip 1truecm \centerline{Here $
\lim_{N\rightarrow+\infty}{\bf{A}}^{\sigma,\Gamma}_N(s_i,r,\delta)
=
 \frac{1}{4}\left(-2-{\sin
2\theta_{i}}-{\sin 2\theta_{i-1}}\right). $} \centerline{ This
limit is negative.} }
\end{figure}
\newpage

\noindent {\bf{Proofs for Jordan curves.}} We consider now the
case of Jordan curves. In order to apply theorem \ref{T2}, we have
to check the condition (\ref{c5bis}). For this, we generalize the
lemma \ref{lem5bis} as follows.
 \begin{lemma}\label{ana}
Let $f$ be a monotone function of class ${\cal C}_1$ defined on
$[a,b]$. Let ${\cal L}_N$ denote one of the two maximal subpaths
of $\BBz_N^2$ covering $f$.
 Let $N_+({\cal L}_N)$ be as
defined in (\ref{def}). Then
$$
\lim_{N\rightarrow +\infty}\frac{ N_+({\cal
L}_N)}{N}=\int_a^b\left(|f'(x)|\wedge 1\right)dx.
$$
\end{lemma}
{\bf{Proof of lemma \ref{ana}.}} We suppose without loss of
generality that  the function $f$ is nondecreasing on $[a,b]$. Let
$(I_i)_{i\in I}$ be the collection of the open intervals where
$f'-1$ is nonzero. Setting $I_i=]x_{i-1},x_i[$ for $i\in I$, we
have
\begin{equation}\label{l7heu1}
\left(\frac{f(x_i)-f(x_{i-1})}{x_i-x_{i-1}}\right)\wedge 1=
\frac{1}{x_i-x_{i-1}}\int_{x_{i-1}}^{x_i}(f'(x)\wedge 1)\,dx.
\end{equation}
We denote by $f_i$ the restriction of $f$ to $[x_{i-1},x_i[$ and
by ${\cal L}_N^{(i)}$ the associated polygonal line. We deduce
from the suitable construction of the intervals $(I_i)_{i\in I}$
and arguing as in the proof of lemma \ref{lem5bis}, that
$$
\lim_{N\rightarrow +\infty}\frac{ N_+({\cal
L}^{(i)}_N)}{N}=(x_i-x_{i-1})\wedge (f(x_i)-f(x_{i-1})).
$$
Hence
$$
\lim_{N\rightarrow +\infty}\frac{ N_+({\cal L}_N)}{N}=\sum_{i\in
I} (x_i-x_{i-1})\wedge (f(x_i)-f(x_{i-1})).
$$
Lemma \ref{ana} is proved by collecting the last bound together
with (\ref{l7heu1}). \ \ \ $\Box$
\medskip

\noindent We define a monotone function $f$, such that the part of
$\gamma$ limited by $x_0$ and $x_0(\delta)$ (where $x_0(\delta)$
is the point of $\gamma\cap \partial B(x_0,\delta)\setminus
B(s,r)$) is equal to the graph $\{(x,y):\, y=f(x)\}$ and we apply
lemma \ref{ana} to the monotone path ${\cal L}_N$ covering the
part of $\gamma$ limited by $x_0$ and $x_0(\delta)$. We deduce,
since $\left|N_+({\cal L}_N)-C_N(x_0,\delta)\right|\leq 1,$ that
\begin{eqnarray*}
\lim_{N\rightarrow +\infty}\frac{C_N(x_0,\delta)}{N\delta}& = &
\frac{1}{\delta}\int_{I_{\delta}}\left(|f'(x)|\wedge
1\right)dx \\
&=&
|\cos\theta(\delta)|\frac{1}{I_{\delta}}\int_{I_{\delta}}\left(|f'(x)|\wedge
1\right)dx,
\end{eqnarray*}
where $I_{\delta}$ is the segment $[x_0\cdot i, x_0(\delta)\cdot
i]$. We obtain, taking the limit over $\delta\rightarrow 0$ in the
last equality,
\begin{eqnarray*}
\lim_{\delta\rightarrow 0}\lim_{N\rightarrow
+\infty}\frac{C_N(x_0,\delta)}{N\delta}\, =\,  |\cos\theta_0|
\left(|f'(x_0\cdot i)|\wedge 1\right) &=& |\cos\theta_0|\wedge
|\sin\theta_0|.
\end{eqnarray*}
 We then obtain from the conclusion of theorem
\ref{T2}, that for $r$ small enough,
  \begin{eqnarray}\label{il}
{\lefteqn{\lim_{\delta\rightarrow 0}\liminf_{N\rightarrow
+\infty}{\bf{A}}^{\sigma, \gamma}_N(s,r,\delta)=
\lim_{\delta\rightarrow 0}
  \limsup_{N\rightarrow +\infty}{\bf{A}}^{\sigma, \gamma}_N(s,r,\delta) }}
\\
&&
  =\frac{1}{4}\sin2\theta_{0}\left(2\BBone_{|\sin\theta_{0}|<|\cos\theta_{0}|}-1\right)
-\frac{1}{4}\sin2\theta_{1}\left(2\BBone_{|\sin\theta_{1}|<|\cos\theta_{1}|}-1\right){\nonumber}\\
&&
+\frac{1}{2}\left(\sg(\tan\theta_1)\BBone_{|\sin\theta_{1}|<|\cos\theta_{1}|}-
\sg(\tan\theta_{0})\BBone_{|\sin\theta_{0}|<|\cos\theta_{0}|}\right)\\
&&+
\BBone_{\sin\theta_{0}\sin\theta_{1}>0}\left(\sg(\theta_{1}-\theta_{0})\BBone_{cos\theta_{0}\cos\theta_{1}>0}+
\sg(\tan\theta_0)\BBone_{\cos\theta_{0}\cos\theta_{1}<0}\right).{\nonumber}
\end{eqnarray}

\noindent{\bf{End of the proof of proposition \ref{theo1heu1}.}}
In order to prove proposition \ref{theo1heu1}, we suppose first
that $\theta$ takes a value different from $(2k+1)\frac{\textstyle
\pi}{\textstyle 4}$, for $k\in \BBn$. Since the curve $\gamma$
admits a tangent at the point $s$, then for $r$ small enough,
$(\theta_{0},\theta_1)$ belongs to $[(2k+1)\frac{\textstyle
\pi}{\textstyle 4},(2k+3)\frac{\textstyle \pi}{\textstyle
4}]\times [(2k+5)\frac{\textstyle \pi}{\textstyle
4},(2k+7)\frac{\textstyle \pi}{\textstyle 4}],$ for some $k\in
\BBn$. We then deduce from (\ref{il}) that,
  \begin{itemize}
  \item if\ \
$(\theta_{0},\theta_1)\in[(2k+1)\frac{\textstyle \pi}{\textstyle
4},(2k+3)\frac{\textstyle \pi}{\textstyle 4}]\times
[(2k+5)\frac{\textstyle \pi}{\textstyle 4},(2k+7)\frac{\textstyle
\pi}{\textstyle 4}],$ with {\bf{$k\in\{0,2\}$}}, then $$
\lim_{\delta\rightarrow 0} \lim_{\varepsilon\rightarrow
0}\limsup_{N\rightarrow+\infty}{\bf{A}}^{\sigma,
\gamma}_N(s,r,\delta) =
 \frac{1}{4}\left({\sin
2\theta_{1}}-{\sin 2\theta_{0}}\right). $$
\item if\ \
$(\theta_{0},\theta_1)\in[(2k+1)\frac{\textstyle \pi}{\textstyle
4},(2k+3)\frac{\textstyle \pi}{\textstyle 4}]\times
[(2k+5)\frac{\textstyle \pi}{\textstyle 4},(2k+7)\frac{\textstyle
\pi}{\textstyle 4}],$ with {\bf{$k\in\{1,3\}$}}, then $$
\lim_{\delta\rightarrow 0}
\limsup_{N\rightarrow+\infty}{\bf{A}}^{\sigma,
\gamma}_N(s,r,\delta) =
 \frac{1}{4}\left({\sin 2\theta_{0}}-{\sin 2\theta_1}\right). $$
 \end{itemize}
 We now need the following lemma.
 \begin{lemma}\label{ana1}
 Let $\gamma$ be a Jordan curve of $\BBr^2$ of class ${\cal C}_2$. Let $s$ be a fixed point of
  $\gamma$. Let $r$ be a positive real number sufficiently small
  such that $\partial B(s,r)\cap \gamma$ contains exactly two
  points
  $x_0$ and $x_1$. Suppose that $x_0$, $s$ and $x_1$ are arranged counterclockwise.
  Let $s'$ be the common point to $T_{x_0}\gamma$ and $T_{x_1}\gamma$. Let $\theta_1\in [0,2\pi]$
(respectively $\theta_{0}\in [0,2\pi]$)
 be the oriented angle between the half horizontal axis $[0,+\infty[$
 and the segment
 $[s',x_{1}[$ (respectively $[s',x_0[$).
  Then
  $$
  \lim_{r\rightarrow 0} \frac{\sin
  \left(\theta_0-
  \theta_1\right)}{2r}=\xi_{\gamma}(s),
  $$
  and
  $$
   \lim_{r\rightarrow 0} {\cos
  \left(\theta_0+
  \theta_1\right)}=-\cos 2\theta,
  $$
  where $\theta$ is the angle between the half horizontal axis $[0,+\infty[$ and
  $T_s\gamma$.
  \end{lemma}
Lemma \ref{ana1}, together with the two equalities just above
Lemma \ref{ana1} and the fact ${\sin 2a}-{\sin
2b}=2\sin(a-b)\cos(a+b)$, gives
 $$
\lim_{r\rightarrow 0}\lim_{\delta \rightarrow 0
}\limsup_{N\rightarrow+\infty}\frac{1}{2r}{\bf{A}}^{\sigma,
\gamma}_N(s,r,\delta) =\left\{
\begin{array}{rl}
\frac{\textstyle 1}{\textstyle 2}(\cos 2\theta)\,\xi_{\gamma}(s) &
if \ \ \theta\in ](1+4k)\frac{\textstyle \pi}{\textstyle
4},(3+4k)\frac{\textstyle \pi}{\textstyle 4}[
\\
\\
-\frac{\textstyle 1}{\textstyle 2}(\cos 2\theta)\,\xi_{\gamma}(s)
 & if \ \theta\in ](3+4k)\frac{\textstyle
\pi}{\textstyle 4}, (5+4k)\frac{\textstyle \pi}{\textstyle 4}[.
\end{array}
\right.
  $$
  which proves theorem \ref{theo1heu1} when $\theta$ is different
  from $(2k+1)\frac{\textstyle
\pi}{\textstyle 4}$, for $k\in \BBn$.  Now, suppose that
$\theta=\frac{\textstyle \pi}{\textstyle 4}$ and that for any $r$
small enough $(\theta_{0},\theta_1)\in[\frac{\textstyle
\pi}{\textstyle 4},3\frac{\textstyle\pi}{\textstyle 4}]\times
[3\frac{\textstyle \pi}{\textstyle 4},5\frac{\textstyle
\pi}{\textstyle 4}]$ (the arguments for the proof for the other
values of $\theta$ and the corresponding values of $\theta_1$,
$\theta_0$ will be similar). We have in that case,
\begin{eqnarray}\label{l61heu1}
\lim_{\delta\rightarrow 0}
\limsup_{N\rightarrow+\infty}{\bf{A}}^{\sigma,
\gamma}_N(s,r,\delta) &=& \frac{\textstyle 1}{\textstyle
4}\left(2-{\sin 2\theta_{1}}-{\sin
2\theta_{0}}\right)\nonumber\\
&=& \frac{\textstyle 1}{\textstyle 2}\left(\sin(\frac{\textstyle
\pi}{\textstyle 4}-\theta_1)\cos(\frac{\textstyle \pi}{\textstyle
4}+\theta_1)+ \sin(\frac{\textstyle \pi}{\textstyle
4}-\theta_0)\cos(\frac{\textstyle \pi}{\textstyle
4}+\theta_0)\right).
\end{eqnarray}
Now the method of the proof of lemma \ref{ana1} gives
$$
\lim_{r\rightarrow 0} \frac{\sin
  \left(\theta-
  \theta_1\right)}{r}=\lim_{r\rightarrow 0} \frac{\sin
  \left(\theta-
  \theta_0\right)}{r}=-\xi_{\gamma}(s).
$$
This fact, together with (\ref{l61heu1}), leads to
$$
\lim_{r\rightarrow 0}\lim_{\delta\rightarrow 0}
\limsup_{N\rightarrow+\infty}\frac{1}{r}{\bf{A}}^{\sigma,
\gamma}_N(s,r,\delta) =0,
$$
which is the conclusion of theorem \ref{theo1heu1} for
$\theta=\frac{\textstyle \pi}{\textstyle 4}$.
\medskip

\noindent
{\bf{Proof of lemma \ref{ana1}.}} We begin by giving the
definition of the curvature of $\gamma$ at any $s\in \gamma$.
\medskip

\noindent
{\it Definition.} Let $\gamma$ be a smooth Jordan curve of
$\BBr^2$. Suppose that $(\phi(t))_{t\in [-1,1]}$ is a
parametrization of the
curve $\gamma$. Let $s=\phi(t)=(x(t),y(t))$ be a fixed point of $\gamma$. The
{\bf{curvature}} of $\gamma$ at the point $s$ is
defined by
$$
\xi_{\gamma}(s)=\frac{x'(t)y''(t)-x''(t)y'(t)}{(x^{'2}(t)+
y^{'2}(t))^{3/2}}.
$$
Let $s$, $x_0$ and $x_1$ be as defined in lemma \ref{ana1}. Let
$t$, $t_0$ and $t_1$ be three real numbers of $[-1,1]$ such that
$s= \phi(t)=(x(t), y(t)),$ and for $i\in\{0,1\},$
$x_i=\phi(t_i)=(x(t_i), y(t_i)).$ We have
$r^2=\left(x(t_i)-x(t)\right)^2+ \left(y(t_i)-y(t)\right)^2$, for
$i\in\{0,1\}$. Hence
$$\lim_{t_0\rightarrow t,\, t_0<t }\frac{r}{t-t_0}  =
\sqrt{x^{'2}(t)+ y^{'2}(t)},\ \ \ \ \ \, \lim_{t_1\rightarrow t,\,
t<t_1 }\frac{r}{t_1-t} = \sqrt{x^{'2}(t)+ y^{'2}(t)}. $$ For any
$\tau\in[-1,1]$, define $f(\tau)=\frac{\textstyle
x'(\tau)}{\textstyle \sqrt{x^{'2}(\tau)+y^{'2}(\tau)}}.$ We have
$$
f'(\tau)=\frac{x''(\tau)}{\sqrt{x^{'2}(\tau)+y^{'2}(\tau)}}-x'(\tau)\frac{x'(\tau)x''(\tau)+y'(\tau)y''(\tau)}{(x^{'2}(\tau)+y^{'2}(\tau))^{3/2}}.
$$
Hence
\begin{eqnarray}\label{lim14}
{\lefteqn{\cos\theta_0  =-\frac{x'(t_0)}{\sqrt{x^{'2}(t_0)+
y^{'2}(t_0)}}}} \nonumber\\ && = -\frac{x'(t)}{\sqrt{x^{'2}(t)+
y^{'2}(t)}}+(t-t_0)f'(t) + o\left(|t-t_0|\right).
\end{eqnarray}
\begin{equation}\label{lim14a}
\cos\theta_1  = \frac{x'(t)}{\sqrt{x^{'2}(t)+
y^{'2}(t)}}+(t_1-t)f'(t) + o\left(|t_1-t|\right).
\end{equation}
We obtain, combining the last two equalities
\begin{eqnarray*}
\lim_{t_1\rightarrow t,\, t_0\rightarrow t,\, t_0<t<t_1}
\frac{\cos\theta_0+ \cos\theta_1}{r} = \frac{2x''(t)}{{x^{'2}(t)+
y^{'2}(t)}}-2x'(t)\frac{x'(t)x''(t)+y'(t)y''(t)}{(x^{'2}(t)+y^{'2}(t))^{2}}.
\end{eqnarray*}
The last limit together with  $$ \lim_{t_1\rightarrow t,\, t<t_1
}{\sin\theta_1}=\frac{y'(t)}{\sqrt{x^{'2}(t)+ y^{'2}(t)}},$$
ensures
$$\displaylines{
\lim_{t_1\rightarrow t,\, t_0\rightarrow t,\, t_0<t<t_1}
\frac{1}{r}\sin\theta_1\left(\cos\theta_0+ \cos\theta_1\right)\, =
\hfill\cr
\frac{2x''(t)y'(t)}{(x^{'2}(t)+
y^{'2}(t))^{3/2}}-2x'(t)y'(t)\frac{x'(t)x''(t)+y'(t)y''(t)}{(x^{'2}(t)+y^{'2}(t))^{5/2}}\,.
}$$
In the same way, we prove that
$$\displaylines{
\lim_{t_1\rightarrow t,\, t_0\rightarrow t,\, t_0<t<t_1}
\frac{1}{r}\cos\theta_1\left(\sin\theta_0+ \sin\theta_1\right)\, =\hfill\cr
\frac{2x'(t)y''(t)}{(x^{'2}(t)+
y^{'2}(t))^{3/2}}-2x'(t)y'(t)\frac{x'(t)x''(t)+y'(t)y''(t)}{(x^{'2}(t)+y^{'2}(t))^{5/2}}\,.
}$$
The last two limits together with
\begin{eqnarray*}
\sin
  \left(\theta_0-
  \theta_1\right) & = & \cos\theta_1\left(\sin\theta_1+\sin\theta_0\right)
  - \sin\theta_1\left(\cos\theta_0+\cos\theta_1\right),
\end{eqnarray*}
prove that $$\lim_{t_1\rightarrow t,\, t_0\rightarrow t,\,
t_0<t<t_1} \frac{1}{2r}\sin
  \left(\theta_0-
  \theta_1\right)=\frac{x'(t)y''(t)-x''(t)y'(t)}{(x^{'2}(t)+ y^{'2}(t))^{3/2}}.$$
 Now the equality
 $$
\cos
  \left(\theta_0+
  \theta_1\right)= \cos\theta_0\cos\theta_1-
  \sin\theta_0\sin\theta_1,
 $$
 together with the limits (\ref{lim14}), (\ref{lim14a}), yields
 $$
  \lim_{t_1\rightarrow t,\, t_0\rightarrow t,\,
t_0<t<t_1} {\cos
  \left(\theta_0+
  \theta_1\right)}=\frac{y^{'2}(t)-x^{'2}(t)}{{x^{'2}(t)+
y^{'2}(t)}}.
  $$
  The last limit is equal to $-\cos2\theta$, where $\theta$ is the angle
between the horizontal axis and $T_s\gamma$. \\
   $\Box$
\subsection{Proof of proposition \ref{P1}}
Our purpose is to apply theorem \ref{T2}. So we have  to check,
the requirements of theorem \ref{T2}. We first prove the condition
(\ref{haus}). We claim that, for all $\varepsilon>0$,
\begin{equation}\label{spr}
\lim_{N\rightarrow \infty}
{\mu_N}\kern-3pt\left(\sup|\Phi_N(x_N)-f(x_N)|\geq \varepsilon
\right)=0,
\end{equation}
 where the supremum is taken over $x_N \in [a,b]\cap
 \frac{\BBz}{N}$.
\\
{\it Proof of (\ref{spr}).} For $l\in\BBz$, we denote by
$\eta(\frac{l}{N})$ the height difference
$\eta(\frac{l}{N})=\Phi_N(\frac{l+1}{N})-\Phi_N(\frac{l}{N})$.
Without loss of generality, we will take $a=0$. We write, for
$\frac{k}{N}\,\in\, [0,b]\cap \frac{\BBz}{N}$,
$$
\Phi_N(\frac{k}{N})-f(\frac{k}{N})=
\sum_{l=0}^{k-1}\left(\eta(\frac{l}{N})-(f(\frac{l+1}{N})-f(\frac{l}{N}))\right)
+ \left(\Phi_N(0)-f(0)\right)\,.
$$
The last equality gives, since
${\mu_N}\kern-3pt\left(\eta(\frac{k}{N})\right)=
\frac{1}{N}|f'|(\frac{k}{N}),$
\begin{eqnarray*}
{\lefteqn{\Phi_N(\frac{k}{N})-f(\frac{k}{N})\,=}}\\
&&
\sum_{l=0}^{k-1}\left(\eta(\frac{l}{N})-\mu_N\left(\eta(\frac{l}{N})\right)\right)
-\sum_{l=0}^{k-1}\left((f(\frac{l+1}{N})-f(\frac{l}{N}))-\frac{1}{N}|f'|(\frac{l}{N})\right)
+ \left(\Phi_N(0)-f(0)\right).
\end{eqnarray*}
We deduce from the last equality, assumption (\ref{(i)}) of
proposition \ref{P1} and the fact
$$
\sum_{l=0}^{k-1}\left|f(\frac{l+1}{N})-f(\frac{l}{N})-\frac{1}{N}|f'|(\frac{l}{N})\right|\leq
\frac{b}{N}\|f''\|_{\infty},
$$
that (\ref{spr}) is proved as soon as,
\begin{equation}\label{sprbis}
\lim_{N\rightarrow \infty}{\mu_N}\kern-4pt\left(\sup_{0\leq
k\leq\,
Nb}\left|\sum_{l=0}^{k-1}\left(\eta(\frac{l}{N})-\mu_N\left(\eta(\frac{l}{N})\right)\right)\right|\geq
\varepsilon \right)=0.
\end{equation}
For this, we use a Markov inequality, the independence of the
random variables $(\eta(\frac{l}{N}))_{l\in \BBz}$ and a Rosenthal
inequality (cf. section 2.6.19 and Theorem 2.9 of Petrov (1995)).
We get, for an universal constant $C$,
\begin{eqnarray*}
{\lefteqn{{\mu_N}\kern-3pt\left(\sup_{0\leq k\leq\,
Nb}\left|\sum_{l=0}^{k-1}\left(\eta(\frac{l}{N})-\mu_N\left(\eta(\frac{l}{N})\right)\right)\right|\geq
\varepsilon \right)}} \\
&&\leq \frac{1}{\varepsilon^3}\mu_N\kern-4pt\left(\sup_{0\leq k\leq\,
Nb}\left|\sum_{l=0}^{k-1}\left(\eta(\frac{l}{N})-\mu_N\left(\eta(\frac{l}{N})\right)\right)\right|\right)^3\\
&& \leq
\frac{C}{\varepsilon^3}\left\{\left(\sum_{l=0}^{Nb}\Var_{\mu_N}\left(\eta(\frac{l}{N})\right)\right)^{3/2}+
\sum_{l=0}^{Nb}\mu_N\left(\left|\eta(\frac{l}{N})-\mu_N\left(\eta(\frac{l}{N})\right)\right|^3\right)\right\}\,.
\end{eqnarray*}
The last estimations and the fact that, for some constant $C$
depending on $\|f'\|_{\infty}$,
$$\Var_{\mu_N}\eta(\frac{l}{N})=\frac{1}{N^2}|f'|(\frac{l}{N})\left(1+|f'|(\frac{l}{N})\right)^2,\ \ \
\mu_N\left(\left|\eta(\frac{l}{N})\right|^3\right)\leq C \frac{1}{N^3} $$
give
\begin{eqnarray*}
{\mu_N}\kern-4pt\left(\sup_{0\leq k\leq\,
Nb}\left|\sum_{l=0}^{k-1}\left(\eta(\frac{l}{N})-\mu_N\left(\eta(\frac{l}{N})\right)\right)\right|\geq
\varepsilon \right)={\cal{O}}\left(
\left(\frac{1}{N}\right)^{3/2}\right),
\end{eqnarray*}
which proves (\ref{sprbis}) and then (\ref{spr}).  Now (\ref{spr})
allows to deduce the condition (\ref{haus}).
\\
We deduce from the definition of $\mu_N$, that for any $N\in
\BBn^*$
$$
\forall{k}\in [Na,Nb]\cap{\BBz}\qquad
{\mu_N}\left(\sg(\eta(\frac{k}{N})\,f'(\frac{k}{N}))<0\right)=0.
$$
Since the graph of the monotone function~$f$ coincides
with the
restriction of $\gamma$ over $[a,b]$, we conclude from the above
formula
that $\partial A^N_{\sigma}\cap {\cal S}(s,r,\delta,\delta)$ and
$\gamma \cap {\cal S}(s,r,\delta,\delta)$ are both nondecreasing
or both nonincreasing.
\\
Our task now is to check the condition~(\ref{c5bis}) and to
precise the value of the corresponding function $C$. Recall that
$f$ and $\Phi_N$ are both increasing or decreasing.  Therefore
$$
C_N(x_0,\delta)=\sum_{x_0\cdot i\leq k/N\leq\,
\delta|\cos\theta_0|+x_0\cdot i }\kern-10pt
\BBone_{\eta(\frac{k}{N})\,\neq\, 0},
$$
where the quantity $C_N(x_0,\delta)$ is defined just before
theorem \ref{T1}. We have
\begin{eqnarray*}
\frac{1}{N\delta}\mu_N\kern-4pt\left(C_N(x_0,\delta)\right)& = &
\frac{1}{N\delta}\sum_{x_0\cdot i\,\leq k/N\leq\,
\delta|\cos\theta_0|+ x_0\cdot i}\kern-10pt
\mu_N({\eta(\frac{k}{N})\,\neq\,
0})\\
&=& \frac{1}{N\delta}\sum_{x_0\cdot i\leq k/N\leq\,
\delta|\cos\theta_0|+x_0\cdot
i}\frac{|f'|(\frac{k}{N})}{1+|f'|(\frac{k}{N})}.
\end{eqnarray*}
The last equality gives
\begin{eqnarray*}
\lim_{N\rightarrow\infty}\frac{1}{N\delta}\mu_N\kern-4pt\left(C_N(x_0,\delta)\right)&=&\frac{1}{\delta}\int_{x_0\cdot
i}^{\delta|\cos\theta_0|+x_0\cdot i}
\frac{|f'|(x)}{1+|f'|(x)}d\,x.
\end{eqnarray*}
Hence
\begin{eqnarray*}
\lim_{\delta\rightarrow 0}
\lim_{N\rightarrow\infty}\frac{1}{N\delta}\mu_N\kern-4pt\left(C_N(x_0,\delta)\right)&
= &|\cos\theta_0| \frac{|f'|(x_0\cdot i)}{1+|f'|(x_0\cdot i)}
= |\cos\theta_0| \frac{|\tan \theta_0|}{1+|\tan\theta_0|}\\
&=& \frac{|\sin (2\theta_0)|}{2(|\sin\theta_0|+
|\cos\theta_0|)}=C(\theta_0).
\end{eqnarray*}
We have assumed in proposition \ref{P1} that the curve $\gamma$ is
monotone in $B(s,r)\cup B(x_0,\delta)\cup B(x_1,\delta)$. This
fact allows to deduce that,
$$
\sin\theta_0\sin\theta_1\leq 0,\ \ \cos\theta_0\cos\theta_1\leq 0.
$$
We use the last constatation together with the conclusion of
theorem \ref{T2} to obtain,
\begin{eqnarray}\label{fins}
{\lefteqn{ \lim_{\delta\rightarrow
0}\limsup_{N\rightarrow\infty}\,\mu_N\kern-4pt\left({\bf{A}}^{{\sigma,\gamma}}_N(s,r,\delta)\right)=
\lim_{\delta\rightarrow 0}\liminf_{N\rightarrow\infty}\mu_N\kern-4pt\left({\bf{A}}^{{\sigma,\gamma}}_N(s,r,\delta)\right)}}\\
&&
=\frac{1}{2}\,\sg(\tan\theta_0)\left(\cos^2\theta_1-\cos^2\theta_0\right){\nonumber}\\
&&  + \frac{1}{2}\,\sg(\tan\theta_0)\bigg(\frac{|\sin
(2\theta_1)|}{2(|\sin\theta_1|+
|\cos\theta_1|)}(\left|\sin\theta_1\right|
-
\left|\cos\theta_1\right|){\nonumber}\\
&&\kern75pt- \frac{|\sin
(2\theta_0)|}{2(|\sin\theta_0|+
|\cos\theta_0|)}\left(\left|\sin\theta_0\right|-
\left|\cos\theta_0\right|\right)\bigg).{\nonumber}
  \end{eqnarray}
  We have
  $$
 \sg(\tan\theta_0)\left(\frac{|\sin
(2\theta_1)|}{2(|\sin\theta_1|+
|\cos\theta_1|)}(\left|\sin\theta_1\right|-
\left|\cos\theta_1\right|)- \frac{|\sin
(2\theta_0)|}{2(|\sin\theta_0|+
|\cos\theta_0|)}\left(\left|\sin\theta_0\right|-
\left|\cos\theta_0\right|\right)\right)
  $$
  $$
=\frac{\sin \theta_1\cos\theta_1}{(|\sin\theta_1|+
|\cos\theta_1|)}(\left|\sin\theta_1\right|-
\left|\cos\theta_1\right|)-
\frac{\sin\theta_0\cos\theta_0}{(|\sin\theta_0|+
|\cos\theta_0|)}\left(\left|\sin\theta_0\right|-
\left|\cos\theta_0\right|\right)
  $$
  $$
  =\frac{-\sin(\theta_0-\theta_1)}{(\left|\sin\theta_1\right|+
\left|\cos\theta_1\right|)(|\sin\theta_0|+ |\cos\theta_0|)}+$$
$$\frac{
(\cos^2\theta_1-\cos^2\theta_0)(\sin\theta_0\cos\theta_1+
\sin\theta_1\cos\theta_0+
\cos\theta_0\cos\theta_1\sg(\tan\theta_0)+
\sin\theta_0\sin\theta_1\sg(\tan\theta_0))}{(\left|\sin\theta_1\right|+
\left|\cos\theta_1\right|)(|\sin\theta_0|+ |\cos\theta_0|)}
  $$
  $$
  =\frac{-\sin(\theta_0-\theta_1)}{(\left|\sin\theta_1\right|+
\left|\cos\theta_1\right|)(|\sin\theta_0|+
|\cos\theta_0|)}-\sg(\tan\theta_0)(\cos^2\theta_1-\cos^2\theta_0).
  $$
  We conclude from (\ref{fins}) together with the last equalities,
\begin{eqnarray*}
{\lefteqn{ \lim_{\delta\rightarrow
0}\limsup_{N\rightarrow\infty}\mu_N\kern-4pt\left({\bf{A}}^{{\sigma,\gamma}}_N(s,r,\delta)\right)=
\lim_{\delta\rightarrow 0}\liminf_{N\rightarrow\infty}\mu_N\kern-4pt\left({\bf{A}}^{{\sigma,\gamma}}_N(s,r,\delta)\right)}}\\
&& = \frac{-\sin(\theta_0-\theta_1)}{2(\left|\sin\theta_1\right|+
\left|\cos\theta_1\right|)(|\sin\theta_0|+ |\cos\theta_0|)}.
\end{eqnarray*}
The last limit together with lemma \ref{ana1} completes the proof
of proposition \ref{P1}.\,\, $\Box$.


\begin{thebibliography}{9}

\bibitem{EDA} Andjel, E. D.: Invariant measures for the zero range
process. Ann. Probab. 10 (1982), no. 3, 525--547.

\bibitem{CSS}
Chayes, L., Schonmann, R. H., Swindle, G.: Lifshitz' law for the
volume of a two-dimensional droplet at zero temperature. J.
Statist. Phys. 79 (1995), no. 5-6, 821--831.

\bibitem{CS}
Chayes, L., Swindle, G.: Hydrodynamic limits for one-dimensional
particle systems with moving boundaries. Ann. Probab. 24 (1996),
no. 2, 559--598.

\bibitem{DOPT1}
De Masi, A., Orlandi, E., Presutti, E., Triolo, L.: Glauber
evolution with the Kac potentials. I. Mesoscopic and macroscopic
limits, interface dynamics. Nonlinearity 7 (1994), no. 3,
633--696.

\bibitem{DOPT2}
De Masi, A., Orlandi, E., Presutti, E., Triolo, L.: Motion by
curvature by scaling nonlocal evolution equations. J. Statist.
Phys. 73 (1993), no. 3-4, 543--570.

\bibitem{KS}
Katsoulakis, M. A., Souganidis, P. E.: Stochastic Ising models and
anisotropic front propagation. J. Statist. Phys. 87 (1997), no.
1-2, 63--89.

\bibitem{KSP}
Katsoulakis, Markos A., Souganidis, Panagiotis E: Generalized
motion by mean curvature as a macroscopic limit of stochastic
Ising models with long range interactions and Glauber dynamics.
Comm. Math. Phys. 169 (1995), no. 1, 61--97.

\item Petrov, V. V. (1995) {Limit theorems of probability theory:
sequences of independent random variables.} Clarendon Press,
Oxford.

\bibitem{So}
Sowers, Richard B.: Hydrodynamical limits and geometric measure
theory: mean curvature limits from a threshold voter model. J.
Funct. Anal. 169 (1999), no. 2, 421--455.

\bibitem{Sp}
Spohn, Herbert: Interface motion in models with stochastic
dynamics. J. Statist. Phys. 71 (1993), no. 5-6, 1081--1132.

\end{thebibliography}
\end{document}